\documentclass[10pt]{article}

\UseRawInputEncoding

\usepackage{fullpage,graphicx,subfigure,mathdots,mathpazo}
\usepackage{amsmath,amscd,tikz,mathrsfs}
\usepackage[normalem]{ulem}
\usepackage{amsmath}
\usepackage{setspace}

\usepackage[latin1]{inputenc}
\usepackage{epsfig,amsmath,amsthm,graphicx,amssymb,overpic}
\usepackage{epsfig,amsmath,graphicx,amssymb,overpic,cite}
\usepackage{listings,diagbox}
\usepackage{epsfig,amsmath,graphicx,amssymb,overpic}
\usepackage{bm}
\usepackage{graphicx}
\usepackage{subfigure}

\usepackage{pdfsync}
\usepackage{cite}
\usepackage{amscd}
\usepackage{amsmath}
\usepackage{latexsym}
\usepackage{amsfonts}
\usepackage{amssymb}
\usepackage{amsthm}
\usepackage{graphicx}
\usepackage{indentfirst}
\usepackage{enumerate}
\usepackage{amstext}
\usepackage[dvipdfm,
            pdfstartview=FitH,
            CJKbookmarks=true,
            bookmarksnumbered=true,
            bookmarksopen=true,
            pdfborder=001,   
            linkcolor=blue,
            ]{hyperref}

 \allowdisplaybreaks[4]

\def\be{\begin{equation}}
\def\ee{\end{equation}}
\def\bee{\begin{eqnarray}}
\def\ene{\end{eqnarray}}
\def\bes{\begin{subequations}}
\def\ees{\end{subequations}}
\def\no{\nonumber}
\def\v{\vspace{0.1in}}

\newtheorem{thm}{\textbf Theorem}[section]
\newtheorem{lem}{\textbf Lemma}[section]
\newtheorem{rem}{\textbf Remark}[section]
\newtheorem{cor}{\textbf Corollary}[section]

\newtheorem{prop}{\textbf Proposition}[section]
\newtheorem{definition}{\textbf Definition}[section]

\numberwithin{equation}{section}

\renewcommand{\theequation}{\arabic{section}.\arabic{equation}}

\begin{document}

\baselineskip=13pt
\renewcommand {\thefootnote}{\dag}
\renewcommand {\thefootnote}{\ddag}
\renewcommand {\thefootnote}{ }

\pagestyle{plain}

\begin{center}
\baselineskip=16pt \leftline{} \vspace{-.3in} {\Large \bf The Cauchy problem and multi-peakons for the mCH-Novikov-CH equation with quadratic and cubic nonlinearities} \\[0.2in]
\end{center}

\begin{center}
{\bf Guoquan Qin$^{\rm a,b}$, Zhenya Yan$^{\rm b,c,*}$, Boling Guo$^{\rm d}$}\footnote{$^{*}$Corresponding author. {\it Email address}: zyyan@mmrc.iss.ac.cn (Z. Yan). }  \\[0.1in]
{\it\small  $^a$National Center for Mathematics and Interdisciplinary Sciences, Academy of Mathematics and Systems Science, Chinese Academy of Sciences, Beijing 100190, China \\
$^b$Key Laboratory of Mathematics Mechanization, Academy of Mathematics and Systems Science, \\ Chinese Academy of Sciences, Beijing 100190, China \\
 $^c$School of Mathematical Sciences, University of Chinese Academy of Sciences, Beijing 100049, China\\
 $^d$Institute of Applied Physics and Computational Mathematics, Beijing 100088, PR China} \\
\end{center}

\vspace{0.4in}

{\baselineskip=13pt


\vspace{-0.28in}

\noindent {\bf Abstract}\, {\small This paper investigates the Cauchy problem of a generalized Camassa-Holm equation with quadratic and cubic nonlinearities (alias the mCH-Novikov-CH equation), which is
a generalization of  some special equations such as the Camassa-Holm (CH) equation, the modified CH (mCH) equation ((alias the Fokas-Olver-Rosenau-Qiao equation), the Novikov equation, the CH-mCH equation, the mCH-Novikov equation, and the CH-Novikov equation. We first show the local well-posedness for the strong solutions of the mCH-Novikov-CH equation in Besov spaces by means of the Littlewood-Paley theory and the transport equations theory. Then, the H\"{o}lder continuity of the data-to-solution map to this equation are exhibited in some Sobolev spaces. After providing the  blow-up criterion and  the precise blow-up quantity in light of the  Moser-type estimate in the Sobolev spaces, we then trace a portion and the whole of the precise blow-up quantity, respectively, along the characteristics associated with this equation, and
obtain two kinds of sufficient conditions on the gradient of the initial data to guarantee the occurance of the
wave-breaking phenomenon. Finally, the non-periodic and periodic peakon and multi-peakon solutions for this equation are also explored. }



\vspace{0.1in} \noindent {\bf Keywords}\, MCH-Novikov-CH equation $\cdot$ Wave breaking $\cdot$ Local well-posedness $\cdot$
H\"older continuity $\cdot$ Non-periodic and periodic peakon and multi-peakon solutions}

\vspace{0.1in} \noindent {\bf Mathematics Subject Classification}\, 35B30; 35G25; 35A01; 35B44; 35Q53; 35Q35

\baselineskip=13pt

\vspace{0.2in}

\section{Introduction}

We are concerned with the Cauchy problem and multi-peakons of the following nonlinear dispersive equation with  quadratic and cubic nonlinearities, which is a combination of the modified Camassa-Holm equation, the Novikov equation,  and
the Camassa-Holm equation (alias the mCH-Novikov-CH equation)
\begin{eqnarray}\begin{cases}\label{gmCH-Novikov}
m_{t}+k_{1}[(u^{2}-u_{x}^{2})m]_{x} +k_{2}(u^{2} m_{x}+3 u u_{x} m)+k_{3}(u m_{x}+2u_{x} m)=0, \quad x\in\mathbb{R},\,\, t>0, \\
u(0, x)=u_{0}(x),\quad x\in\mathbb{R},
\end{cases}
\end{eqnarray}
where $k_{1}, k_{2}$ and $k_{3}$ are all real-valued parameters, and $m=u-u_{xx}$.

Eq.~(\ref{gmCH-Novikov}) with $k_{1,2}=0$ and $k_{3}=1$ reduces to the celebrated Camassa-Holm (CH) equation~\cite{CamassaHolm1993PRL}
\begin{eqnarray}\label{CH}
m_{t}
+u m_{x}+2u_{x} m=0, \quad m=u-u_{xx},
\end{eqnarray}
where $u(t, x)$ denotes the fluid velocity and $m(t, x)$ stands for the corresponding potential density,
which can describe the unidirectional propagation of shallow water waves over a flat bottom \cite{CamassaHolm1993PRL,
DullinGottwaldHolm2001PRL,Johnson2002JFM} and the propagation of axially symmetric waves in hyperelastic rods
\cite{ConstantinStrauss2000PLA, Dai1998AM}. Equation (\ref{CH}) can  be deduced  from the well-known Korteweg-de Vries (KdV) equation by means of the tri-Hamiltonian duality~\cite{OlverRosenau1996PRE} and shares some alike properties with the KdV equation.
For instance, it is  completely integrable and has infinitely many conservation laws
\cite{CamassaHolm1993PRL,
FisherSchiff1999PLA,
FokasFuchssteiner1981PD}.
It can be written as the
 following  bi-Hamilton structure~\cite{CamassaHolm1993PRL}
 \begin{eqnarray*}
m_{t}=J_1\frac{\delta H_{CH,1}}{\delta m}=K_1\frac{\delta H_{CH,2}}{\delta m},\quad J_1=-\frac12 (m \partial+\partial m),\quad K_1=\partial^3-\partial
 \end{eqnarray*}
where
\begin{eqnarray*}
H_{CH,1}=\int_{\mathbb{R}}mu \mbox{d}x,\qquad
H_{CH,2}=\frac{1}{2}\int_{\mathbb{R}}(u^{3}+uu_{x}^{2}) \mbox{d}x,
 \end{eqnarray*}
and the sum of the two Hamiltonian
operators $J_1$
and $K_1$
is still a Hamiltonian operators.
It can be derived from the compatibility condition of the following
system for $\psi$ \cite{CamassaHolm1993PRL}:
\begin{eqnarray}\label{introdution1}
\psi_{x x}
=\left[\frac{1}{4}-\frac{m(x, t)}{2 \lambda}\right] \psi, \quad
\psi_{t}=-(\lambda+u) \psi_{x}+\frac{1}{2} u_{x} \psi.
\end{eqnarray}
System (\ref{introdution1}) enables one to solve the Cauchy problem of the CH equation
by the inverse scattering transform (IST)~\cite{ConstantinGerdjikovIvanov2006IP}.
Also, the   action angle variables of (\ref{CH}) can be constructed by using the IST\cite{BealsSattingerSzmigielski1998AM,
BoutetKostenkoShepelskyTeschl2009SJMA,  ConstantinGerdjikovIvanov2006IP,ConstantinMcKean1999CPAM}.
Despite these similarities,
the CH equation (\ref{CH}) has some
different properties from the KdV equation.
For example, there is an interesting phenomenon called  wave-breaking  in nature,
namely, the wave itself  remains bounded, while its slope becomes unbounded in finite time
\cite{BealsSattingerSzmigielski2000AM,ConstantinEscher1998AM,Whitham1980}.
This kind of singularities can not be derived from the KdV equation.
However, it can be modelled by Eq.~(\ref{CH}).
This is a significant difference between the CH and KdV equations.
Now that the wave-breaking phenomenon can occur in the solution of (\ref{CH}),
a natural problem is how the solution develops after the wave-breaking ?
This kind of problem  has been considered in \cite{BressanConstantin2007ARMA,
HoldenRaynaud2007CPDE}
and \cite{BressanConstantin2007AA,
HoldenRaynaud2009DCDS},
where the authors construct the global conservative solutions
and global dissipative solutions,
respectively.
On the other hand,
the CH equation (\ref{CH}) also admits a kind of solutions called the peakons
\cite{BealsSattingerSzmigielski2000AM,
CamassaHolm1993PRL,
Constantin2006IM,
ConstantinEscher2007BAMC,
ConstantinEscher1996AM,
ConstantinStrauss2000CPAM,
Lenells2005JDE}
whose shapes  are stable
under small perturbations
\cite{ConstantinStrauss2000CPAM,
Lenells2004IMRN}. The CH equation (\ref{CH}) is also related to geometry. Firstly, it represents
the geodesic flows. In the aperiodic case and when some asymptotic conditions are satisfied at
infinity,
equation (\ref{CH}) can be viewed as  the geodesic flow on
a manifold of diffeomorphim of the line
\cite{Constantin2000AIF}.
In the periodic case,
it can represent
the geodesic flow on the diffeomorphism group of the circle
\cite{ConstantinKolev2003CMH,
Kouranbaeva1999JMP}.
Secondly,
it can represent the families of
 pseudo-spherical surfaces
 \cite{GorkaReyes2011IMRN,
 Reyes2002LMP,
 GuiLiuOlverQu2013CMP}.
Thirdly, it appears  from a
non-stretching invariant planar curve flow in the
 centro-equiaffine geometry
  \cite{ChouQu2002PD,
  GuiLiuOlverQu2013CMP}.
The geometric illustration
of equation (\ref{CH}) yields
 the least action principle.
Other mathematical facts about equation
(\ref{CH}) include the existence of a
recursion operator
\cite{CamassaHolm1993PRL,
CamassaHolm1994AAM,
Fokas1995PD,
FokasFuchssteiner1981PD,
Fuchssteiner1996PD,
Fuchssteiner1982PTP},
a generalized Fourier transform
via squared eigenfunctions
\cite{ConstantinGerdjikovIvanov2007IP}
and so on. The local-in-time well-posedness for the initial value problem of (\ref{CH}) has been
achieved~\cite{ConstantinEscher1998ASNSPCS,
ConstantinEscher1998CPAM,
Danchin2001DIE,
Rodr2001NA}.
Equation (\ref{CH}) also admits global strong solutions~\cite{Constantin2000AIF,
ConstantinEscher1998ASNSPCS,
ConstantinEscher1998CPAM}, finite time blow-up strong solutions
\cite{Constantin2000AIF,
ConstantinEscher1998ASNSPCS,
ConstantinEscher1998CPAM,
ConstantinEscher1998AM}, as well as
global weak solutions
\cite{ConstantinMolinet2000CMP,
XinZhang2000CPAM,
XinZhang2002CPDE}.

Equation (\ref{CH}) has the quadratic nonlinearity,
while, there do exist other CH-type equations with cubic nonlinearity.
The first one of such equations is the modified CH (mCH) equation (alias the Fokas-Olver-Rosenau-Qiao (FORQ) equation)~\cite{Fokas1995PD,Fuchssteiner1996PD,OlverRosenau1996PRE,Qiao2006JMP}
\begin{eqnarray}\label{mCH}
m_{t}+[(u^{2}-u_{x}^{2}) m]_{x}=0, \quad m=u-u_{xx},
\end{eqnarray}
which corresponds to  Eq.~(\ref{gmCH-Novikov}) with $k_{2,3}=0$ and $k_{1}=1$. Eq.~(\ref{mCH}) was first derived
  by applying the tri-Hamiltonian duality to the
  bi-Hamiltonian representation
   of the modified KdV equation
\cite{Fuchssteiner1996PD,
  OlverRosenau1996PRE}.
  Later, it was rededuced
from the two-dimensional Euler equations
\cite{Qiao2006JMP}.
Similar to the CH equation (\ref{CH}),
the mCH equation (\ref{mCH})
 is also completely integrable.
Its  bi-Hamiltonian structure reads
\cite{OlverRosenau1996PRE,
QiaoLi2011TMP,
GuiLiuOlverQu2013CMP}
\begin{eqnarray*}
m_{t}=J_2 \frac{\delta H_{mCH,1}}{\delta m}=K_2 \frac{\delta H_{mCH,2}}{\delta m},\quad
J_2=-\partial_{x} m \partial_{x}^{-1} m \partial_{x},\quad K_2=\partial_{x}^{3}-\partial_{x},
\end{eqnarray*}
 where the  Hamiltonians are
 \begin{eqnarray*}
H_{mCH,1}=\int_{\mathbb{R}} m u \mbox{d} x, \quad H_{mCH,2}=\frac{1}{4} \int_{\mathbb{R}}(u^{4}+2 u^{2} u_{x}^{2}-\frac{1}{3} u_{x}^{4}) \mbox{d} x.
\end{eqnarray*}
One can solve equation (\ref{mCH})
using the inverse scattering transform since
it admits  the following Lax pair \cite{QiaoLi2011TMP, GuiLiuOlverQu2013CMP}:
$\psi_{t}=U\psi,\quad \psi_{x}=V\psi$ with $\psi=(\psi_{1}, \psi_{2})^{T}$ and
with
\begin{equation*}
U=\frac{1}{2}\left[\!\!\begin{array}{cc}
-1 & \lambda m \\
-\lambda m & 1\end{array}\!\!\right],\,\,
V=\left[\!\!\begin{array}{cc}
\lambda^{-2} +\frac{1}{2} \left(u^{2}-u_{x}^{2}\right) & -\lambda^{-1}\left(u- u_{x}\right)-\frac{1}{2} \lambda\left(u^{2}-u_{x}^{2}\right) m \\
\lambda^{-1}\left(u+ u_{x}\right)+\frac{1}{2} \lambda\left(u^{2}-u_{x}^{2}\right) m & -\lambda^{-2} -\frac{1}{2} \left(u^{2}-u_{x}^{2}\right)
\end{array}\!\!\right].
\end{equation*}
From the viewpoint of geometry,
equation (\ref{mCH}) arises from an intrinsic
(arc-length preserving) invariant planar curve flow in Euclidean geometry
\cite{GuiLiuOlverQu2013CMP}.
The local-wellposedness, blow-up and
wave breaking problems   of equation (\ref{mCH})
were considered in
\cite{ChenLiuQuZhang2015AM,
FuGuiLiuQu2013JDE,
GuiLiuOlverQu2013CMP,
LiuQuZhang2014AA}.
The explicit form of the multipeakons of  equation (\ref{mCH})
has been addressed in
\cite{ChangSzmigielski2018CMP}.
The peakon solution of equation (\ref{mCH})
is orbital stability
 \cite{LiuLiuQu2014AM,
 QuLiuLiu2013CMP}.
 The H\"{o}lder continuity of equation (\ref{mCH})
has been discussed in
\cite{HimonasMantzavinos2014JNS}.

Another celebrated CH-type equation with cubic nonlinearity is the Novikov equation~\cite{Novikov2009JPA}
\begin{eqnarray}\label{Novikov}
m_{t}
+u^{2} m_{x}+3 u u_{x} m=0,\,\, \,\, m=u-u_{xx},
\end{eqnarray}
which was discovered by Novikov when classifying the nonlocal partial differential equations in the light of symmetry
\cite{Novikov2009JPA}.
Like the CH  and mCH equations, the Novikov equation (\ref{Novikov}) is also completely
integrable, and its
Lax pair reads \cite{Novikov2009JPA}
\begin{eqnarray*}
\begin{cases}
 \psi_{x x x}= \psi_{x}+\lambda m^{2} \psi+2  \dfrac{m_{x}}{m} \psi_{x x}+\dfrac{m m_{x x}-2 m_{x}^{2}}{m^{2}} \psi_{x}, \v\\
\psi_{t}=\dfrac{1}{\lambda} \dfrac{u}{m} \psi_{x x}-\dfrac{1}{\lambda} \dfrac{(m u)_{x}}{m^{2}} \psi_{x}-u^{2} \psi_{x}.
\end{cases}
\end{eqnarray*}
or  \cite{HoneWang2008JPAMT}
$\Psi_{x}=A \Psi, \quad \Psi_{t}=B\Psi$, where $\Psi=(\Psi_{1}, \Psi_{2}, \Psi_{3})^{T}$ and
the two matrices $A$ and $B$ are defined as
\begin{equation*}
A=\left[\begin{array}{ccc}
0 & m \lambda & 1 \\
0 & 0 & m \lambda \\
1 & 0 & 0
\end{array}\right], \quad B=\left[\begin{array}{cccc}
\frac{1}{3 \lambda^{2}}-u u_{x} & \frac{u_{x}}{\lambda}-u^{2} m \lambda & u_{x}^{2} \v\\
\frac{u}{\lambda} & -\frac{2}{3 \lambda^{2}} & -\frac{u_{x}}{\lambda}-u^{2} m \lambda \v\\
-u^{2} & \frac{u}{\lambda} & \frac{1}{3 \lambda^{2}}+u u_{x}
\end{array}\right].
\end{equation*}
The bi-Hamiltonian structure of equation (\ref{Novikov}) was  provided
by  \cite{HoneWang2008JPAMT}:
\begin{eqnarray*}
m_{t_{k}}=J_3 \frac{\delta H_{k+6}}{\delta m}=K_3\frac{\delta H_{k}}{\delta m},
\end{eqnarray*}
where the Hamiltonian operator pairs are
\begin{eqnarray*}
J_3=-2\left(3 m \partial_{x}+2 m_{x}\right)\left(4 \partial_{x}-\partial_{x}^{3}\right)^{-1}\left(3 m \partial_{x}+m_{x}\right),\,\,
K_3=\left(1-\partial_{x}^{2}\right) m^{-1} \partial_{x} m^{-1}\left(1-\partial_{x}^{2}\right)
\end{eqnarray*}
and the first few Hamiltonians are given by
\begin{eqnarray*}
H_{1}=\int_{\mathbb{R}} \frac{1}{8}\left(u^{4}+2 u^{2} u_{x}^{2}-\frac{u_{x}^{4}}{3}\right) \mbox{d} x,\quad
H_{5}=\int_{\mathbb{R}} m^{2 / 3}  \mbox{d}  x,\quad
 H_{7}=\int_{\mathbb{R}} \frac{1}{3}\left(m^{-8 / 3} m_{x}^{2}+9 m^{-2 / 3}\right)  \mbox{d}  x
\end{eqnarray*}
with $k=(\pm1 \mod 6).$

The periodic  initial value  problem
for equation (\ref{Novikov}) in $H^{s}$
with $s>5/2$
was   considered  in
\cite{Tiglay2011IMRN}
  based on Arnold's geometric framework.
The index $s>5/2$ was improved to
$s>3/2$ in
\cite{HimonasHolliman2012Non}
based on a Galerkin-type approximation.
Moreover, the continuous  dependence
on the initial data is showed to be optimal in \cite{HimonasHolliman2012Non}.
The initial value problem  on the line was proved to be local well-posedness
in $B_{2,r}^{s}$ in
\cite{NiZhou2011JDE}
in light of  Kato's semigroup theory.
\cite{YanLiZhang2012JDE} also considered the initial value problem
in the Besov circumstances.
The explicit form of  multipeakons for
(\ref{Novikov})
 was computed in
\cite{HoneLundmarkSzmigielski2009DPDE}.
The wave breaking phenomenon for (\ref{Novikov})
was discussed in
\cite{ChenGuoLiuQu2016JFA,
JiangNi2012JMAA,
YanLiZhang2013NDEA}
under  the sign-changing or non-sign-changing conditions.
(\ref{Novikov})  admits
global weak solution
\cite{Lai2013JFA,
Lai2013JMAA,
WuYin2011JPA}.
The measure of momentum support with upper and lower bounds
after time $t$ was researched in
\cite{Guo2017NARWA}. Concerning the CH type equations with analytic initial data, Barostichi, {\it et al}~\cite{BarostichiHimonasPetronilho2016JFA} established  an  Ovsyannikov type theorem for an autonomous abstract Cauchy problem.
Then they apply it to some CH type equations (e.g, CH, mCH, DP, and the Novikov equations)   with
initial data in spaces of analytic functions
to obtain the analytic lifespan of
the solution as well as the
continuity of the data-to-solution map in spaces
of analytic functions.

Eq.~(\ref{gmCH-Novikov}) with $k_2=0$ becomes the following mCH-CH equation with quadratic and cubic nonlinearities~\cite{Fokas1995PD}
\begin{eqnarray}\label{mCH-CH}
m_{t}+k_{1}[(u^{2}-u_{x}^{2}) m]_{x}+k_3(um_{x}+2u_xm)=0,\,\,\,\, m=u-u_{xx},
\end{eqnarray}
which  is completely integrable and admits the Lax pair and bi-Hamiltonian structure~\cite{qiao2012,xia2018}, which is a linear combination of those of the  CH and mCH equations, that is,
 \begin{eqnarray*}
m_{t}=J_4 \frac{\delta H_{mCH-CH,1}}{\delta m}=K_4 \frac{\delta H_{mCH-CH,2}}{\delta m},\quad
J_4=k_1J_2+k_1J_1,\quad K_4=K_2=K_1,
\end{eqnarray*}
 where the  Hamiltonians are
 \begin{eqnarray*}
H_{mCH-CH,1}=H_{mCH,1}=H_{CH,1}, \quad
H_{mCH-CH,2}=k_1H_{mCH,2}+k_2H_{CH,2}.
\end{eqnarray*}
The  Cauchy problem and blow-up  of Eq.~(\ref{mCH-CH}) were studied in ~\cite{Liu2014CP,ChenGuoLiuQu2016JFA,liu2016NA}.

Eq.~(\ref{gmCH-Novikov}) with $k_3=0$ reduces to  the mCH-Novikov equation with cubic nonlinearity derived recently in~\cite{ChenHuLiu},
\begin{eqnarray}\label{mCH-Novikov}
m_{t}
+k_{1}((u^{2}-u_{x}^{2}) m)_{x}
+k_{2}(u^{2} m_{x}+3 u u_{x} m)=0,\,\,\,\, m=u-u_{xx}.
\end{eqnarray}
It admits two conserved quantities
\begin{eqnarray*}
H_{1}(u)=\int_{\mathbb{R}}mu \mbox{d} x, \quad H_{2}(u)=\int_{\mathbb{R}}\left(u^{4}+2 u^{2} u_{x}^{2}-\frac{1}{3} u_{x}^{4}\right) \mbox{d} x
\end{eqnarray*}
and a Hamiltonian structure
\begin{eqnarray*}
m_t=J \frac{\delta H_{1}}{\delta m},\qquad  J=k_{1}J_2+k_2J_3.
\end{eqnarray*}
Recently, Mi {\it et al}~\cite{MiLiuHuangGuo2020JDE} obtained the local well-posedness
 of Eq.~(\ref{mCH-Novikov})  in noncritical or critical Besov spaces and
provided the blow-up criterion in Sobolev spaces
as well as the wave breaking condition, the persistence property
and the analyticity of the solution.

As far as we know,
the initial value problem (\ref{gmCH-Novikov})
has not been investigated  yet.
 Employing the Littlewood-Paley theory
and the  transport equations theory,
we will first show  the local well-posedness for the strong
solutions of equation  (\ref{gmCH-Novikov})
 in Besov spaces following the spirit of
\cite{Danchin2001DIE,
Danchin2003,
Danchin2003JDE}.
Then,
the H\"{o}lder continuity of the
data-to-solution map of this equation will
be shown.
A Moser-type estimate
in the Sobolev spaces will help us
to establish a blow-up criterion
and  the precise blow-up quantity
for (\ref{gmCH-Novikov}).
We will also
 provide two sufficient conditions on the initial
data  to ensure the occurance of the wave breaking phenomenon.
The first condition is obtained by
only consider a portion of the blow-up quantity,
while the second one is derived by
inspecting the whole  blow-up quantity.
We will finally explore the non-periodic peakons, periodic peakons and
multi-peakon solutions for this equation.

Note that equation (\ref{gmCH-Novikov})
admits a conserved quantity $H_{1}=\int_{\mathbb{R}} m u\mbox{d} x$
and a Hamiltonian structure
\begin{eqnarray}
m_t=J^{\prime} \frac{\delta H_{1}}{\delta m},\qquad J^{\prime}=k_1J_2+k_2J_3+k_3J_1.
\end{eqnarray}

For convenience,
the solution spaces $E_{p, r}^{s}(T)$ are defined as follows
\begin{equation*}
E_{p, r}^{s}(T) \triangleq\left\{\begin{array}{l}
C\left([0, T) ; B_{p, r}^{s}\right) \cap C^{1}\left([0, T) ; B_{p, r}^{s-1}\right), \text { if } r<\infty, \vspace{0.1in}\\
C_{w}\left([0, T) ; B_{p, \infty}^{s}\right) \cap C^{0,1}\left([0, T) ; B_{p, \infty}^{s-1}\right), \text { if } r=\infty.
\end{array}\right.
\end{equation*}

We next introduce some notations to be used in this paper.

\textbf{Notation.}
Let $p(x)=\frac{1}{2}e^{-|x|},\, x\in\mathbb{R}$ be the fundamental solution of $1-\partial_{x}^{2}$
on $\mathbb{R}$, and two convolution operators
 $p_{\pm}$ be~\cite{ChenGuoLiuQu2016JFA}
\begin{eqnarray}
p_{\pm} * f(x)=\frac{1}{2}e^{\mp x} \int_{-\infty}^{\pm x} e^{y} f(\pm y) d y,\label{pPlus}
\end{eqnarray}
where the star denotes the spatial convolution.
Note that $p$ and $p_{\pm}$ satisfy
\begin{equation}\label{pPlusMinus}
p=p_{+}+p_{-}, \quad p_{x}=p_{-}-p_{+},\quad
p_{xx}*f=p*f-f.
\end{equation}
Let  $\mathcal{S}$ stands for the Schwartz space
and
$\mathcal{S}^{\prime}$ represents the spaces of  temperate distributions.
Let $L^{p}(\mathbb{R})$ be the Lebesgue space  equipped with the norm $\|\cdot\|_{L^{p}}$ for $1 \leq p \leq \infty$ and
  $H^{s}(\mathbb{R})$ be the Sobolev space equipped with the
norm $\|\cdot\|_{H^{s}}$ for
 $s \in \mathbb{R}$.\\

Our first result is the local well-posedness of
the Cauchy problem (\ref{gmCH-Novikov}) in $E_{p, r}^{s}(T)$.

\begin{thm}\label{WaveBreakingThm1}
Suppose $u_{0} \in B_{p, r}^{s}$ with $s>\max \{5/2,\, 2+1/p\}$ and $p, r \in[1, \infty]$. Then there exists a time $T>0$
and a unique solution $u \in E_{p, r}^{s}(T)$ of the Cauchy problem (\ref{gmCH-Novikov})  such that the data-to-solution
map $u_{0} \mapsto u:$ $B_{p, r}^{s} \mapsto C([0, T] ; B_{p, r}^{s^{\prime}}) \cap C^{1}([0, T] ; B_{p, r}^{s^{\prime}-1})$ is continuous for every $s^{\prime}<s\, (=s)$ as $r=\infty\, (<\infty).$
\end{thm}
Theorem \ref{WaveBreakingThm1}
implies the following corollary:

\begin{cor}\label{cor1}
Let $s>1/2.$
Suppose $m_{0}=(1-\partial_{x}^{2}) u_{0} \in H^{s}(\mathbb{R}).$
 Then there exists a time $T^{*}>0$ such that the Cauchy problem (\ref{gmCH-Novikov}) admits a unique strong solution
$m \in C([0, T^{*}] ; H^{s}) \cap C^{1}([0, T^{*}] ; H^{s-1})$ and the map $m_{0} \mapsto m$ is continuous from a neighborhood of $m_{0}$ in $H^{s}$ into $C([0, T^{*}] ; H^{s}) \cap C^{1}([0, T] ; H^{s-1})$.
\end{cor}

The above corollary has shown that
the data-to-solution map
$u(0)\in H^{s}\rightarrow u(t)\in H^{s}$
for the Cauchy problem (\ref{gmCH-Novikov}) is continuous.
Using similar method as in  \cite{HimonasMantzavinos2014NA} for the mCH equation,
one can show the non-uniform dependence of this map.
Nevertheless, this map is H\"{o}lder continuous in $H^{s},$
as will be established in the following result

\begin{thm}\label{HolderContinuityThm} Let $0 \leq r<s$ with $s>5 / 2$, then the data-to-solution map for the  Cauchy  problem
(\ref{gmCH-Novikov}) is H\"{o}lder continuous
 in $H^{s}$ under the $H^{r}$ norm.
Precisely speaking,
for initial data $u_{0}, v_{0}$
with $\|u_{0}\|_{H^{s}}\leq \rho$
and $\|v_{0}\|_{H^{s}}\leq \rho,$
 the corresponding solutions $u, v$ of Eq.~(\ref{gmCH-Novikov}) satisfy
\begin{eqnarray}\label{HolderContinuityThm1}
\|u-v\|_{C\left([0, T] ; H^{r}\right)}
\leq C\left\|u_{0}-v_{0}\right\|_{H^{r}}^{\beta},
\end{eqnarray}
where  the constant $C=C(s, r, k_{1}, k_{2}, k_{3},\rho)$
and
the exponent   $\beta$ is given by
\begin{equation}\label{HolderContinuityThm2}
\beta=\left\{\begin{array}{ll}
1, & (s, r) \in D_{1} \\
(2 s-3)/(s-r), & (s, r) \in D_{2}, \\
(s-r)/2, & (s, r) \in D_{3}, \\
s-r, & (s, r) \in D_{4},
\end{array}\right.
\end{equation}
where the regions $D_j\, (j=1,2,3,4)$ in the $(s, r)$-plane are defined  by
$$\left\{
\begin{array}{l}
D_{1}=\{(s, r)\,| \,0 \leq r \leq 1.5,\,\, 3-s \leq r \leq s-2\} \cup\{(s, r)\,| \, 1.5< r \leq s-1\}, \v\\
D_{2}=\{(s, r)\,| \, 2.5<s<3,\,\, 0 \leq r \leq-s+3\}, \v\\
D_{3}=\{(s, r)\,| \, 2.5<s,\,\, s-2 \leq r \leq 1.5\}, \v\\
D_{4}=\{(s, r)\,| \, 2.5<s,\,\,  s-1 \leq r<s\}.
\end{array}\right.
$$
\end{thm}

Having  established the  local well-posedness result, it is a natural problem whether the local solution can be extended
 globally  or not. Our third result is about the   blow-up criterion.
\begin{thm}\label{WaveBreakingThm2}
Assume that $m_{0} \in H^{s}(\mathbb{R})$ with $s>\frac{1}{2}.$
Let  $T^{*}>0$ be the maximal existence time of the solution $m$ to the Cauchy problem $(\ref{gmCH-Novikov}).$ Then
\begin{equation}\label{blowup00}
  T^{*}<\infty\quad \Rightarrow\quad \int_{0}^{T^{*}}\|m(t)\|_{L^{\infty}}^{2} \mbox{d} t=\infty.
\end{equation}
\end{thm}

\begin{rem}\label{rem1}
Suppose $m_{0}\in H^{s}$  with $s>\frac{1}{2}$.
Let $m_{s}$ be the  maximal $H^{s}$ solution corresponding to
$m_{0}$ with lifespan $T_{s}.$
Let $m_{s^{\prime}}$ be the
maximal $H^{s^{\prime}}$ solution corresponding to
 the same initial data with lifespan $T_{s^{\prime}}$
for  $\frac{1}{2}<s^{\prime} <s$.
Since $H^{s}\hookrightarrow H^{s^{\prime}}$,
the uniqueness ensures that $T_{s}\leq T_{s^{\prime}}$
and $m_{s}=m_{s^{\prime}}$ on $[0,T_{s})$.
If $T_{s}<  T_{s^{\prime}},$
then we can conclude
$m_{s^{\prime}}\in C([0,T_{s}]; H^{s^{\prime}})$
and consequently
$m_{s^{\prime}}\in L^{2}([0,T_{s}]; L^{\infty}),$
contradicting the above blow up criterion (\ref{blowup00}).
That is to say, the lifespan $T^{*}$ does not depend on the regularity index $s$ of the initial data $m_{0}.$
\end{rem}

\begin{rem} In fact, Theorems \ref{WaveBreakingThm1}-\ref{WaveBreakingThm2} also hold true for the mCH-Novikov-CH equation with the linear dispersive term
\begin{eqnarray}\begin{cases}\label{gmCH-Novikov-d}
m_{t}+\kappa u_x+k_{1}[(u^{2}-u_{x}^{2})m]_{x} +k_{2}(u^{2} m_{x}+3 u u_{x} m)+k_{3}(u m_{x}+2u_{x} m)=0, \quad x\in\mathbb{R},\,\, t>0, \\
u(0, x)=u_{0}(x),\quad x\in\mathbb{R},
\end{cases}
\end{eqnarray}
where $\kappa,\, k_{1}, k_{2}$ and $k_{3}$ are all real-valued parameters, and $m=u-u_{xx}$.

\end{rem}

The following theorem shows the precise blow-up scenario for sufficiently regular solutions
to the Cauchy problem (\ref{gmCH-Novikov}).
\begin{thm}\label{blowupQuantityThm}
Suppose $u_{0} \in H^{s}(\mathbb{R})$
with $s>\frac{5}{2}$.
Let  $T^{*}>0$ be the maximal existence time of the solution
 $u$ to the Cauchy problem (\ref{gmCH-Novikov}).
  Then the solution $u$ blows up in finite time if and only if
\begin{eqnarray}\label{blowupQuantity}
  \liminf _{t \rightarrow T^{*}}\left(\inf _{x \in \mathbb{R}}\left(2 k_{1} u_{x} m(t, x)+3 k_{2} u u_{x}(t, x)+2k_{3}u_{x}(t,x)\right)\right)=-\infty.
\end{eqnarray}
\end{thm}

From Theorem \ref{blowupQuantityThm},
we find that the solution blows up in  finite time if and only if  the quantity
$M(t,x)=2 k_{1} u_{x} m(t, x)+3 k_{2} u u_{x}(t, x)+2k_{3}u_{x}(t,x)$
is unbounded below.
We next  prove that
this quantity
admits an uniform upper bound
 as long as
$m(t, x)$  exists
under the assumption of Theorem \ref{blowupQuantityThm}.
\begin{thm}\label{finalThm}
Let the assumption of Theorem \ref{blowupQuantityThm}
be satisfied.
Suppose also that $m_{0}(x)=\left(1-\partial_{x}^{2}\right) u_{0} \geq 0$ for all $x \in \mathbb{R},$ and $m_{0}\left(x_{0}\right)>0$ at some point $x_{0} \in \mathbb{R}.$
 Then there holds
\begin{eqnarray}\label{MuniformUpperBound}
\sup _{x \in \mathbb{R}}M(t,x)
\leq
2 k_{1}\|u_{0}\|_{H^{1}} \sup_{x\in \mathbb{R}}m_{0}(x)
+3 k_{2} \|u_{0}\|_{H^{1}}^{2}
+2k_{3}\|u_{0}\|_{H^{1}}
\end{eqnarray}
for all $t \in[0, T^{*}).$
\end{thm}

With  the precise blow-up quantity (\ref{blowupQuantity})
at hand,
 we next display  some sufficient conditions
  for the occurance of the wave-breaking phenomenon
to the Cauchy problem (\ref{gmCH-Novikov})
  in the case of a non-sign-changing  momentum.
The  result reads:

\begin{thm}\label{WaveBreakingThm3}
Suppose $k_1>0,\, k_2,\, k_3\geq 0,\, u_{0}(x) \in H^{s}(\mathbb{R})$ with
$s>\frac{5}{2}$ and  $m_{0}(x) \geq 0.$
Let one of the following four cases be satisfied:

{\it Case 1:} For $k_{1},\, k_{2},\, k_{3}>0$, there exists some   $x_{0} \in \mathbb{R}$ and $\alpha>2+\frac{4k_{1}}{3k_{2}}$ such that
$m_{0}(x_{0})>0$
and
\begin{eqnarray}\label{initialGradient}
u_{0, x}(x_{0})
<
-\frac{\alpha}{\sqrt{2}}\sqrt{\frac{2k_{1}+3k_{2}}{2k_{1}}}
\bigg[u(x_{0})+\frac{3k_{2}\gamma+3k_{3}}{2(2k_{1}+3k_{2})}\bigg],
\end{eqnarray}
where
\begin{eqnarray*}
\gamma=\frac{(3k_{2}k_{3}+4k_{1}k_{3})\pm 2k_{3}
 \sqrt{4k_{1}^{2}+6k_{1}k_{2}}}{3k_{2}^{2}}.
\end{eqnarray*}

{\it Case 2:} For $k_{1},\, k_3>0,\, k_{2}=0$, there exists some   $x_{0} \in \mathbb{R}$ and $m_{0}(x_{0})>0$
such that~\cite{ChenGuoLiuQu2016JFA}
\begin{eqnarray}\label{initialGradient-1}
u_{0, x}(x_{0})
<
-\frac{\alpha}{\sqrt{2}}
\bigg[u(x_{0})+\frac{3k_{3}}{4k_{1}}\bigg].
\end{eqnarray}

{\it Case 3:} For $k_{1},\, k_{2}>0, k_{3}=0$, there exists some   $x_{0} \in \mathbb{R}$ and $m_{0}(x_{0})>0$
such that
\begin{eqnarray}\label{initialGradient-2}
u_{0, x}(x_{0})
< -\frac{1}{\sqrt{2}}
\sqrt{\frac{2k_{1}+3k_{2}}{2k_{1}}}
u(x_{0}).
\end{eqnarray}

{\it Case 4:} For $k_{1}>0, k_{2}=k_{3}=0$, there exists some   $x_{0} \in \mathbb{R}$ and $m_{0}(x_{0})>0$
such that~\cite{ChenGuoLiuQu2016JFA}
\begin{eqnarray}\label{initialGradient-3}
u_{0, x}(x_{0})
< -\frac{1}{\sqrt{2}}
u(x_{0}).
\end{eqnarray}
Then the solution  $u(t, x)$ to the Cauchy problem (\ref{gmCH-Novikov})
will blow up at time   $T^{*}$  with
$$
T^{*} \leq-\frac{1}{2 k_{1} m_{0}\left(x_{0}\right) u_{0, x}\left(x_{0}\right)}.
$$
\end{thm}

\begin{rem}\label{rem2}
Similar to the proof of the wave-breaking of the generalized mCH equation considered in \cite{ChenGuoLiuQu2016JFA}, for Case 1: $k_1,\,k_2,\,k_3>0$, we trace the dynamics of the ratio $\widehat{u}_{x}/(\widehat{u}+\gamma)$ along
the characteristics with some
nonnegative parameter $\gamma$ to be settled to make the
constant term in the brackets of Eq.~(\ref{character10}) to be zero in the proof of Theorem \ref{WaveBreakingThm3}.
Having determined the parameter $\gamma$ as in Eq.~(\ref{character11}), we need the two conditions (\ref{character13}) and (\ref{character14})
to be satisfied simultaneously.
In \cite{ChenGuoLiuQu2016JFA},
the initial state of the corresponding
condition (\ref{character13}) can ensure
the corresponding
condition (\ref{character14}) exactly.
However, we have to multiply
some factor $\alpha$
on the condition of the initial gradient
to  guarantee  the condition (\ref{character14})
in the present paper, as we see
 on  the right hand side of (\ref{initialGradient}).
\end{rem}

Theorem \ref{WaveBreakingThm3} has provided some sufficient conditions on the initial data
for the breaking of waves, however, it  contains no  information of  the blow-up rate.
Using a different method, one can establish another wave-breaking Theorem which includes the estimation
of the blow-up rate.
The result reads
\begin{thm}\label{WaveBreakingThm4}
Suppose $k_{1}, k_{2}, k_{3}>0,\,  s>\frac{5}{2},\, u_{0} \in H^{s}(\mathbb{R})$ and  $m_{0} \geq 0 .$
There exists some   $x_{1} \in \mathbb{R}$
and
$C_{2}>0$
such that if
$m_{0}(x_{1})>0,$
$u_{0}(x_{1})+1\leq C_{2}m_{0}(x_{1})$
and
\begin{eqnarray}\label{finalCondition}
&&\left[
2k_{1}+\frac{3k_2u_{0}(x_{1})+2k_3}{m_{0}(x_{1})}\right]
\partial_{x}u_{0}(x_{1})\nonumber\\
&&\quad <-\left[
\frac{(2k_{1}+(3k_{2}+2k_{3})C_{2})
(\frac{38}{3}k_{1}+23k_{2}+\frac{25}{2}k_{3})
(\|u_{0}\|_{H^{1}}^{2}+\|u_{0}\|_{H^{1}}^{3})}
{m_{0}(x_{1})}
\right]^{1/2},
\end{eqnarray}
then, the solution  $u(t, x)$ to the Cauchy problem (\ref{gmCH-Novikov})
will blow up at time   $T_{*}$  with
$$
T_{*}
\leq t_-:=\frac{C_{0}}{2 C_{3}}-\frac{1}{2} \sqrt{\left(\frac{C_{0}}{C_{3}}\right)^{2}-\frac{2}{C_{3} m_{0}\left(x_{1}\right)}},
$$
where $C_{0},\, C_1$ and $C_3$ are given by
\begin{eqnarray}\no
C_{0}=-\partial_{x}u_{0}(x_{1})\left(2k_{1}\!+\!\frac{3k_{2}u_{0}(x_{1})+2k_{3}}{m_{0}(x_{1})}\right),
\quad
C_{1}=\left(\frac{38}{3}k_{1}+23k_{2}+\frac{25}{2}k_{3}\right)\left(\|u_{0}\|_{H^{1}}^{2}+\|u_{0}\|_{H^{1}}^{3}\right),
\end{eqnarray}
and $C_{3}=\frac{1}{2}[2k_{1}+(3k_{2}+2k_{3})C_{2}] C_{1}$,
respectively.

Furthermore,  when
$T_{*}=t_-$,
 one has the following estimation  of the blow-up rate
\begin{eqnarray}\label{blowupRate}
\liminf _{t \rightarrow T_{*}^{-}}\left(\left(T_{*}-t\right) \inf _{x \in \mathbb{R}} M(t, x)\right) \leq -\frac{1}{2},
\end{eqnarray}
where $M(t, x)$ is the blow-up quantity defined in (\ref{M}).
\end{thm}
\begin{rem}\label{rem3}
The result obtained in Theorem \ref{WaveBreakingThm4}
is inspired by \cite{GuiLiuOlverQu2013CMP,MiLiuHuangGuo2020JDE}.
However, we need the additional condition
$u_{0}(x_{1})+1\leq C_{2}m_{0}(x_{1})$
as stated in  Theorem \ref{WaveBreakingThm4}.
This condition  essentially comes from the
terms associated with the
coefficients  $k_{2}$ and $k_{3}$
in equation (\ref{gmCH-Novikov}) and is used  to
make the differential inequality
 (\ref{mrecast}) to be integrated not that hard.
\end{rem}

According to the relation
$u=\left(1-\partial_{x}^{2}\right)^{-1} m=p * m$
with $p$ defined in the \textbf{Notation} part, Eq.~(\ref{gmCH-Novikov}) can be recast as the weak form
\begin{eqnarray} \label{gmchn}
&&u_{t}
+\frac{k_{1}+k_{2}}{3}\partial_{x}(u^{3})
-\frac{k_1}{3}u_{x}^{3}+\frac{k_3}{2}\partial_{x}(u^{2})\nonumber\\
&&+p_{x}*\left[\left(\frac{2k_1}{3}+k_{2}\right)u^{3}
+\left(k_{1}+\frac{3k_2}{2}\right)uu_{x}^{2}+k_{3}u^{2}+\frac{k_3}{2}u_{x}^{2}\right]
+p*\left[\left(\frac{k_1}{3}+\frac{k_2}{2}\right)u_{x}^{3}\right]=0,\qquad
\end{eqnarray}
which  allows us to define the  weak solution
for (\ref{gmCH-Novikov}) as follows:
\begin{definition}\label{defOfWeakSolution}
For  $u_{0} \in W^{1,3}(\mathbb{R})$, the function $u(t,x) \in L_{l o c}^{\infty}\left([0, T), W_{l o c}^{1,3}(\mathbb{R})\right)$ is said to be a weak solution to the Cauchy  problem (\ref{gmCH-Novikov}) if it satisfies the following
integral equality
\begin{align}\label{defOfWeakSolution1}
&\int_{0}^{T} \int_{\mathbb{R}}
\bigg\{u \varphi_{t}+\frac{k_{1}+k_{2}}{3} u^{3} \varphi_{x}+\frac{k_1}{3}u_{x}^{3} \varphi
+\frac{k_3}{2}u^{2}\varphi_{x}-p*\left[\left(\frac{k_1}{3}+\frac{k_2}{2}\right)u_{x}^{3}\right]\varphi\nonumber\\
&\quad+p*\left[\left(\frac{2k_1}{3}+k_{2}\right)u^{3}+\left(k_{1}+\frac{3k_2}{2}\right)uu_{x}^{2}
+k_{3}u^{2}+\frac{k_3}{2}u_{x}^{2}\right] \varphi_x
\bigg\}\mbox{d} x \mbox{d} t
+\int_{\mathbb{R}} u_{0}(x) \varphi(0, x) \mbox{d} x=0
\end{align}
for any smooth test function $\varphi(t, x) \in C_{c}^{\infty}([0, T) \times \mathbb{R}).$
Moreover, $u(t,x)$ is called a global weak solution if
$T$ can be taken arbitrarily large.
\end{definition}

In the following four propositions we would like to give the non-periodic and periodic peakon and multi-peakon solutions of Eq.~(\ref{gmCH-Novikov}).

\begin{prop}\label{exsitenceOfPeakon}
For $a\not=0$ and $x_0\in\mathbb{R}$,
 Eq.~(\ref{gmCH-Novikov}) admits a global weak peakon solution, in the sense of Definition \ref{defOfWeakSolution}, of the form
\begin{eqnarray}\label{singlePeakon}
u_{a}(t, x)=a e^{-|x-c t-x_0|},
\end{eqnarray}
with the velocity $c$  satisfying
\begin{eqnarray}\label{solu-c}
c=\left(\frac{2k_1}{3}+k_{2}\right)a^{2}+k_{3}a.
\end{eqnarray}
\end{prop}

\begin{rem} \label{peakon-r1}
For  real parameters $a,\,c$ satisfying  Eq.~(\ref{solu-c}), the peakon solution (\ref{singlePeakon}) is a bright right-going travelling wave solution for $a>0$ and $c>0$, a dark right-going travelling wave solution for $a<0$ and $c>0$, a bright left-going travelling wave solution for $a>0$ and $c<0$, and a dark left-going travelling wave solution for $a<0$ and $c<0$.
\end{rem}

\begin{rem} \label{peakon-r2}
In particular, one can reduce the peakon solution given by Eqs.~(\ref{singlePeakon})-(\ref{solu-c})  to those of the CH equation for $k_1=k_2=0,\,k_3=1$, where $a=c$~\cite{CamassaHolm1993PRL}; the mCH equation for $k_2=k_3=0,\,k_1=1$, where $a=\pm\sqrt{3c/2}$~\cite{GuiLiuOlverQu2013CMP}; the Novikov equation for $k_1=k_3=0,\,k_2=1$, where $a=\pm \sqrt{c}$~\cite{HoneWang2008JPAMT}; the mCH-CH equation for $k_2=0,\,k_1k_3\not=0$, where $a=(-3k_3\pm\sqrt{9k_3^2+24ck_1})/(4k_1)$~\cite{Liu2014CP,qiao2012}; the mCH-Novikov equation for $k_1k_2\not=0,\,k_3=0$, where $a=\pm \sqrt{3c/(2k_1+3k_2)}$~\cite{MiLiuHuangGuo2020JDE} (Notice that we here correct it); or the Novikov-CH equation for $k_1=0,\, k_2k_3\not=0$, where $a=(-k_3\pm\sqrt{k_3^2+4ck_2})/(2k_2)$.
\end{rem}

\begin{rem} \label{peakon-r3} For the given real velocity $c$
and for $k_1k_2k_3\not=0$,
 one finds from Eq.~(\ref{solu-c}) that
\begin{eqnarray}
a=\left\{\begin{array}{ll}
 \dfrac{-3k_3 \pm \sqrt{9k_3^2+12c(2k_1+3k_2)}}{2(2k_1+3k_2)}, &  {\rm as}\,\,\, 2k_1+3k_2\not=0, \\
 c/k_3, & {\rm as}\,\,\, 2k_1+3k_2=0,\, k_1k_2k_3\not=0.
 \end{array}
\right.
\end{eqnarray}
Therefore, (i) if $2k_1+3k_2\not=0,\, k_1k_2k_3\not=0$, and $9k_3^2+12c(2k_1+3k_2)<0$, then $a$ is complex such that the solution given by Eqs.~(\ref{singlePeakon})-(\ref{solu-c}) is a complex peakon solution; (ii) If $2k_1+3k_2=0,\, k_1k_2k_3\not=0$, then the peakon solution of Eq.~(\ref{gmCH-Novikov}) is of the form $u_a(t,x)=\frac{c}{k_3}e^{-|x-ct-x_0|}$, which is independent of $k_1,\, k_2$.
\end{rem}

Similarly, the non-periodic multi-peakon solutions of Eq.~(\ref{gmCH-Novikov}) in the sense of Definition \ref{defOfWeakSolution} is
\begin{prop}\label{multipeakon} The multi-peakon solutions of  Eq.~(\ref{gmCH-Novikov}) in the sense of Definition \ref{defOfWeakSolution} are given by
\begin{eqnarray}\label{mp}
u_{npm}(t, x)=\sum_{i=1}^{N} p_{i}(t) e^{-\left|\xi_i\right|},\quad \xi_i=\xi_i(x,t)=x-q_i(t),
\end{eqnarray}
where the time-dependent position functions $q_i(t)$ (without loss of generality $q_1(t)<q_2(t)<\cdots<q_N(t)$ is required) and amplitudes $p_i(t)$ satisfy the dynamical system
\begin{eqnarray}\label{pDot}
\left\{\begin{array}{l}
\dot{p}_{i}(t)= \displaystyle p_{i}(t)\sum_{j=1}^{N}p_{j}(t)\operatorname{sgn}(q_{i}-q_{j})e^{-|q_{i}-q_{j}|}
\left[k_2\left(\sum_{k=1}^{N}p_{k}(t)e^{-|q_{i}-q_{k}|}\right)+k_3\right], \v\\
\dot{q}_{i}(t)= \displaystyle -\frac{k_1}{3}p_i^2(t)-k_{1}\!\left(\sum_{j=1}^{N}p_{j}(t)
 {\rm sgn}(q_i-q_j)e^{-|q_{i}-q_{j}|}\right)^2+(k_{1}+k_{2})\!\left(\sum_{j=1}^{N}p_{j}(t)e^{-|q_{i}-q_{j}|}\right)^2 \v\\
\qquad\qquad \displaystyle+k_{3}\sum_{j=1}^{N}p_{j}(t)e^{-|q_{i}-q_{j}|}
\end{array}\right.
\end{eqnarray}
for $i=1,2,\cdots,N.$
\end{prop}

\begin{rem} When $N=1$, the multi-peakon solution
(\ref{mp})-(\ref{pDot}) reduces to the single non-periodic peakon solution given in Proposition~\ref{exsitenceOfPeakon}.
\end{rem}

\begin{rem} For $N=2$, let $P_+(t)=p_1(t)+p_2(t),\, P_-(t)=p_1(t)-p_2(t),\, Q_+(t)=q_1(t)+q_2(t),\, Q_-(t)=q_1(t)-q_2(t)$, then it follows from system (\ref{pDot}) that we have
\begin{eqnarray}
\left\{\begin{array}{l}
\dot{P}_+(t)=\dfrac{k_2}{2}\left[P_+^2+P_-^2+(P_+^2-P_-^2)e^{-|Q_-|}\right], \v\\
\dot{Q}_+(t)=(P_+^2+P_-^2)\left[\dfrac{k_1}{3}+\dfrac{k_2}{2}\left(1+e^{-2|Q_-|}\right)\right]+(k_1+k_2)(P_+^2-P_-^2)e^{-|Q_-|},\v\\
\dot{P}_-(t)=k_2P_+P_-+\dfrac{k_3}{2}(P_-^2-P_+^2)e^{-|Q_-|}, \v\\
\dot{Q}_-(t)=P_+P_-\left[\dfrac{2k_1}{3}+k_2\left(1-e^{-2|Q_-|}\right)\right]+k_3P_-\left(1-e^{-|Q_-|}\right),
\end{array}\right.
\end{eqnarray}
whose explicit solution may be hard to be find, but can be numerically studied.

\end{rem}

\begin{rem}\label{remMultipeakon} When $k_{1}=k_{2}=0$ and $k_{3}=1,$
   system (\ref{pDot}) reduces to the  dynamical system for the CH equation \cite{CamassaHolm1993PRL};
  when $k_{1}=k_{3}=0$ and $k_{2}=1,$   system (\ref{pDot}) becomes the dynamical system for the Novikov  equation \cite{HoneWang2008JPAMT}; when $k_{2}=k_{3}=0$ and $k_{1}=1,$   system (\ref{pDot}) reduces to the
  dynamical system for the mCH equation~\cite{GuiLiuOlverQu2013CMP}, when $k_2=0$,  system (\ref{pDot}) becomes the dynamical system for the mCH-CH equation \cite{xia2018};  when $k_3=0$,  system (\ref{pDot}) becomes the dynamical system for the mCH-Novikov equation~\cite{MiLiuHuangGuo2020JDE}.   System (\ref{pDot}) was also given in \cite{anco2019} by using the distribution theory~\cite{gelfand1964}.
\end{rem}

Proposition~\ref{exsitenceOfPeakon} provides the
explicit peakon formula for equation (\ref{gmCH-Novikov})
on the line. Our next Proposition gives its explicit peakon formula
in the periodic case.

\begin{prop}\label{exsitenceOfPeriodicPeakon} For any $a\neq 0,$ the periodic function of the form
\begin{eqnarray}\label{singlePeriodicPeakon}
u_{c}(t, x)=a \cosh(\xi),\quad \xi=\xi(x,t)=\frac{1}{2}-(x-ct)+\lfloor x-ct\rfloor,
\end{eqnarray}
where the notation $ \lfloor\cdot\rfloor$ denotes the floor function or the greatest integer function  and
\begin{eqnarray}\label{c-con}
c=\left[\frac{k_1}{3}+\cosh^{2}(1/2)\left(\frac{2k_1}{3}+k_2\right)\right]a^{2}+\cosh(1/2) k_3a,
\end{eqnarray}
is a global weak periodic peakon solution to Eq.~(\ref{gmCH-Novikov})
in the sense of Definition \ref{defOfWeakSolution} of the periodic case.
\end{prop}

\begin{rem}
The definition of weak solution for the periodic case
is similar to Definition \ref{defOfWeakSolution},
while $\mathbb{R}$ and $p(x)$ should be replaced by  $\mathbb{S}=[0, 1)$ and $G(x)=\frac{1}{2}{\rm csch}(1/2)\cosh\left(x-1/2\right)$, respectively, there.
\end{rem}

\begin{rem} \label{p-peakon-r2}
We also have a similar remark about the periodic peakon solutions as {\bf Remark}~\ref{peakon-r1}.
In particular, one can reduce the periodic peakon solution given by Eqs.~(\ref{singlePeriodicPeakon})-(\ref{c-con}) to those of the CH equation for $k_1=k_2=0,\,k_3=1$~\cite{QuLiuLiu2013CMP,GrayHimon2013}, where $a={\rm sech}(1/2)c$; the mCH equation for $k_2=k_3=0,\,k_1=1$~\cite{QuLiuLiu2013CMP}, where $a=\pm\sqrt{3c/(2+\cosh(1))}$; the Novikov equation for $k_1=k_3=0,\,k_2=1$~\cite{Grayshan2013,GrayHimon2013}, where $a=\pm {\rm sech}(1/2)\sqrt{c}$; the mCH-CH equation for $k_2=0,\,k_1k_3\not=0$~\cite{wangyan2020}, where $a=[-3\cosh(1/2)k_3\pm\sqrt{
12ck_1+3\cosh^2(1/2)(3k_3^2+8ck_1)}]/[(4+2\cosh(1))k_1]$; the mCH-Novikov equation for $k_1k_2\not=0,\,k_3=0$~\cite{wangyan2020}, where $a=\pm \sqrt{3c/[k_1+\cosh^2(1/2))(2k_1+3k_2)]}$; or the Novikov-CH equation for $k_1=0,\, k_2k_3\not=0$, where $a={\rm sech}(1/2)(-k_3\pm\sqrt{k_3^2+4ck_2})/(2k_2)$.
\end{rem}

\begin{rem} \label{p-peakon-r3} For the given real velocity $c$
 and for $k_1k_2k_3\not=0$,
 one finds from Eq.~(\ref{c-con}) that
\begin{eqnarray}
a=\left\{\begin{array}{ll}
 \dfrac{-3\cosh(1/2)k_3 \pm \sqrt{9\cosh^2(1/2)k_3^2+12cf(k_1,k_2)}}
  {2f(k_1,k_2)}, &  {\rm as}\,\,\, f(k_1,k_2)\not=0, \\
 \dfrac{{\rm sech}(1/2)c}{k_3}, & {\rm as}\,\,\, f(k_1,k_2)=0,\, k_1k_2k_3\not=0,
 \end{array}
\right.
\end{eqnarray}
where $f(k_1,k_2)=[2+\cosh(1)]k_1+3\cosh^2(1/2)k_2$.
Therefore, (i) if $f(k_1,k_2)\not=0,\, k_1k_2k_3\not=0$ and $9\cosh^2(1/2)k_3^2+12cf(k_1,k_2)<0$, then $a$ is complex so that the solution given by Eqs.~(\ref{singlePeakon})-(\ref{solu-c}) is a complex peakon solution; (ii) If $f(k_1,k_2)=0,\, k_1k_2k_3\not=0$, then the periodic peakon solution of Eq.~(\ref{gmCH-Novikov}) is of the form $u_a(t,x)=\frac{{\rm sech}(1/2)c}{k_3}\cosh(\xi)$, which is independent of $k_1,\, k_2$.
\end{rem}

Similarly, the periodic multi-peakon solutions of Eq.~(\ref{gmCH-Novikov}) in the sense of Definition \ref{defOfWeakSolution} for the periodic case are presented as follows.

\begin{prop}\label{mppeakon} The periodic multi-peakon solutions of  Eq.~(\ref{gmCH-Novikov}) in the sense of Definition \ref{defOfWeakSolution} for the periodic case are given by
\begin{eqnarray}\label{mpp}
u_{pm}(t, x)=\sum_{i=1}^{N} p_{i}(t) \cosh\left(\frac{1}{2}-(x-q_i(t))+\lfloor x-q_i(t)\rfloor\right),
\end{eqnarray}
where where the notation $ \lfloor\cdot\rfloor$ denotes the floor function or the greatest integer function, the time-dependent position functions $q_i(t)$ (without loss of generality, the condition, $0<q_1(t)<q_2(t)<\cdots<q_N(t)<1$, is required) and amplitudes $p_i(t)$ satisfy the dynamical system ($\xi_{jm}:=q_j-q_m$)
\begin{align}\label{multi-pp}
\dot{p}_{m}=
&\frac{k_2p_m}{2\operatorname{sinh}(1/2)}\bigg\{\frac{1}{2}\sum_{j\neq m}p_j\operatorname{sgn}(\xi_{jm})
\Big[p_{m}\Big(\operatorname{cosh}(1/2+|\xi_{jm}|)-\operatorname{cosh}(3/2-|\xi_{jm}|)\Big) \no\\
&+p_j\!\Big(\!\operatorname{cosh}(1/2\!-\!2|\xi_{jm}|)\!-\!\operatorname{cosh}(3/2\!-\!2|\xi_{jm}|)\!\Big)\Big]
\!+\!2\sinh(1/2)\!\!\!\sum_{j<m<k}\!\!p_{j}p_{k}\sinh(\xi_{jm}\!+\!\xi_{km})\no\\
&
+\Big(\!\!\sum_{m<j<k}\!\!-\!\!\sum_{j<k<m}\!\!\Big)p_{j}p_{k}\Big(\operatorname{cosh}(1/2\!-\!|\xi_{jm}|\!-\!|\xi_{km}|)
\!-\!\operatorname{cosh}(3/2\!-\!|\xi_{jm}|\!-\!|\xi_{km}|)\Big)\bigg\}\no\\
&+\frac{k_3p_m}{2\operatorname{sinh}(1/2)}
\sum_{j=1}^{N}p_j\operatorname{sgn}(\xi_{jm})\Big(\operatorname{cosh}(|\xi_{jm}|)-\operatorname{cosh}(1-|\xi_{jm}|)\Big), \no\\
\dot{q}_{m}=&\Big[\Big(\frac{2k_1}{3}\!+\!k_2\Big)\sinh^2(1/2)\!+\!\frac{k_2}{2}\!\Big]p_m^2\!+\!
\frac{(k_1\!+\!k_2)p_m}{2\operatorname{sinh}(1/2)}\!\sum_{j\neq m}  \!p_{j}\Big(\sinh(1/2\!+\!|\xi_{jm}|)\!+\!\operatorname{sinh}(3/2\!-\!|\xi_{jm}|)\Big) \no\\
& +\frac{2k_1+k_2}{2\operatorname{sinh}(1/2)}
  \bigg\{\sum_{j<m<k}p_{j}p_{k}\Big(\operatorname{sinh}(3/2+\xi_{jm}-\xi_{km})-\operatorname{sinh}(1/2+\xi_{jm}-\xi_{km})\Big)\no\\
&\qquad +\sinh(1/2)\Big[\sum_{i=1}^Np_{i}^{2}+2\Big(\!\!\sum_{m<j<k}\!\!+\!\!\sum_{j<k<m}\!\!\Big)p_{j}p_{k}\cosh(\xi_{jm}-\xi_{km})\Big]\bigg\}\no\\
& +\frac{k_2}{2\operatorname{sinh}(1/2)}\bigg\{\!\!\sum_{j\neq
m}\!\!\frac{p_{j}^{2}}{2}\!\Big(\!\operatorname{sinh}(3/2\!-\!2|\xi_{jm}|)\!-\!
\operatorname{sinh}(1/2\!-\!2|\xi_{jm}|)\Big)\!+\!2\sinh(1/2)\!\!\!\!
\sum_{j<m<k}\!\!\!\!\cosh(\xi_{jm}\!+\!\xi_{km})\no\\
&\qquad +\Big(\!\!\sum_{m<j<k}\!\!+\!\!\sum_{j<k<m}\!\!\Big)p_{j}p_{k}\Big(\operatorname{sinh}(3/2-|\xi_{km}|-|\xi_{jm}|)
-\operatorname{sinh}(1/2-|\xi_{km}|-|\xi_{jm}|)\Big)\bigg\}\no\\
&+\frac{k_3}{2\operatorname{sinh}(1/2)}\sum_{j=1}^{N}p_{j}\Big(
\operatorname{sinh}(|\xi_{jm}|)+\operatorname{sinh}(1-|\xi_{jm}|)\Big)
\end{align}
for $m=1,2,\cdots,N.$
\end{prop}

\begin{rem} When $N=1$, the periodic multi-peakon solution
(\ref{mpp}) with (\ref{multi-pp}) reduces to the single periodic peakon solution given in Proposition~\ref{exsitenceOfPeriodicPeakon}.
\end{rem}



\begin{rem}
We have investigated some properties of the Cauchy
problem (\ref{gmCH-Novikov}),
including the local well-posedness of its strong solution,
the H\"{o}lder continuity of its data-to-solution map,
the blow-up criterion, the precise blow-up quantity
and the wave-breaking phenomenon and so on.
Some other properties associated with  the Cauchy
problem (\ref{gmCH-Novikov})
can also be concerned with, for instance, the analyticity
the persistence properties of its solution,
which can be finished using similar methods as that
of \cite{HimonasMisiolek2003MA,HimonasMisiolekPonceZhou2007CMP,MiLiuHuangGuo2020JDE,NiZhou2011JDE}.
However, we will not consider these problems in the current paper
in detail.
\end{rem}

The rest of this paper is arranged as follows. The proof of Theorem \ref{WaveBreakingThm1} will be presented in Section 2, where the uniqueness and continuous dependence parts
will be proved first and then the existence part. Section 3 will deal with the H\"{o}lder continuity of the data-to-solution map and prove  Theorem \ref{HolderContinuityThm}. In Section 4, we will invoke the Moser-type estimates in Sobolev spaces to
prove the blow-up criterion, i.e., Theorem \ref{WaveBreakingThm2}. A global conservative property will be displayed
in Section 5, where the precise blow-up quantity will also be established. The uniform upper bound of the blow-up quantity (\ref{MuniformUpperBound}) stated in Theorem \ref{finalThm} will be proved in Section 6.  Section 7 will
provide the proof of the wave-breaking Theorem \ref{WaveBreakingThm3} by only considering a portion of the blow-up quantity.
 Section 8 will give the  proof of Theorem \ref{WaveBreakingThm4} by inspecting the whole blow-up quantity.
The non-periodic peakon solution stated in Proposition \ref{exsitenceOfPeakon} and  multi-peakon solutions of this equation in Proposition~\ref{multipeakon} will be considered in Sections 9 and 10, respectively. The periodic peakon solution stated
in Proposition \ref{exsitenceOfPeriodicPeakon} and multi-peakon in Proposition~\ref{mppeakon} will be considered in Sections
11 and 12, respectively. In Appendixes A and B, some useful facts and Lemmas are given for the proofs
of our main results, including the properties of Besov spaces and some estimates about the transport equation theory.



\section{Proof of Theorem \ref{WaveBreakingThm1}}

Since the local well-posedness for the  Cauchy problem  (\ref{gmCH-Novikov}) will be established in  Besov-type spaces, we firstly recall some basic properties about the Littlewood-Paley theory (see {\bf Appendix A} and \cite{BahouriCheminDanchin2011,
Chemin2004CRM} for more details) and some useful lemmas  of the transport equation  theory (see {\bf Appendix B} and \cite{BahouriCheminDanchin2011,Danchin2001DIE} for more details).

For the proof of Theorem \ref{WaveBreakingThm1}, we first prove the uniqueness and the continuous dependent
on the initial data parts. Suppose $u,\, v\in B_{p, r}^{s}$ are two solutions of Eq.~(\ref{gmchn}), the equivalent form of Eq.~(\ref{gmCH-Novikov}).
Let $w=u-v.$ Then $w$ satisfies the following transport
equation
\begin{eqnarray}\label{difference}
&&w_{t}
+\left[(k_{1}+k_2)v^{2}-\frac{k_1}{3}(u_{x}^{2}+u_{x}v_{x}+v_{x}^{2})+k_{3}v\right]w_{x}\nonumber\\
&&=-[(k_{1}+k_2)(u+v)+k_{3}]wu_{x}\nonumber\\
&&\quad -\left(\frac{k_1}{3}+\frac{k_2}{2}\right)p*[w_{x}(u_{x}^{2}+u_{x}v_{x}+v_{x}^{2})]
-\left(\frac{2k_1}{3}+k_2\right)p_{x}*[w(u^{2}+uv+v^{2})]\nonumber\\
&&\quad -\left(k_1+\frac{3k_2}{2}\right)p_{x}*[wu_{x}(u+v)+v^{2}w_{x}]
-k_{3}p_{x}*\left[w(u+v)+\frac{1}{2}w_{x}(u_{x}+v_{x})\right] \nonumber\\
&& \equiv f(w,u,v)
\end{eqnarray}
with  the initial condition $w(0,x)=w_0(x)=u_0(x)-v_0(x)$.

 Applying (\ref{priori1}) in Lemma \ref{22lem1} to Eq.~(\ref{difference}) yields
\begin{eqnarray}\label{equality1}
\|w\|_{B^{s-1}_{p,r}}
&\leq& C\int_{0}^{t}\left\|(k_{1}+k_2)v^{2}-\frac{k_1}{3}(u_{x}^{2}+u_{x}v_{x}+v_{x}^{2})+k_{3}v\right\|_{B^{s-1}_{p,r}}\|w(\tau)\|_{B^{s-1}_{p,r}}\mbox{d}\tau\nonumber\\
&&+\int_{0}^{t}\|f(w,u,v)\|_{B^{s-1}_{p,r}}\mbox{d}\tau
+\|w(0)\|_{B^{s-1}_{p,r}}.
\end{eqnarray}
Thanks to the product law in the Besov spaces
stated in Lemma \ref{21lem3}
and the embedding relation in Lemma \ref{21lem1},
one obtains
\begin{eqnarray}\label{inequalityStart}
\begin{array}{ll}
\left\|(k_{1}+k_2)v^{2}-\dfrac{k_1}{3}(u_{x}^{2}+u_{x}v_{x}+v_{x}^{2})+k_{3}v\right\|_{B^{s-1}_{p,r}}
&\leq C(\|v\|_{B^{s}_{p,r}}^{2}+\|u\|_{B^{s}_{p,r}}^{2}+\|v\|_{B^{s}_{p,r}}) \v\\
&\leq C(\|v\|_{B^{s}_{p,r}}^{2}+\|u\|_{B^{s}_{p,r}}^{2}+1).
\end{array}
\end{eqnarray}

A similar procedure yields
\begin{eqnarray}
\begin{array}{ll}
\|-(k_{1}+k_{2})wu_{x}(u+v)\|_{B^{s-1}_{p,r}} &\leq C\|w\|_{B^{s-1}_{p,r}}\|u\|_{B^{s-1}_{p,r}}\|u_{x}\|_{B^{s-1}_{p,r}}
+C\|w\|_{B^{s-1}_{p,r}}\|v\|_{B^{s-1}_{p,r}}\|u_{x}\|_{B^{s-1}_{p,r}}\v\\
&\leq C\|w\|_{B^{s-1}_{p,r}}(\|u\|_{B^{s}_{p,r}}^{2}+\|v\|_{B^{s}_{p,r}}^{2})
\end{array}\end{eqnarray}
and
\begin{eqnarray}
\|-k_{3}wu_{x}\|_{B^{s-1}_{p,r}}
\leq C\|w\|_{B^{s-1}_{p,r}}\|u\|_{B^{s}_{p,r}}
\leq C\|w\|_{B^{s-1}_{p,r}}(\|u\|_{B^{s}_{p,r}}^{2}+1).
\end{eqnarray}

Using the Moser-type estimates in Lemma \ref{21lem4},
one finds
\begin{eqnarray}
&&\|-(k_1/3+k_2/2)p*[w_{x}(u_{x}^{2}
+u_{x}v_{x}+v_{x}^{2})]\|_{B^{s-1}_{p,r}}\leq C\|w_{x}(u_{x}^{2}+u_{x}v_{x}
+v_{x}^{2})\|_{B^{s-3}_{p,r}}\nonumber\\
&&\leq C\|w_{x}\|_{B^{s-3}_{p,r}}
\|u_{x}^{2}+u_{x}v_{x}
+v_{x}^{2}\|_{B^{s-2}_{p,r}}\leq C\|w\|_{B^{s-1}_{p,r}}
(\|u\|_{B^{s}_{p,r}}^{2}+\|v\|_{B^{s}_{p,r}}^{2}),
\end{eqnarray}
\begin{eqnarray}
&&\|-(2k_1/3+k_2)p_{x}*[w(u^{2}
+uv+v^{2})]\|_{B^{s-1}_{p,r}}\leq C\|[w(u^{2}+uv+v^{2})]\|_{B^{s-2}_{p,r}}\nonumber\\
&&\leq C\|w\|_{B^{s-2}_{p,r}}
\|u^{2}+uv+v^{2}\|_{B^{s-1}_{p,r}}\leq C\|w\|_{B^{s-1}_{p,r}}
(\|u\|_{B^{s}_{p,r}}^{2}+\|v\|_{B^{s}_{p,r}}^{2}),
\end{eqnarray}
\begin{eqnarray}
&&\|-(k_{1}+3k_2/2)p_{x}*[wu_{x}(u+v)+v^{2}w_{x}]\|_{B^{s-1}_{p,r}}\leq C\|wuu_{x}+wvu_{x}
+v^{2}w_{x}\|_{B^{s-2}_{p,r}}\nonumber\\
&&\leq C(\|w\|_{B^{s-2}_{p,r}}+
\|w_{x}\|_{B^{s-2}_{p,r}})
(\|u^{2}\|_{B^{s-1}_{p,r}}
+\|vu_{x}\|_{B^{s-1}_{p,r}}
+\|v^{2}\|_{B^{s-1}_{p,r}})\nonumber\\
&&\leq C\|w\|_{B^{s-1}_{p,r}}
(\|u\|_{B^{s}_{p,r}}^{2}+\|v\|_{B^{s}_{p,r}}^{2}),
\end{eqnarray}
\begin{eqnarray}
\begin{array}{rl}
\|-k_{3}p_{x}*[w(u+v)]\|_{B^{s-1}_{p,r}}
\leq C\|w\|_{B^{s-1}_{p,r}}(\|u\|_{B^{s}_{p,r}}+\|v\|_{B^{s}_{p,r}})\leq C\|w\|_{B^{s-1}_{p,r}}(\|u\|_{B^{s}_{p,r}}^{2}+\|v\|_{B^{s}_{p,r}}^{2}+1)
\end{array}
\end{eqnarray}
and
\begin{eqnarray}\label{inequalityEnd}
\|-k_3/2p_{x}*[w_{x}(u_{x}+v_{x})]\|_{B^{s-1}_{p,r}}
\leq C\|w\|_{B^{s-1}_{p,r}}(\|u\|_{B^{s}_{p,r}}^{2}+\|v\|_{B^{s}_{p,r}}^{2}+1).
\end{eqnarray}

Substituting Eqs.~(\ref{inequalityStart})-(\ref{inequalityEnd}) into Eq.~(\ref{equality1}) leads to
\begin{eqnarray*}
\|w\|_{B^{s-1}_{p,r}}
\leq \|w(0)\|_{B^{s-1}_{p,r}}
+C\int_{0}^{t}
(\|u\|_{B^{s}_{p,r}}^{2}+\|v\|_{B^{s}_{p,r}}^{2}+1)
\|w\|_{B^{s-1}_{p,r}}\mbox{d}\tau.
\end{eqnarray*}
The Gronwall's inequality then yields
\begin{eqnarray*}
\|w\|_{B^{s-1}_{p,r}}
\leq \|w(0)\|_{B^{s-1}_{p,r}}
\exp\left\{
C\int_{0}^{t}
(\|u\|_{B^{s}_{p,r}}^{2}+\|v\|_{B^{s}_{p,r}}^{2}+1)
\mbox{d}\tau
\right\},
\end{eqnarray*}
which is suffice to  illustrate the
uniqueness and the continuous dependence
parts of Theorem \ref{WaveBreakingThm1}.

Motivated by the proofs of local existence  of the CH equation\cite{Danchin2001DIE} and mCH-Novikov equation\cite{MiLiuHuangGuo2020JDE}, we next prove the local existence part of Theorem \ref{WaveBreakingThm1}. Let us start with the construction of approximate solutions for  Eq.~(\ref{gmCH-Novikov}) in the spirit
of Friedrichs. Let $m^{(0)}=0$, and $m^{(n+1)}(n=0,1, \cdots.)$ being the solutions of
the following  linear transport equations
\begin{eqnarray}\label{linearTrans}
\begin{cases}
\partial_{t}m^{(n+1)}
+\left[(k_{1}+k_2)(u^{(n)})^{2}-k_{1}(u_{x}^{(n)})^{2}+k_{3}u^{(n)}\right]\partial_{x}m^{(n+1)} \v\\
\qquad =-u_{x}^{(n)}m^{(n)}\left(2k_{1}m^{(n)}+3k_{2}u^{(n)}+2k_{3}\right),\v\\
u^{(n+1)}(0,x)=u^{(n+1)}_0(x)=S_{n+1}u_{0}(x).
\end{cases}
\end{eqnarray}
Suppose $m^{(n)} \in L^{\infty}(0, T ; B_{p, r}^{s-2}), s-2>\max \{\frac{1}{p}, \frac{1}{2}\}.$
The condition on $s$ ensures that $B_{p, r}^{s-2}$ is an algebra.
So the right hand side of
Eq.~(\ref{linearTrans}) is in
$L^{\infty}\left(0, T ; B_{p, r}^{s-2}\right).$
Hence, by Lemma \ref{22lem4} and the high regularity of $u,$
Eq.~(\ref{linearTrans})
has a global solution $m^{(n+1)}\in E_{p, r}^{s-2}$ for all positive $T$.

Applying (\ref{priori2}) in Lemma \ref{22lem1} to Eq.~(\ref{linearTrans}) leads to
\begin{eqnarray}\label{induction1}
&&\|m^{(n+1)}\|_{B^{s-2}_{p,r}}\leq \|S_{n+1}u_{0}\|_{B^{s}_{p,r}}\exp\left[C\int_{0}^{t}\|(k_{1}+k_2)(u^{(n)})^{2}-k_{1}(u_{x}^{(n)})^{2}
+k_{3}u^{(n)}\|_{B^{s-1}_{p,r}}(\tau)\mbox{d}\tau\right]\nonumber\\
&&\qquad\qquad\qquad\qquad+C\int_{0}^{t}
\exp\left[C\int_{\tau}^{t}\|(k_{1}+k_2)(u^{(n)})^{2}-k_{1}(u_{x}^{(n)})^{2}
+k_{3}u^{(n)}\|_{B^{s-1}_{p,r}}(\tau^{\prime})
\mbox{d}\tau^{\prime}\right]\nonumber\\
&&\qquad\qquad\qquad\qquad\times \|-2k_{1}u_{x}^{(n)}(m^{(n)})^{2}
-3k_{2}u^{(n)}u_{x}^{(n)}m^{(n)}
-2k_{3}u_{x}^{(n)}m^{(n)}\|_{B^{s-2}_{p,r}}(\tau)\mbox{d}\tau
\end{eqnarray}
for  $n=0,1,2,\cdots.$

 Lemma \ref{21lem3}  yields
\begin{eqnarray}\label{induction2}
\|(k_{1}+k_2)(u^{(n)})^{2}-k_{1}(u_{x}^{(n)})^{2}+k_{3}u^{(n)}\|_{B^{s-1}_{p,r}}\leq C(\|u^{(n)}\|_{B^{s}_{p,r}}^{2}+\|u^{(n)}\|_{B^{s}_{p,r}})\leq C(\|u^{(n)}\|_{B^{s}_{p,r}}^{2}+1),
\end{eqnarray}
\begin{eqnarray}\label{induction3}
&&\|-2k_{1}u_{x}^{(n)}(m^{(n)})^{2}
-3k_{2}u^{(n)}u_{x}^{(n)}m^{(n)}
-2k_{3}u_{x}^{(n)}m^{(n)}\|_{B^{s-2}_{p,r}}\leq C(\|u^{(n)}\|_{B^{s}_{p,r}}^{3}+\|u^{(n)}\|_{B^{s}_{p,r}}^{2})\qquad \nonumber\\
&&\leq C(\|u^{(n)}\|_{B^{s}_{p,r}}^{3}+\|u^{(n)}\|_{B^{s}_{p,r}}).
\end{eqnarray}

Substituting (\ref{induction2})-(\ref{induction3}) into (\ref{induction1}), one finds
\begin{eqnarray}\label{induction4}
&&\|u^{(n+1)}\|_{B^{s}_{p,r}}=\|m^{(n+1)}\|_{B^{s-2}_{p,r}}\nonumber\\
&&\qquad\leq \|S_{n+1}u_{0}\|_{B^{s}_{p,r}}\exp\left[C\int_{0}^{t}(\|u^{(n)}\|_{B^{s}_{p,r}}^{2}+1)(\tau)\mbox{d}\tau\right] \nonumber\\
&&\qquad\quad +C\int_{0}^{t}\exp\left[C\int_{\tau}^{t}(\|u^{(n)}\|_{B^{s}_{p,r}}^{2}
+1)(\tau^{\prime})
\mbox{d}\tau^{\prime}\right]\left(\|u^{(n)}\|_{B^{s}_{p,r}}^{3}+\|u^{(n)}\|_{B^{s}_{p,r}}\right)(\tau)\mbox{d}\tau.
\end{eqnarray}

To obtain the uniform upper bound of the solution sequence $\{u^{(n)}\},$ we assume
$\|u^{(n)}\|_{B^{s}_{p,r}}\leq a(t)$ with
 $a(t)$ some unknown function of $t$.
This combined with (\ref{induction4}) yields
\begin{eqnarray}\label{induction5}
&&\|u^{(n+1)}\|_{B^{s}_{p,r}}\leq \|u_{0}\|_{B^{s}_{p,r}}\exp\left[C\int_{0}^{t}(a^2(\tau)+1)\mbox{d}\tau\right]\nonumber\\
&&\qquad\qquad\qquad\quad + C\int_{0}^{t}\exp\left[C\int_{\tau}^{t}\left((a^2(\tau^{\prime}))+1\right)\mbox{d}\tau^{\prime}\right]
\left(a^3(\tau)+a(\tau)\right)\mbox{d}\tau.
\end{eqnarray}
Let the right hand side  of (\ref{induction5}) be equal to $a(t).$
Then $a(t)$ satisfies the following ordinary differential equation
\begin{eqnarray}\label{aeq}
\begin{cases}
\dot{a}(t)=2C[a^3(t)+a(t)],\\
a(0)=\|u_{0}\|_{B^{s}_{p,r}}.
\end{cases}
\end{eqnarray}

Solving  (\ref{aeq}), we deduce that a time $T$ can be chosen to satisfy
\begin{eqnarray*}
0<T< \frac{1}{4C}\ln\bigg(1+\frac{1}{a^{2}(0)}\bigg)
\end{eqnarray*}
so that the solution sequence $\{u^{(n)}\}$ of the linear
transport equation (\ref{linearTrans})
has the uniform bound
\begin{eqnarray}\label{at}
a(t)=\frac{a(0) e^{2Ct}}{\sqrt{a^{2}(0)(1- e^{4Ct})+1}}
\end{eqnarray}
 for $t\in [0,T].$
 Therefore, $\{u^{(n)}\}$ is uniformly bounded
 in $L^{\infty}([0,T]; B^{s}_{p, r}).$


In the following,
we will prove that
the solution sequence $\{u^{(n)}\}$
is  Cauchy in $L^{\infty}([0,T]; B^{s-1}_{p, r}),$
or equivalently $\{m^{(n)}\}$
is  Cauchy in $L^{\infty}([0,T]; B^{s-3}_{p, r})$.
Direct computation leads to
\begin{eqnarray}\label{cauchy1}
&&\partial_{t}(m^{(n+1+l)}-m^{(n+1)})+[(k_{1}+k_2)(u^{(n+l)})^{2}-k_{1}(u_{x}^{(n+l)})^{2}+k_{3}u^{(n+l)}]
\partial_{x}(m^{(n+1+l)}-m^{(n+1)})\nonumber\\
&&=\bigg[(k_{1}+k_2)((u^{(n)})^{2}-(u^{(n+l)})^{2})
-k_{1}((u_{x}^{(n)})^{2}-(u_{x}^{(n+l)})^{2})\nonumber\\
&&\quad+k_{3}(u^{(n)}-u^{(n+l)})\bigg]\partial_{x}(m^{(n+1+l)}-m^{(n+1)})\nonumber\\
&&\quad -2k_{1}[u_{x}^{(n+l)}(m^{(n+l)})^{2}
-u_{x}^{(n)}(m^{(n)})^{2}]-3k_{2}[u^{(n+l)}u_{x}^{(n+l)}m^{(n+l)}
-u^{(n)}u_{x}^{(n)}m^{(n)}]\nonumber\\
&&\quad -2k_{3}[u_{x}^{(n+l)}m^{(n+l)}-u_{x}^{(n)}m^{(n)}]
\equiv g, \qquad n,l=0,1,\cdots.
\end{eqnarray}

The right hand side of (\ref{cauchy1})
can be rewritten as
\begin{eqnarray*}
&&g=\bigg[(u^{(n)}-u^{(n+l)})[(k_{1}+k_2)(u^{(n)}+u^{(n+l)})+k_3]
-k_{1}(u_{x}^{(n)}-u_{x}^{(n+l)})(u_{x}^{(n)}+u_{x}^{(n+l)})\bigg]\partial_{x}m^{(n+1)}\nonumber\\
&&\quad\quad-2k_{1}[(u_{x}^{(n+l)}-u_{x}^{(n)})(m^{(n+1)})^{2}
+u_{x}^{(n)}(m^{(n+l)}-m^{(n)})(m^{(n+l)}+m^{(n)})]\nonumber\\
&&\quad\quad-3k_{2}\bigg[(u^{(n+l)}-u^{(n)})u_{x}^{(n+l)}m^{(n+1)}
+u^{(n)}(u_{x}^{(n+l)}-u_{x}^{(n)})m^{(n+1)}+u^{(n)}u_{x}^{(n)}(m^{(n+l)}-m^{(n)})\bigg]\nonumber\\
&&\quad\quad-2k_{3}[(u_{x}^{(n+l)}-u_{x}^{(n)})m^{(n+1)}
+u_{x}^{(n)}(m^{(n+l)}-m^{(n)})].
\end{eqnarray*}
According to Lemma \ref{22lem1},
one finds
\begin{eqnarray}\label{cauchy2}
&&\|m^{(n+1+l)}-m^{(n+1)}\|_{B^{s-3}_{p,r}}\nonumber\\
&&\leq \exp\left\{C\int_{0}^{t}
\|(k_{1}+k_2)(u^{(n+l)})^{2}-k_{1}(u_{x}^{(n+l)})^{2}+k_{3}u^{(n+l)}\|_{B^{s-3}_{p,r}}(\tau)\mbox{d}\tau
\right\}\bigg[\|m_{0}^{(n+1+l)}-m_{0}^{(n+1)}\|_{B^{s-3}_{p,r}} \,\, \nonumber\\
&&\quad+\int_{0}^{t}\exp\left(-C\int_{0}^{\tau}\|(k_{1}+k_2)(u^{(n+l)})^{2}-k_{1}(u_{x}^{(n+l)})^{2}+k_{3}u^{(n+l)}
\|_{B^{s-3}_{p,r}}(\tau^{\prime})\mbox{d}\tau^{\prime}\right)\|g\|_{B^{s-3}_{p,r}}\mbox{d}\tau\bigg].\,\,
\end{eqnarray}

Invoking Lemma \ref{21lem3}, we obtain
\begin{eqnarray}\label{cauchy3}
&&\bigg\|\left\{[(k_{1}+k_2)(u^{(n)}+u^{(n+l)})+k_3](u^{(n)}-u^{(n+l)})
-k_{1}(u_{x}^{(n)}-u_{x}^{(n+l)})(u_{x}^{(n)}+u_{x}^{(n+l)})\right\}\partial_{x}m^{(n+1)}\bigg\|_{B^{s-3}_{p,r}}\nonumber\\
&&\quad \leq C\|\partial_{x}m^{(n+1)}\|_{B^{s-3}_{p,r}}
\|u^{(n)}-u^{(n+l)}\|_{B^{s-1}_{p,r}}
\|u^{(n)}+u^{(n+l)}\|_{B^{s-1}_{p,r}}\nonumber\\
&&\quad \leq C\|m^{(n)}-m^{(n+l)}\|_{B^{s-3}_{p,r}}
(\|u^{(n)}\|_{B^{s}_{p,r}}^{2}
+\|u^{(n+l)}\|_{B^{s}_{p,r}}^{2}
+\|u^{(n+1)}\|_{B^{s}_{p,r}}^{2}),
\end{eqnarray}
\begin{eqnarray}\label{cauchy4}
&&\|-2k_{1}[(u_{x}^{(n+l)}-u_{x}^{(n)})(m^{(n+1)})^{2}
+u_{x}^{(n)}(m^{(n+l)}-m^{(n)})(m^{(n+l)}+m^{(n)})
]\|_{B^{s-3}_{p,r}}\nonumber\\
&&\quad \leq C\|m^{(n+l)}\|_{B^{s-2}_{p,r}}^{2}
\|u_{x}^{(n)}-u_{x}^{(n+l)}\|_{B^{s-3}_{p,r}}+\|m^{(n)}-m^{(n+l)}\|_{B^{s-3}_{p,r}}\|u_{x}^{(n)}\|_{B^{s-2}_{p,r}}
\|m^{(n)}+m^{(n+l)}\|_{B^{s-2}_{p,r}}\nonumber\\
&&\quad \leq C\|m^{(n)}-m^{(n+l)}\|_{B^{s-3}_{p,r}}
(\|u^{(n)}\|_{B^{s}_{p,r}}^{2}
+\|u^{(n+l)}\|_{B^{s}_{p,r}}^{2}),
\end{eqnarray}
\begin{eqnarray}\label{cauchy5}
&&\bigg\|-3k_{2}[(u^{(n+l)}-u^{(n)})u_{x}^{(n+l)}m^{(n+1)}
+u^{(n)}(u_{x}^{(n+l)}-u_{x}^{(n)})m^{(n+1)}+u^{(n)}u_{x}^{(n)}(m^{(n+l)}-m^{(n)})]\bigg\|_{B^{s-3}_{p,r}}\nonumber\\
&&\quad \leq C\|u^{(n)}-u^{(n+l)}\|_{B^{s-3}_{p,r}}
\|u_{x}^{(n+l)}\|_{B^{s-2}_{p,r}}
\|m^{(n+l)}\|_{B^{s-2}_{p,r}}+C\|u^{(n)}\|_{B^{s-2}_{p,r}}
\|u_{x}^{(n)}-u_{x}^{(n+l)}\|_{B^{s-3}_{p,r}}\nonumber\\
&&\quad \leq C\|m^{(n)}-m^{(n+l)}\|_{B^{s-3}_{p,r}}
(\|u^{(n)}\|_{B^{s}_{p,r}}^{2}
+\|u^{(n+l)}\|_{B^{s}_{p,r}}^{2}+1)
\end{eqnarray}
and
\begin{eqnarray}\label{cauchy6}
&&\|-2k_{3}[(u_{x}^{(n+l)}-u_{x}^{(n)})m^{(n+1)}
+u_{x}^{(n)}(m^{(n+l)}-m^{(n)})]\|_{B^{s-3}_{p,r}}\nonumber\\
&&\quad \leq C\|m^{(n+l)}\|_{B^{s-3}_{p,r}}
\|u_{x}^{(n)}-u_{x}^{(n+l)}\|_{B^{s-2}_{p,r}}
+C\|m^{(n)}-m^{(n+l)}\|_{B^{s-3}_{p,r}}
\|u_{x}^{(n)}\|_{B^{s-2}_{p,r}}\nonumber\\
&&\quad \leq  C\|m^{(n)}-m^{(n+l)}\|_{B^{s-3}_{p,r}}
(\|u^{(n)}\|_{B^{s}_{p,r}}^{2}
+\|u^{(n+l)}\|_{B^{s}_{p,r}}^{2}+1).
\end{eqnarray}

The substitution of (\ref{cauchy3})-(\ref{cauchy6}) into (\ref{cauchy2}) derives
\begin{eqnarray}
&&\|m^{(n+1+l)}-m^{(n+1)}\|_{B^{s-3}_{p,r}}\nonumber\\
&&\leq \exp\bigg\{C\int_{0}^{t}
\|(k_{1}+k_2)(u^{(n+l)})^{2}-k_{1}(u_{x}^{(n+l)})^{2}+k_{3}u^{(n+l)}
\|_{B^{s-3}_{p,r}}(\tau)d\tau
\bigg\}\bigg[\|m_{0}^{(n+1+l)}-m_{0}^{(n+1)}\|_{B^{s-3}_{p,r}}\nonumber\\
&&\quad+\int_{0}^{t}\exp\left(-C\int_{0}^{\tau}
\|(k_{1}+k_2)(u^{(n+l)})^{2}-k_{1}(u_{x}^{(n+l)})^{2}+k_{3}u^{(n+l)}
\|_{B^{s-3}_{p,r}}(\tau^{\prime})\mbox{d}\tau^{\prime}\right)\nonumber\\
&&\quad\times\|m^{(n)}-m^{(n+l)}\|_{B^{s-3}_{p,r}}
(\|u^{(n+1)}\|_{B^{s}_{p,r}}^{2}
+\|u^{(n)}\|_{B^{s}_{p,r}}^{2}
+\|u^{(n+l)}\|_{B^{s}_{p,r}}^{2}+1)\mbox{d}\tau\bigg].
\end{eqnarray}

Note that
\begin{eqnarray}
\begin{array}{rl}
\|m_{0}^{(n+1+l)}-m_{0}^{(n+1)}\|_{B^{s-3}_{p,r}}
&=\|S_{n+1+l}m_{0}-S_{n+1}m_{0}\|_{B^{s-3}_{p,r}} \v\\
&=\|\sum_{q=n+1}^{n+l}\Delta_{q}m_{0}\|_{B^{s-3}_{p,r}}\leq C2^{-n}\|m_{0}\|_{B^{s-3}_{p,r}}
\end{array}
\end{eqnarray}
and $\{u^{(n)}\}$ is bounded in $L^{\infty}([0,T]; B_{p,r}^{s})$,
we thus deduce
\begin{eqnarray}
\|m^{(n+1+l)}-m^{(n+1)}\|_{B^{s-3}_{p,r}}
\leq C_{T}\left(2^{-n}
+\int_{0}^{t}\|m^{(n+l)}-m^{(n)}\|_{B^{s-3}_{p,r}}\mbox{d}\tau\right).
\end{eqnarray}
Then from induction one concludes
\begin{eqnarray}
\|m^{(n+1+l)}-m^{(n+1)}\|_{L^{\infty}(0,T; B^{s-3}_{p,r})}
\leq \frac{C_{T}}{2^n}\sum_{k=0}^{n}\frac{(2TC_{T})^{k}}{k!}+\frac{(TC_{T})^{n+1}}{(n+1)!}\|m^{(l)}-m^{(0)}\|_{L^{\infty}(0,T; B^{s-3}_{p,r})}.
\end{eqnarray}
Since $\{m^{(n)}\}$ is uniformly bounded
in ${L^{\infty}(0,T; B^{s-3}_{p,r})},$
one can find a new constant $C_{T}^{\prime}$
so that
\begin{equation*}
\|m^{(n+l+1)}-m^{(n+1)}\|_{L^{\infty}(0, T ; B_{p, r}^{s-3})} \leq \frac{C_{T}^{\prime}}{2^{n}}.
\end{equation*}
Accordingly, $\{m^{(n)}\}$ is  Cauchy  in $L^{\infty}(0, T ; B_{p, r}^{s-3})$ and   converges to some limit function $m \in L^{\infty}(0, T ; B_{p, r}^{s-3}).$

We are now in a position to prove the existence of the solution to equation (\ref{gmCH-Novikov}).
We will prove the limit function $m$ satisfies equation
(\ref{gmCH-Novikov}) in the sense of distribution and belongs to
$E^{s}_{p,r}.$

Firstly, the Fatou property for the Besov spaces in Lemma \ref{besovProperty}(iii)
and the uniform boundedness of $\{m^{(n)}\}$ in
$L^{\infty}(0,T; B^{s-2}_{p,r})$ ensure
that $m\in L^{\infty}(0,T; B^{s-2}_{p,r}).$

Secondly, we claim that $\{m^{(n)}\}$
converges to $m$ in $ L^{\infty}(0,T; B^{s^{\prime}}_{p,r})$
for all $s^{\prime} < s-2.$
In fact, this results from the following statement: $\|m_{n}-m\|_{B_{p, r}^{s^{\prime}}} \leq C\|m_{n}-m\|_{B_{p, r}^{s-3}}$
when $s^{\prime} \leq s-3$
and $
\|m_{n}-m\|_{B_{p, r}^{\prime}}
\leq C\|m_{n}-m\|_{B_{p, r}^{s-3}}^{\theta}(\|m_{n}\|_{B_{p, r}^{s-2}}+\|m\|_{B_{p, r}^{s-2}})^{1-\theta}
$
with $\theta=s-2-s^{\prime}$
when $s-3<s^{\prime} \leq s-2$
on account of the interpolation between Besov spaces stated in
Lemma \ref{21lem2}.
This claim enables us to take limit in equation (\ref{linearTrans})
to conclude that the limit function $m$ does satisfy equation (\ref{gmCH-Novikov}).

Note that equation (\ref{gmCH-Novikov}) can be rewritten as a transport equation
\begin{eqnarray}\label{meq}
\partial_{t}m
+[k_{1}(u^{2}-u_{x}^{2})
+k_{2}u^{2}+k_{3}u]\partial_{x}m
=-2k_{1}u_{x}m^{2}
-3k_{2}uu_{x}m
-2k_{3}u_{x}m.
\end{eqnarray}
Since $m\in L^{\infty}(0,T; B^{s-2}_{p,r}),$
it is easy to deduce that the right hand side of the above
equation  also belongs to $L^{\infty}(0,T; B^{s-2}_{p,r})$
in view of the product law in Besov spaces and the
Sobolev embedding.
Consequently, Lemma \ref{22lem4} implies
$m\in C([0, T) ; B_{p, r}^{s-2})$ when $r<\infty$
or
$m\in C_{w}([0, T) ; B_{p, r}^{s-2})$ when $r=\infty.$
On the other hand, from the Moser-type estimates
in Lemma \ref{21lem4},
we infer that $[k_{1}u^{2}-k_{1}u_{x}^{2}
+k_{2}u^{2}+k_{3}u]\partial_{x}m$
is bounded in $L^{\infty}(0,T; B^{s-3}_{p,r}).$
Therefore, one knows
$\partial_{t} m\in C([0, T) ; B_{p, r}^{s-3})$
when $r<\infty$ according to the  high regularity of $u$
and equation (\ref{gmCH-Novikov}) and thus
$m\in E^{s-2}_{p,r}.$ We thus complete the proof of Theorem \ref{WaveBreakingThm1}.

\section{Proof of Theorem \ref{HolderContinuityThm}}

In this section, we will prove Theorem \ref{HolderContinuityThm} in the spirit of \cite{HimonasMantzavinos2014JNS} for the mCH equation.
We first of all prove the Lipschitz continuity, namely,
$\beta=1$ in the region $A_{1}.$ Applying $\partial_{x}$ to the equivalent, nonlocal form (\ref{gmchn}) of Eq.~(\ref{gmCH-Novikov}) and using (\ref{pPlusMinus}) produce
\begin{eqnarray}\label{equivalent-1-gradient}
&&\partial_{t}u_{x}
=k_{1}u_{x}^{2}u_{xx}-\left(k_{1}+\frac{k_2}{2}\right)uu_{x}^{2}-\left(k_{1}+k_{2}\right)u^{2}u_{xx}
+\left(\frac{2k_1}{3}+k_{2}\right)u^{3}-\frac{k_3}{2}u_{x}^{2}-k_{3}uu_{xx}+k_{3}u^{2}\nonumber\\
&&\quad\quad\quad \quad  -p_{x}*\left[\left(\frac{k_1}{3}+\frac{k_2}{2}\right)u_{x}^{3}\right]
-p*\left[\left(\frac{2k_1}{3}+k_{2}\right)u^{3}+\left(k_{1}+\frac{3k_2}{2}\right)uu_{x}^{2}
+k_{3}u^{2}
+\frac{k_3}{2}u_{x}^{2}
\right].
\end{eqnarray}

Let $w=u_{x}$. Then Eqs.~(\ref{gmchn}) and (\ref{equivalent-1-gradient}) lead to
\begin{eqnarray}\label{coupled}
\begin{cases}
 u_{t}=\displaystyle\frac{k_1}{3}w^{3}-(k_{1}+k_{2})u^{2}w-k_{3}uw-F(u,w),\\
w_t=\displaystyle k_{1}w^{2}w_{x}-\left(k_{1}+\frac{k_2}{2}\right)uw^{2}-\left(k_{1}+k_{2}\right)u^{2}w_{x}+\left(\frac{2k_1}{3}+k_{2}\right)u^{3}\\
\quad\quad\quad \displaystyle-\frac{k_3}{2}w^{2}-k_{3}uw_{x}+k_{3}u^{2}-G(u,w), \\
u(0,x)=u_0(x),\quad w(0,x)=\partial_x u_0(x)=w_0(x),
\end{cases}
\end{eqnarray}
where the nonlocal terms $F$ and $G$ are given by
\begin{eqnarray}
&&F(u,w)=p*\left[\left(\frac{k_1}{3}+\frac{k_2}{2}\right)w^{3}\right]+p_{x}*\left[\left(\frac{2k_1}{3}+k_{2}\right)u^{3}
+\left(k_{1}+\frac{3k_2}{2}\right)uw^{2}+k_{3}u^{2}+\frac{k_3}{2}w^{2}\right],\qquad \label{coupledF}\\
&&G(u,w)=p_{x}*\left[\left(\frac{k_1}{3}+\frac{k_2}{2}\right)w^{3}\right]+p*\left[\left(\frac{2k_1}{3}+k_{2}\right)u^{3}
+\left(k_{1}+\frac{3k_2}{2}\right)uw^{2}+k_{3}u^{2}+\frac{k_3}{2}w^{2}\right].\qquad \label{coupledG}
\end{eqnarray}

Using  similar method as in \cite{HimonasMantzavinos2014NA},
one can show that for $u_0(x)\in H^s\, (s>5/2)$ the solution of system
(\ref{coupled}) corresponding to  initial data
$(u_{0}(x),w_{0}(x))$ satisfies $(u,w)\in C([0, T]; H^{s-1})$ and the following size estimates
in the lifespan of the solution
\begin{eqnarray}\label{sizeEstimate}
\|(u,w)\|_{H^{s-1}}
\leq C\|u_{0}\|_{H^{s}},\quad (s>5/2)
\end{eqnarray}
with $C$  a generic constant.

Let $(v,z)\in C([0, T]; H^{s-1})$ be another solution to system
(\ref{coupled}) corresponding to  initial data $(v_{0}(x),z_{0}(x))$
Set $\varphi=u-v,\, \psi=w-z$. Then $\varphi$ and $\psi$ satisfy
\begin{eqnarray}\label{coupledDifference1}
\partial_{t}\varphi=
\left[\frac{k_1}{3}g-(k_{1}+k_{2})u^{2}\right] \psi
-(k_{1}+k_{2})(u+v) z \varphi
-k_{3}u\psi
-k_{3}\varphi z
-F(u, w)+F(v, z)
\end{eqnarray}
and
\begin{eqnarray}\label{coupledDifference2}
&&\partial_{t}\psi=
\left[k_{1}w^{2}-(k_{1}+k_{2})u^{2}-k_{3}u\right]\partial_{x}\psi+\left[k_{1}(w+z) \partial_{x} z-
\left(k_{1}+\frac{k_2}{2}\right)u(w+z)-\frac{k_3}{2}(w+z)\right] \psi \nonumber\\
&&\qquad\quad +\left[\left(\frac{2k_1}{3}+k_{2}\right) f-\left(k_{1}+\frac{k_2}{2}\right)z^{2}-\left(k_{1}+k_{2}\right)(u+v) \partial_{x} z-k_{3}\partial_{x} z+k_{3}(u+v)\right] \varphi\nonumber\\
&&\qquad\quad-G(u, w)+G(v, z)
\end{eqnarray}
respectively, where $f=u^{2}+u v+v^{2}$  and $g=w^{2}+w z+z^{2}.$

Applying the energy method to Eqs.~(\ref{coupledDifference1})-(\ref{coupledDifference2})
gives
\begin{subequations}\label{energy1}
\begin{align}
\frac{1}{2} \frac{\mathrm{d}}{\mathrm{d}t}\|\varphi\|_{H^{r}}^{2}
=&
 \int_{\mathbb{R}} D^{r}\left[\left(\frac{k_1}{3} g-(k_{1}+k_{2})u^{2}\right)\psi\right] \cdot D^{r} \varphi \mathrm{d} x \label{energy1-1}\\
&-\int_{\mathbb{R}} D^{r}[(k_{1}+k_{2})(u+v) z \varphi] \cdot D^{r} \varphi \mathrm{d} x \label{energy1-2}\\
&-\int_{\mathbb{R}} D^{r}[k_{3}u\psi] \cdot D^{r} \varphi \mathrm{d} x \label{energy1-3}\\
&-\int_{\mathbb{R}} D^{r}[k_{3}\varphi z] \cdot D^{r} \varphi \mathrm{d} x \label{energy1-4}\\
&-\int_{\mathbb{R}} D^{r}[F(u, w)-F(v, z)] \cdot D^{r} \varphi \mathrm{d} x\label{energy1-5}
\end{align}
\end{subequations}
and
\begin{subequations}\label{energy2}
\begin{align}
\frac{1}{2} \frac{\mathrm{d}}{\mathrm{d} t}\|\psi\|_{H^{r}}^{2}
=& \int_{\mathbb{R}} D^{r}\left[\left(k_{1}w^{2}-(k_{1}+k_{2})u^{2}-k_{3}u\right) \partial_{x}\psi\right] \cdot D^{r} \psi \mathrm{d} x \label{energy2-1}\\
&+\int_{\mathbb{R}} D^{r}\bigg\{\bigg[k_{1}(w+z) \partial_{x} z-
\left(k_{1}+\frac{k_2}{2}\right)u(w+z)-\frac{k_3}{2}(w+z)
\bigg] \psi\bigg\} \cdot D^{r} \psi
 \mathrm{d}x \label{energy2-2}\\
&+\int_{\mathbb{R}} D^{r}\left\{\left[\frac{2}{3} f-z^{2}-(u+v) \partial_{x} z\right] \varphi\right\} \cdot D^{r} \psi
\mathrm{d} x \label{energy2-3}\\
&+\int_{\mathbb{R}} D^{r}\bigg\{\bigg[\left(\frac{2k_1}{3}+k_{2}\right) f
-\left(k_{1}+\frac{k_2}{2}\right)z^{2}-[\left(k_{1}+k_{2}\right)(u+v)+k_3] \partial_{x} z\nonumber\\
&\qquad\qquad\qquad
+k_{3}(u+v)\bigg] \varphi\bigg\} \cdot D^{r} \psi
\mathrm{d} x\label{energy2-4}\\
&-\int_{\mathbb{R}} D^{r}[G(u, w)-G(v, z)] \cdot D^{r} \psi \mathrm{d} x.\label{energy2-6}
\end{align}
\end{subequations}

We first  estimate the terms in (\ref{energy1}). As $1/2<r\leq s-1$,
one can invoke
the algebra property of $H^{r}$
 and (\ref{sizeEstimate}) to deduce
 \begin{eqnarray}\label{energy1-1Esimate1}
 &&|(\ref{energy1-1})|
 \leq C\left\|\left(\frac{k_1}{3}g
 -(k_{1}+k_{2})u^{2}\right)\psi\right\|_{H^{r}}
 \|\varphi\|_{H^{r}}
 \leq C
\left(\|g\|_{H^{r}}+\|u\|_{H^{r}}^{2}\right)
\|\varphi\|_{H^{r}}\|\psi\|_{H^{r}} \nonumber\\
&&\,\,\,\quad\quad\quad
\leq
C(\|w\|_{H^{s-1}}^{2}+\|w\|_{H^{s-1}}
\|z\|_{H^{s-1}}+\|z\|_{H^{s-1}}^{2}+\|u\|_{H^{s-1}}^{2})
\|\varphi\|_{H^{r}}\|\psi\|_{H^{r}} \nonumber\\
&&\,\,\,\quad\quad\quad
\leq C(\|u_{0}\|_{H^{s}}^{2}
+\|u_{0}\|_{H^{s}}\|v_{0}\|_{H^{s}}
+\|v_{0}\|_{H^{s}}^{2}
+\|u_{0}\|_{H^{s}}^{2})
\|\varphi\|_{H^{r}}\|\psi\|_{H^{r}} \nonumber\\
&&\,\,\,\quad\quad\quad
\leq C\rho^{2}\|\varphi\|_{H^{r}}\|\psi\|_{H^{r}}
 \end{eqnarray}
recalling $\|u_{0}\|_{H^{s}}\leq \rho$
and $\|v_{0}\|_{H^{s}}\leq \rho$
in the assumption of Theorem \ref{HolderContinuityThm}.

As $-1/2<r\leq 1/2$ and $r\leq s-2$,
using Lemma \ref{21lem4} and (\ref{sizeEstimate}) yields
\begin{eqnarray}\label{energy1-1Esimate2}
 &&|(\ref{energy1-1})|
 \leq C
\left(\|g\|_{H^{r+1}}+\|u^{2}\|_{H^{r+1}}\right)
\|\varphi\|_{H^{r}}\|\psi\|_{H^{r}} \nonumber\\
&&\,\,\,\quad\quad\quad
\leq
C(\|w\|_{H^{r+1}}^{2}+\|w\|_{H^{r+1}}
\|z\|_{H^{r+1}}+\|z\|_{H^{r+1}}^{2}+\|u\|_{H^{r+1}}^{2})
\|\varphi\|_{H^{r}}\|\psi\|_{H^{r}} \nonumber\\
&&\,\,\,\quad\quad\quad
\leq C\rho^{2}\|\varphi\|_{H^{r}}\|\psi\|_{H^{r}}.
 \end{eqnarray}

When $-1\leq r\leq -1/2$,
employing the inequality
\begin{eqnarray}\label{Lemma2Inequality}
\|fg\|_{H^{r}}\leq
 C\|f\|_{H^{s-1}}\|g\|_{H^{r}}\quad (-1\leq r\leq 0, s>3/2)
\end{eqnarray}
provided by Lemma 2 in \cite{HimonasMantzavinos2014JNS},
one finds
\begin{eqnarray}\label{energy1-1Esimate3}
 &&|(\ref{energy1-1})|
 \leq C\left(\|g\|_{H^{s-1}}+\|u^{2}\|_{H^{s-1}}\right)
\|\varphi\|_{H^{r}}\|\psi\|_{H^{r}} \nonumber\\
&&\,\,\,\quad\quad\quad
\leq
C(\|w\|_{H^{s-1}}^{2}+\|w\|_{H^{s-1}}
\|z\|_{H^{s-1}}+\|z\|_{H^{s-1}}^{2}+\|u\|_{H^{s-1}}^{2})
\|\varphi\|_{H^{r}}\|\psi\|_{H^{r}} \nonumber\\
&&\,\,\,\quad\quad\quad
\leq C\rho^{2}\|\varphi\|_{H^{r}}\|\psi\|_{H^{r}}.
\end{eqnarray}
Using similar procedure
as above, we can
estimate  (\ref{energy1-2})-(\ref{energy1-4})
as
\begin{eqnarray}\label{energy1-234Estimate}
&&|(\ref{energy1-2})|\leq C\rho^{2}\|\varphi\|_{H^{r}}^{2},\quad
|(\ref{energy1-3})| \leq C\rho\|\varphi\|_{H^{r}}\|\psi\|_{H^{r}},\quad
|(\ref{energy1-4})| \leq
 C\rho\|\varphi\|_{H^{r}}\|\psi\|_{H^{r}}
\end{eqnarray}
for  $r \in\{-1 \leq r \leq-1 / 2\}
\cup\{-1 / 2<r \leq 1 / 2, \, r \leq s-2\} \cup\{1/2<r \leq s-1\}$.

For the nonlocal terms (\ref{energy1-5}), a direct calculation yields
\begin{subequations}\label{energy1-5Estimate0}
\begin{align}
|(\ref{energy1-5})|
\leq&
 \bigg| \int_{\mathbb{R}}
  D^{r-2}\left[\left(\frac{k_1}{3}+\frac{k_2}{2}\right)g \psi\right] \cdot D^{r} \varphi \mathrm{d} x\bigg|\label{energy1-5Estimate0-1}\\
&+\bigg| \int_{\mathbb{R}}
 D^{r-2} \partial_{x}
 \bigg[\left(\left(\frac{2k_1}{3}+k_{2}\right) f
 +\left(k_{1}+\frac{3k_2}{2}\right)z^{2}
 +k_{3}(u+v)\right) \varphi \nonumber \\
&\quad+\left(k_{1}+\frac{3k_2}{2}\right)u(w+z)\psi+\frac{k_3}{2}(w+z) \psi\bigg] \cdot D^{r}
\varphi \mathrm{d} x \bigg|.\label{energy1-5Estimate0-2}
\end{align}
\end{subequations}
Processing similarly as the local terms, one has
\begin{eqnarray}\label{energy1-5Estimate0-1-1}
&&(\ref{energy1-5Estimate0-1})
\leq
C\|g\psi\|_{{H^{r}}}\|\varphi\|_{{H^{r}}}\nonumber\\
&&\quad\quad\quad\leq
\begin{cases}
 C\|g\|_{{H^{r}}}\|\psi\|_{{H^{r}}}\|\varphi\|_{{H^{r}}},\quad (1/2<r\leq s-1);\nonumber\\
C\|g\|_{{H^{r+1}}}\|\psi\|_{{H^{r}}}\|\varphi\|_{{H^{r}}},\quad
(-1/2<r\leq 1/2, r\leq s-2);\nonumber\\
C\|g\|_{{H^{s-1}}}\|\psi\|_{{H^{r}}}\|\varphi\|_{{H^{r}}},\quad
(-1\leq r\leq -1/2)
\end{cases}\nonumber\\
&&\quad\quad\quad\leq C\|g\|_{{H^{s-1}}}\|\psi\|_{{H^{r}}}\|\varphi\|_{{H^{r}}}
\leq C\rho^{2}\|\psi\|_{{H^{r}}}\|\varphi\|_{{H^{r}}}
\end{eqnarray}
and
\begin{eqnarray}\label{energy1-5Estimate0-1-2}
&&(\ref{energy1-5Estimate0-2})
\leq C\|[f +z^{2}+(u+v)]\varphi+u(w+z)\psi+(w+z) \psi\|_{{H^{r}}}\|\varphi\|_{{H^{r}}}\nonumber\\
&&\quad\quad\quad\leq
\begin{cases}
 C(\|f\|_{{H^{r}}}+\|z\|_{{H^{r}}}^{2} +\|u\|_{{H^{r}}}+\|v\|_{{H^{r}}}) \|\varphi\|_{{H^{r}}}^{2}\nonumber\\
 \quad+C(\|u\|_{{H^{r}}}+1)(\|w\|_{{H^{r}}}+\|z\|_{{H^{r}}}) \|\psi\|_{{H^{r}}}\|\varphi\|_{{H^{r}}},
 \quad (1/2<r\leq s-1),\nonumber\\
 C(\|f\|_{{H^{r+1}}}+\|z\|_{{H^{r+1}}}^{2} +\|u\|_{{H^{r+1}}}+\|v\|_{{H^{r+1}}}) \|\varphi\|_{{H^{r}}}^{2}\nonumber\\
 \quad+C(\|u\|_{{H^{r+1}}}+1)(\|w\|_{{H^{r+1}}}+\|z\|_{{H^{r+1}}})  \|\psi\|_{{H^{r}}}\|\varphi\|_{{H^{r}}},\,\,\,
(-1/2<r\leq 1/2, r\leq s-2),\nonumber\\
C(\|f\|_{{H^{s-1}}}+\|z\|_{{H^{s-1}}}^{2}
 +\|u\|_{{H^{s-1}}}+\|v\|_{{H^{s-1}}})
 \|\varphi\|_{{H^{r}}}^{2}\nonumber\\
 \quad+C(\|u\|_{{H^{s-1}}}+1)(\|w\|_{{H^{s-1}}}+\|z\|_{{H^{s-1}}})
  \|\psi\|_{{H^{r}}}\|\varphi\|_{{H^{r}}},\quad (-1\leq r\leq -1/2)
\end{cases}\nonumber\\
&&\quad\quad\quad\leq C(\|f\|_{{H^{s-1}}}+\|z\|_{{H^{s-1}}}^{2}
 +\|u\|_{{H^{s-1}}}+\|v\|_{{H^{s-1}}})
 \|\varphi\|_{{H^{r}}}^{2}\nonumber\\
&&\quad\quad\quad\quad+C(\|u\|_{{H^{s-1}}}+1)(\|w\|_{{H^{s-1}}}+\|z\|_{{H^{s-1}}})
  \|\psi\|_{{H^{r}}}\|\varphi\|_{{H^{r}}}\nonumber\\
&&\quad\quad\quad\leq C(\rho^{2}+\rho)(\|\varphi\|_{{H^{r}}}^{2}
+\|\psi\|_{{H^{r}}}\|\varphi\|_{{H^{r}}}).
\end{eqnarray}
Combining (\ref{energy1}) and
(\ref{energy1-1Esimate1})-(\ref{energy1-5Estimate0-1-2}),
one has arrived at
\begin{eqnarray}\label{HrOfVarphi}
\frac{\mbox{d}}{\mbox{d}t}\|\varphi\|_{{H^{r}}}
\leq
C(\rho^{2}+\rho)(\|\varphi\|_{{H^{r}}}
+\|\psi\|_{{H^{r}}})
\end{eqnarray}
for  $r \in\{-1 \leq r \leq-1 / 2\}
\cup\{-1 / 2<r \leq 1 / 2,\, r \leq s-2\} \cup\{1 / 2<r \leq s-1\}$. \\

We are now in a position to estimate (\ref{energy2}). The following Calderon-Coifman-Meyer
commutator estimate \cite{Taylor2003PAMS,HimonasMantzavinos2014JNS}
\begin{eqnarray}\label{commutator}
\left\|\left[D^{r} \partial_{x}, f\right] g\right\|_{L^{2}}
 \leq C\|f\|_{H^{s-1}}\|g\|_{H^{r}},(0\leq r+1\leq s-1;s-1>3/2)
\end{eqnarray}
will be used, where $[A,B]=AB-BA$ represents the commutator.

To deal with (\ref{energy2-1}), we first recast it as
\begin{subequations}\label{energy2-1Recast}
\begin{align}
(\ref{energy2-1})=&
\int_{\mathbb{R}}\left[D^{r} \partial_{x},\left(k_{1}w^{2}-(k_{1}+k_{2})u^{2}-k_{3}u\right)\right] \psi \cdot D^{r} \psi \mathrm{d} x\label{energy2-1Recast1}\\
&+\int_{\mathbb{R}}\left(k_{1}w^{2}-(k_{1}+k_{2})u^{2}-k_{3}u\right) D^{r} \partial_{x} \psi \cdot D^{r} \psi \mathrm{d} x \label{energy2-1Recast2}\\
&-\int_{\mathbb{R}} D^{r}\left[\psi \partial_{x}\left(k_{1}w^{2}-(k_{1}+k_{2})u^{2}-k_{3}u\right)\right] \cdot D^{r} \psi \mathrm{d} x.\label{energy2-1Recast3}
\end{align}
\end{subequations}
Invoking (\ref{commutator}), one finds
\begin{eqnarray}\label{energy2-1Recast1Estimate}
&&|(\ref{energy2-1Recast1})|
\leq C\|[D^{r} \partial_{x},(k_{1}w^{2}-(k_{1}+k_{2})u^{2}-k_{3}u)]
\psi\|_{L^{2}}
\|\psi\|_{H^{r}}\nonumber\\
&&\quad\quad\quad\quad\leq C\|k_{1}w^{2}-(k_{1}+k_{2})u^{2}-k_{3}u\|_{H^{s-1}}
\|\psi\|_{H^{r}}^{2}\leq C(\rho^{2}+\rho)\|\psi\|_{H^{r}}^{2}
\end{eqnarray}
for  $r \in\{-1 \leq r \leq-1 / 2, r+s \geq 2\} \cup\{-1 / 2<r \leq 1 / 2, r \leq s-3\} \cup\{r>$
$1 / 2, r \leq s-2\}$.

Integration by parts and Sobolev embedding inequality yields
\begin{eqnarray}\label{energy2-1Recast2Estimate}
&&|(\ref{energy2-1Recast2})|
\leq C\| \partial_{x}[k_{1}w^{2}-(k_{1}+k_{2})u^{2}-k_{3}u]\|_{L^{\infty}}
\|\psi\|_{H^{r}}^{2}\nonumber\\
&&\quad\quad\quad\quad\leq C\|k_{1}w^{2}-(k_{1}+k_{2})u^{2}-k_{3}u\|_{H^{s-1}}
\|\psi\|_{H^{r}}^{2}\leq C(\rho^{2}+\rho)\|\psi\|_{H^{r}}^{2}.
\end{eqnarray}
(\ref{energy2-1Recast3}) can be estimated as
\begin{eqnarray}\label{energy2-1Recast3Estimate}
&&|(\ref{energy2-1Recast3})|
\leq C\|\psi\partial_{x}(k_{1}w^{2}-(k_{1}+k_{2})u^{2}-k_{3}u)\|_{H^{r}}
\|\psi\|_{H^{r}}\nonumber\\
&&\quad\quad\quad\quad\leq
\begin{cases}
C\|\partial_{x}(k_{1}w^{2}-(k_{1}+k_{2})u^{2}-k_{3}u)\|_{H^{r}}
\|\psi\|_{H^{r}}^{2}, (1/2<r\leq s-2),\\
C\|\partial_{x}(k_{1}w^{2}-(k_{1}+k_{2})u^{2}-k_{3}u)\|_{H^{r+1}}
\|\psi\|_{H^{r}}^{2}\\
\quad \leq C\|k_{1}w^{2}-(k_{1}+k_{2})u^{2}-k_{3}u\|_{H^{r+2}}
\|\psi\|_{H^{r}}^{2}, (-1/2< r\leq 1/2; r\leq s-3),\\
C\|\partial_{x}(k_{1}w^{2}-(k_{1}+k_{2})u^{2}-k_{3}u)\|_{H^{s-2}}
\|\psi\|_{H^{r}}^{2}, (-1\leq r\leq -1/2;r+s\geq 2)
\end{cases}\nonumber\\
&&\quad\quad\quad\quad
\leq C\|k_{1}w^{2}-(k_{1}+k_{2})u^{2}-k_{3}u\|_{H^{s-1}}
\|\psi\|_{H^{r}}^{2}\leq C(\rho^{2}+\rho)\|\psi\|_{H^{r}}^{2},
\end{eqnarray}
where the inequality \cite{HimonasMantzavinos2014JNS}
\begin{eqnarray}\label{Lemma4Inequality}
\|f g\|_{H^{r}}
\leq c_{r, s}\|f\|_{H^{s-2}}\|g\|_{H^{r}}\quad (-1\leq r\leq 0,\, r+s\geq 2,\, s>5/2)
\end{eqnarray}
has been used to handle the third case in the brace.

Combining (\ref{energy2-1Recast1Estimate})-(\ref{energy2-1Recast3Estimate})
produces
\begin{eqnarray}
|(\ref{energy2-1})|
\leq C(\rho^{2}+\rho)\|\psi\|_{H^{r}}^{2}
\end{eqnarray}
for $r \in\{-1 \leq r \leq-1 / 2, r+s \geq 2\} \cup\{-1 / 2<r \leq 1 / 2,\, r \leq s-3\} \cup\{r>$
$1 / 2, r \leq s-2\}$.

We next estimate (\ref{energy2-2}).
In fact, we can derive
\begin{eqnarray}\label{energy2-2Estimate}
&&|(\ref{energy2-2})|
\leq C\bigg\|\bigg[k_{1}(w+z) \partial_{x} z-\left(k_{1}+\frac{k_2}{2}\right)u(w+z)
-\frac{k_3}{2}(w+z)\bigg] \psi\bigg\|_{H^{r}}
\|\psi\|_{H^{r}}\nonumber\\
&&\quad\quad\quad\quad\leq
\begin{cases}
C[\|(w+z) \partial_{x} z\|_{H^{r}}
+\|u(w+z)\|_{H^{r}}
+\|(w+z)\|_{H^{r}}]
\|\psi\|_{H^{r}}^{2}\\
\quad\leq C[\|(w+z) \partial_{x} z\|_{H^{s-2}}
+\|u(w+z)\|_{H^{s-2}}
+\|(w+z)\|_{H^{s-2}}]
\|\psi\|_{H^{r}}^{2}, \\
\quad\quad\quad\quad\quad\quad\quad\quad
\quad\quad\quad\quad\quad\quad\quad\quad
\quad\quad\quad\quad\quad\quad\quad
(1/2<r\leq s-2);\\
C[\|(w+z) \partial_{x} z\|_{H^{r+1}}
+\|u(w+z)\|_{H^{r+1}}
+\|(w+z)\|_{H^{r+1}}]
\|\psi\|_{H^{r}}^{2}\\
\quad\leq C[\|(w+z) \partial_{x} z\|_{H^{s-2}}
+\|u(w+z)\|_{H^{s-2}}
+\|(w+z)\|_{H^{s-2}}]
\|\psi\|_{H^{r}}^{2}, \\
\quad\quad\quad\quad\quad\quad\quad\quad
\quad\quad\quad\quad\quad\quad\quad\quad
\quad\quad\quad
(-1/2< r\leq 1/2,\, r\leq s-3),\\
C[\|(w+z) \partial_{x} z\|_{H^{s-2}}
+\|u(w+z)\|_{H^{s-2}}
+\|(w+z)\|_{H^{s-2}}]
\|\psi\|_{H^{r}}^{2},\\
\quad\quad\quad\quad\quad\quad\quad\quad
\quad\quad\quad\quad\quad\quad\quad\quad
\quad\quad\quad
 (-1\leq r\leq -1/2,\,r+s\geq 2)
\end{cases}\nonumber\\
&&\quad\quad\quad\quad
\leq C[(\|w\|_{H^{s-2}}+\|z\|_{H^{s-2}})
(\|\partial_{x} z\|_{H^{s-2}}
+\|u\|_{H^{s-2}}
+1)]
\|\psi\|_{H^{r}}^{2}\nonumber\\
&&\quad\quad\quad\quad\leq C(\rho^{2}+\rho)\|\psi\|_{H^{r}}^{2}.
\end{eqnarray}

Using similar method as estimating
(\ref{energy2-2}),
one deduces
\begin{eqnarray}\label{energy2-34Estimate}
&&|(\ref{energy2-3})|
\leq  C(\rho^{2}+\rho)\|\psi\|_{H^{r}}^{2}, \qquad |(\ref{energy2-4})|
\leq  C(\rho^{2}+\rho)\|\psi\|_{H^{r}}^{2}
\end{eqnarray}
for $r \in\{-1 \leq r \leq-1 / 2, r+s \geq 2\} \cup\{-1 / 2<r \leq 1 / 2, r \leq s-3\} \cup\{r>$
$1 / 2, r \leq s-2\}$.

Finally, the estimation of the
nonlocal term (\ref{energy2-6})
is similar to (\ref{energy1-5}) and
we thus obtain
\begin{eqnarray}\label{energy2-6Estimate}
&&|(\ref{energy2-6})|
\leq  C(\rho^{2}+\rho)(\|\psi\|_{H^{r}}^{2}
+\|\psi\|_{H^{r}}\|\varphi\|_{H^{r}})
\end{eqnarray}
for  $r \in\{-1 \leq r \leq-1 / 2\}\cup\{-1 / 2<r \leq 1 / 2,\, r \leq s-2\}
\cup\{1/2<r \leq s-1\}$.

Gathering (\ref{energy2-1Recast})-(\ref{energy2-6Estimate})
into (\ref{energy2}), one finds
\begin{eqnarray}\label{HrOfPsi}
\frac{\mbox{d}}{\mbox{d}t}\|\psi\|_{{H^{r}}}
\leq
C(\rho^{2}+\rho)(\|\varphi\|_{{H^{r}}}
+\|\psi\|_{{H^{r}}})
\end{eqnarray}
for $r \in\{-1 \leq r \leq-1 / 2,\, r+s \geq 2\} \cup\{-1 / 2<r \leq 1 / 2,\, r \leq s-3\}
\cup\{1 / 2<r\leq s-2\}$.

Combining (\ref{HrOfVarphi}) and (\ref{HrOfPsi}) yields
\begin{eqnarray}\label{HrOfVarphiAndPsi1}
\frac{\mbox{d}}{\mbox{d}t}
(\|\varphi\|_{{H^{r}}}+\|\psi\|_{{H^{r}}})
\leq
C(\rho^{2}+\rho)(\|\varphi\|_{{H^{r}}}
+\|\psi\|_{{H^{r}}})
\end{eqnarray}
or
\begin{eqnarray}\label{HrOfVarphiAndPsi2}
\|\varphi\|_{{H^{r}}}+\|\psi\|_{{H^{r}}}
\leq
C(\|\varphi(0)\|_{{H^{r}}}
+\|\psi(0)\|_{{H^{r}}})e^{(\rho^{2}+\rho)t},\quad (t\in [0,T))
\end{eqnarray}
for $r \in\{-1 \leq r \leq-1 / 2, r+s \geq 2\} \cup\{-1 / 2<r \leq 1 / 2, r \leq s-3\} \cup\{r>$
$1 / 2, r \leq s-2\}$.

Therefore, we deduce
\begin{eqnarray}\label{HrOfVarphiAndPsi3}
\|u-v\|_{H^{r+1}}
\leq C\|u_{0}-v_{0}\|_{H^{r+1}}e^{(\rho^{2}+\rho)T}.
\end{eqnarray}
Replacing $r+1$ with $r$ in (\ref{HrOfVarphiAndPsi3})
leads to
\begin{eqnarray}\label{HrOfVarphiAndPsi4}
\|u-v\|_{H^{r}}
\leq C\|u_{0}-v_{0}\|_{H^{r}}e^{(\rho^{2}+\rho)T}
\end{eqnarray}
for $r \in\{0 \leq r \leq 1 / 2,\, r+s \geq 3\}\cup\{1 / 2<r \leq 3 / 2,\, r \leq s-2\}
\cup\{3 / 2<r \leq s-1\}$. Thus the Lipschitz continuity in the region
$D_{1}$ has been established.

We next show the H\"{o}lder continuity in the region $D_{2}\cup D_{3} \cup D_{4}.$
The method is interpolation based on the Lipschitz continuity proved previously.
The inequality \cite{HimonasMantzavinos2014JNS}
\begin{eqnarray}\label{Lemma5Inequality}
\|f\|_{H^{\sigma}} \leq\|f\|_{H^{\sigma_{1}}}^{\frac{\sigma_{2}-\sigma}{\sigma_{2}-\sigma_{1}}}
\|f\|_{H^{\sigma_{2}}}^{\frac{\sigma-\sigma_{1}}{\sigma_{2}-\sigma_{1}}},
\quad (\sigma_{1}<\sigma< \sigma_{2})
\end{eqnarray}
will be used frequently in the remaining parts of this section.

As $(s,r)\in D_{2}$,
we have
\begin{eqnarray*}
\|u-v\|_{H^{r}}
&\leq&\|u-v\|_{H^{3-s}} \nonumber\\
&\leq& C\|u_{0}-v_{0}\|_{H^{3-s}} e^{(\rho^{2}+\rho) T}
\quad(\text{by}\,\, (\ref{HrOfVarphiAndPsi4}))\nonumber\\
&\leq& C\|u_{0}-v_{0}\|_{H^{r}}^{\frac{2 s-3}{s-r}}
\|u_{0}-v_{0}\|_{H^{s}}^{\frac{3-s-r}{s-r}} e^{(\rho^{2}+\rho)T}
\quad(\text{by}\,\, (\ref{Lemma5Inequality}))\nonumber\\
&\leq& C \rho^{\frac{3-s-r}{s-r}}
\|u_{0}-v_{0}\|_{H^{r}}^{\frac{2 s-3}{s-r}}e^{(\rho^{2}+\rho) T}.
\end{eqnarray*}

As $(s,r)\in D_{3}$,
one finds $s-2\leq r<s,$
consequently
\begin{eqnarray*}
\|u-v\|_{H^{r}}
&\leq&\|u-v\|_{H^{s-2}}^{\frac{s-r}{2}}\|u-v\|_{H^{s}}^{\frac{r-s+2}{2}} \quad(\text{by}\,\, (\ref{Lemma5Inequality}))\nonumber\\
&\leq& C\rho^{\frac{r-s+2}{2}}\left\|u_{0}-v_{0}\right\|_{H^{s-2}}^{\frac{s-r}{2}} e^{(\rho^{2}+\rho) T}
\quad(\text{by}\,\, (\ref{HrOfVarphiAndPsi4}))\nonumber\\
&\leq& C\rho^{\frac{r-s+2}{2}}\left\|u_{0}-v_{0}\right\|_{H^{r}}^{\frac{s-r}{2}} e^{(\rho^{2}+\rho)T}.
\end{eqnarray*}

For $(s,r)\in D_{4}$, there holds
\begin{eqnarray*}
\|u-v\|_{H^{r}}
&\leq&\|u-v\|_{H^{s-1}}^{s-r}\|u-v\|_{H^{s}}^{r-s+1} \quad(\text{by}\,\, (\ref{Lemma5Inequality}))\nonumber\\
&\leq& C\rho^{r-s+1}\left\|u_{0}-v_{0}\right\|_{H^{s-1}}^{s-r} e^{(\rho^{2}+\rho) T}
\quad(\text{by}\,\, (\ref{HrOfVarphiAndPsi4}))\nonumber\\
&\leq& C\rho^{r-s+1}\left\|u_{0}-v_{0}\right\|_{H^{r}}^{s-r} e^{(\rho^{2}+\rho)T}.
\end{eqnarray*}
We thus complete the proof of Theorem \ref{HolderContinuityThm}.

\section{Proof of Theorem \ref{WaveBreakingThm2}}

We would like to prove Theorem \ref{WaveBreakingThm2} by an induction argument on $s>1/2.$ Let us first of all consider the case $s\in (1/2,1).$  Applying Lemma \ref{22lem5} to Eq.~(\ref{meq})
yields
\begin{eqnarray}\label{blowup0}
\|m\|_{H^{s}}
\leq \|m_{0}\|_{H^{s}}
+C\int_{0}^{t}U^{\prime}(\tau)\|m(\tau)\|_{H^{s}}\mbox{d}\tau+C\int_{0}^{t}
\|2k_{1}u_{x}m^{2}+3k_{2}uu_{x}m+2k_{3}u_{x}m\|_{H^{s}}
\mbox{d}\tau,
\end{eqnarray}
where $U(t)=\int_{0}^{t}\|\partial_{x}(k_{1}(u^{2}-u_{x}^{2})+k_{2}u^{2}+k_{3}u)\|_{L^{\infty}}\mbox{d}\tau$.
Using
 $u=(1-\partial_{x}^{2})^{-1} m=p * m$
and $\|p\|_{L^{1}}=\|\partial_{x} p\|_{L^{1}}=1$,
one finds after employing the Young inequality
that for   $s \in \mathbb{R}$
\begin{eqnarray}\label{blowup1}
\begin{array}{l}
\|u\|_{L^{\infty}}+\left\|u_{x}\right\|_{L^{\infty}}+\left\|u_{x x}\right\|_{L^{\infty}} \leq C\|m\|_{L^{\infty}},\v\\
\|u\|_{H^{s}}+\left\|u_{x}\right\|_{H^{s}}+\left\|u_{x x}\right\|_{H^{s}} \leq C\|m\|_{H^{s}}.
\end{array}
\end{eqnarray}

Invoking (\ref{blowup1}), we derive
\begin{eqnarray}\label{blowup3}
U^{\prime}(t)
&=&\|\partial_{x}(k_{1}(u^{2}-u_{x}^{2})+k_{2}u^{2}+k_{3}u)
\|_{L^{\infty}}\nonumber\\
&\leq& C(\|u_{x}m\|_{L^{\infty}}+\|u^{2}\|_{L^{\infty}}
+\|u\|_{L^{\infty}}\leq C(\|m\|_{L^{\infty}}+\|m\|_{L^{\infty}}^{2})
\end{eqnarray}
and
\begin{eqnarray}\label{blowup4}
&&\|2k_{1}u_{x}m^{2}+3k_{2}uu_{x}m+2k_{3}u_{x}m\|_{H^{s}}\nonumber\\
&&\leq C(\|u_{x}\|_{L^{\infty}}\|m^{2}\|_{H^{s}}
+\|u_{x}\|_{H^{s}}\|m^{2}\|_{L^{\infty}})+C(\|u_{x}\|_{L^{\infty}}\|um\|_{H^{s}}
+\|u_{x}\|_{H^{s}}\|um\|_{L^{\infty}})\nonumber\\
&&\quad+C(\|u_{x}\|_{L^{\infty}}\|m\|_{H^{s}}
+\|u_{x}\|_{H^{s}}\|m\|_{L^{\infty}})\nonumber\\
&&\leq C(\|u_{x}\|_{L^{\infty}}\|m\|_{L^{\infty}}\|m\|_{H^{s}}
+\|u_{x}\|_{H^{s}}\|m\|_{L^{\infty}}^{2})\nonumber\\
&&\quad+C[\|u_{x}\|_{L^{\infty}}
(\|u\|_{L^{\infty}}\|m\|_{H^{s}}+\|m\|_{L^{\infty}}\|u\|_{H^{s}})
+\|u_{x}\|_{H^{s}}\|m\|_{L^{\infty}}\|u\|_{L^{\infty}}
]\nonumber\\
&&\quad+C[\|u_{x}\|_{L^{\infty}}\|m\|_{H^{s}}
+\|m\|_{L^{\infty}}\|u_{x}\|_{H^{s}}]\nonumber\\
&&\leq C\|m\|_{H^{s}}(\|m\|_{L^{\infty}}^{2}+\|m\|_{L^{\infty}}).
\end{eqnarray}

Substituting Eqs.~(\ref{blowup3}) and (\ref{blowup4}) into Eq.~(\ref{blowup0}) leads to
\begin{eqnarray*}
\|m\|_{H^{s}}
\leq \|m_{0}\|_{H^{s}}
+C\int_{0}^{t}\|m\|_{H^{s}}
(\|m\|_{L^{\infty}}^{2}+\|m\|_{L^{\infty}})\mbox{d}\tau,
\end{eqnarray*}
which further generates
\begin{eqnarray}\label{blowup5}
\|m\|_{H^{s}}\leq \|m_{0}\|_{H^{s}}\exp\left\{C\int_{0}^{t}(\|m\|_{L^{\infty}}^{2}+\|m\|_{L^{\infty}})\mbox{d}\tau\right\}
\end{eqnarray}
by means of the Gronwall inequality.

Accordingly, if $\int_{0}^{T^{*}}\|m(\tau)\|_{L^{\infty}}^{2} \mbox{d} \tau<\infty$ for the maximal existence time $T^{*}<\infty$,
then the inequality (\ref{blowup5}) indicates that
$\underset{t \rightarrow T^{*}}{\limsup }\|m(t)\|_{H^{s}}<\infty,$
which contradicts our assumption on $T^{*}$. Thus we complete the proof of Theorem \ref{WaveBreakingThm2}
for the case $s\in (1/2,1).$

We are next concerned with the case $s\in [1,2).$  Applying $\partial_{x}$ to Eq.~(\ref{meq}) yields
\begin{eqnarray}
&&\partial_{t}m_{x}
+[k_{1}(u^{2}-u_{x}^{2})+k_{2}u^{2}+k_{3}u]\partial_{x}m_{x}\nonumber\\
&&\qquad =-\partial_{x}[k_{1}(u^{2}-u_{x}^{2})+k_{2}u^{2}+k_{3}u]\partial_{x}m
-2k_{1}u_{xx}m^{2}
-4k_{1}u_{x}mm_{x}\nonumber\\
&&\qquad\quad-3k_{2}u_{x}^{2}m
-3k_{2}uu_{xx}m
-3k_{2}uu_{x}m_{x}
-2k_{3}u_{xx}m
-2k_{3}u_{x}m_{x},
\end{eqnarray}
which combined with Lemma \ref{22lem5} enables us to conclude
\begin{eqnarray}\label{blowup6}
\|\partial_{x}m\|_{H^{s-1}}
&\leq& \|\partial_{x}m_{0}\|_{H^{s-1}}
+C\int_{0}^{t}U^{\prime}(\tau)
\|\partial_{x}m(\tau)\|_{H^{s-1}}\mbox{d}\tau\nonumber\\
&&+C\int_{0}^{t}
\|-\partial_{x}[k_{1}(u^{2}-u_{x}^{2})
+k_{2}u^{2}+k_{3}u]\partial_{x}m\|_{H^{s-1}}
\mbox{d}\tau\nonumber\\
&&+C\int_{0}^{t}
\|2k_{1}u_{xx}m^{2}
+4k_{1}u_{x}mm_{x}
+3k_{2}u_{x}^{2}m
+3k_{2}uu_{xx}m
+3k_{2}uu_{x}m_{x}\|_{H^{s-1}}
\mbox{d}\tau\nonumber\\
&&+C\int_{0}^{t}
\|-2k_{3}u_{xx}m
-2k_{3}u_{x}m_{x}\|_{H^{s-1}}
\mbox{d}\tau.
\end{eqnarray}
For the integrand of the second integral on the right hand side of (\ref{blowup6}), one finds
\begin{eqnarray}\label{blowup7}
&&\|-\partial_{x}[k_{1}(u^{2}-u_{x}^{2})
+k_{2}u^{2}+k_{3}u]\partial_{x}m\|_{H^{s-1}}\nonumber\\
&&\quad \leq C\|\partial_{x}(u^{2}-u_{x}^{2})\partial_{x}m\|_{H^{s-1}}
+C\|uu_{x}\partial_{x}m\|_{H^{s-1}}
+C\|u_{x}\partial_{x}m\|_{H^{s-1}}.
\end{eqnarray}
Employing Lemma \ref{21lem5}  and (\ref{blowup1}),
  there holds
\begin{eqnarray}\label{blowup8}
\|u_{x}\partial_{x}m\|_{H^{s-1}}
&\leq& \|u_{x}m\|_{H^{s}}
+\|u_{xx}m\|_{H^{s-1}}\nonumber\\
&\leq& \|u_{x}\|_{L^{\infty}}\|m\|_{H^{s}}
+\|m\|_{L^{\infty}}\|u_{x}\|_{H^{s}}+\|u_{xx}\|_{L^{\infty}}\|m\|_{H^{s}}
+\|m\|_{L^{\infty}}\|u_{xx}\|_{H^{s-1}}\nonumber\\
&\leq& \|m\|_{H^{s}}\|m\|_{L^{\infty}}.
\end{eqnarray}
Similar argument leads to
\begin{eqnarray}\label{blowup9}
\|uu_{x}\partial_{x}m\|_{H^{s-1}}
&\leq& \|uu_{x}m\|_{H^{s}}
+\|u_{x}^{2}m\|_{H^{s-1}}
+\|uu_{xx}m\|_{H^{s-1}}\nonumber\\
&\leq& \|u\|_{L^{\infty}}\|u_{x}m\|_{H^{s}}
+\|u_{x}m\|_{L^{\infty}}\|u\|_{H^{s}}
+\|u_{x}^{2}\|_{L^{\infty}}\|m\|_{H^{s-1}}\nonumber\\
&&+\|m\|_{L^{\infty}}\|u_{x}^{2}\|_{H^{s-1}}
+\|uu_{xx}\|_{L^{\infty}}\|m\|_{H^{s-1}}
+\|m\|_{L^{\infty}}\|uu_{xx}\|_{H^{s-1}}\nonumber\\
&\leq& \|u\|_{L^{\infty}}
(\|u_{x}\|_{L^{\infty}}\|m\|_{H^{s}}
+\|m\|_{L^{\infty}}\|u_{x}\|_{H^{s}})
+\|u\|_{H^{s}}\|m\|_{L^{\infty}}^{2}\nonumber\\
&&+2\|m\|_{L^{\infty}}^{2}\|m\|_{H^{s}}
+\|m\|_{L^{\infty}}\|u_{x}\|_{L^{\infty}}\|u_{x}\|_{H^{s-1}}\nonumber\\
&&+\|m\|_{L^{\infty}}(\|u\|_{L^{\infty}}\|u_{xx}\|_{H^{s-1}}
+\|u_{xx}\|_{L^{\infty}}\|u\|_{H^{s-1}})\nonumber\\
&\leq& C\|m\|_{H^{s}}\|m\|_{L^{\infty}}^{2}
\end{eqnarray}
and
\begin{align}\label{blowup10}
\|\partial_{x}(u^{2}\!-\!u_{x}^{2})m_{x}\|_{H^{s\!-\!1}}
\leq& C\|u_{x}mm_{x}\|_{H^{s-1}}\leq C\|u_{x}(m^{2})_{x}\|_{H^{s-1}} \no\\
\leq& C\|\partial_{x}(u_{x}m^{2})\|_{H^{s-1}}
+C\|u_{xx}m^{2}\|_{H^{s-1}}\leq C\|u_{x}m^{2}\|_{H^{s}}
+C\|u_{xx}m^{2}\|_{H^{s-1}}\nonumber\\
\leq& C\|u_{x}\|_{L^{\infty}}\|m^{2}\|_{H^{s}}
\!\!+\!C\|m^{2}\|_{L^{\infty}}\|u_{x}\|_{H^{s}}\!\!+\!C\|u_{xx}\|_{H^{s\!-\!1}}\|m^{2}\|_{L^{\infty}}
\!\!+\!C\|m^{2}\|_{H^{s\!-\!1}}\|u_{xx}\|_{L^{\infty}}\nonumber\\
\leq& C\|m\|_{H^{s}}\|m\|_{L^{\infty}}^{2}.
\end{align}

Substituting Eqs.~(\ref{blowup8})-(\ref{blowup10}) into Eq.~(\ref{blowup7}), we derive
\begin{eqnarray}\label{blowup11}
&&\quad\|-\partial_{x}[k_{1}(u^{2}-u_{x}^{2})
+k_{2}u^{2}+k_{3}u]\partial_{x}m\|_{H^{s-1}}\leq C\|m\|_{H^{s}}\|m\|_{L^{\infty}}^{2}.
\end{eqnarray}
On the other hand, the  procedure of handling Eqs.~(\ref{blowup8})-(\ref{blowup10}) enables one to easily
deal with the integrands of the third and fourth
integrals  on the right hand side of Eq.~(\ref{blowup6})
and obtain the corresponding  results in analogy with Eq.~(\ref{blowup11}).

Plugging these estimates into Eq.~(\ref{blowup6}) yields
\begin{eqnarray*}
\|\partial_{x}m\|_{H^{s-1}}
\leq \|\partial_{x}m_{0}\|_{H^{s-1}}
+C\int_{0}^{t}
\|m\|_{H^{s}}(\|m\|_{L^{\infty}}^{2}
+\|m\|_{L^{\infty}})
\mbox{d}\tau.
\end{eqnarray*}
Then arguing similarly as the case for $s\in (1/2,1),$
one can prove Theorem \ref{WaveBreakingThm2} in the case of
$s\in [1,2).$

We finally consider the case $s\geq 2.$
Assume $2 \leq k \in \mathbb{N}.$
Suppose (\ref{blowup00}) holds when $k-1 \leq$ $s<k.$
Using induction, we should prove the validity of (\ref{blowup00})
 for $k\leq s<k+1.$ Applying $\partial_{x}^{k}$ to Eq.~(\ref{meq}) yields
\begin{eqnarray}\label{blowupk0}
&&\partial_{t}\partial_{x}^{k}m
+[k_{1}(u^{2}-u_{x}^{2})
+k_{2}u^{2}+k_{3}u]\partial_{x}\partial_{x}^{k}m\nonumber\\
&&\qquad =-\sum_{l=0}^{k-1}\!C_{k}^{l}\partial_{x}^{k-l}[k_{1}\!(u^{2}\!-\!u_{x}^{2})
\!+\!k_{2}u^{2}\!+\!k_{3}u]\partial_{x}^{l+1}m\!-\!2k_{1}\partial_{x}^{k}(u_{x}m^{2})
\!-\!3k_{2}\partial_{x}^{k}(uu_{x}m)
\!-\!2k_{3}\partial_{x}^{k}(u_{x}m).\qquad\,\,
\end{eqnarray}
Applying Lemma \ref{22lem5} again yields
\begin{eqnarray}\label{blowupk1}
\|\partial_{x}^{k}m\|_{H^{s-k}}
&\leq& \|\partial_{x}^{k}m_{0}\|_{H^{s-k}}
+C\int_{0}^{t}U^{\prime}(\tau)
\|\partial_{x}^{k}m(\tau)\|_{H^{s-k}}
\mbox{d}\tau\nonumber\\
&&+C\int_{0}^{t}
\left\|\sum_{l=0}^{k-1}C_{k}^{l}\partial_{x}^{k-l}
[k_{1}(u^{2}-u_{x}^{2})
+k_{2}u^{2}+k_{3}u]\partial_{x}^{l+1}m\right\|_{H^{s-k}}
\mbox{d}\tau\nonumber\\
&&+C\int_{0}^{t}
\|-2k_{1}\partial_{x}^{k}(u_{x}m^{2})
-3k_{2}\partial_{x}^{k}(uu_{x}m)
-2k_{3}\partial_{x}^{k}(u_{x}m)\|_{H^{s-k}}
\mbox{d}\tau.
\end{eqnarray}
Invoking Lemma \ref{21lem5} and the Sobolev embedding inequality
produces
\begin{eqnarray}\label{blowupk2}
&&\left\|\sum_{l=0}^{k-1}C_{k}^{l}\partial_{x}^{k-l}
[k_{1}(u^{2}-u_{x}^{2})
+k_{2}u^{2}+k_{3}u]\partial_{x}^{l+1}m\right\|_{H^{s-k}}\nonumber\\
&&\leq \sum_{l=0}^{k-1}C_{k}^{l}\|\partial_{x}^{k-l}
[k_{1}(u^{2}-u_{x}^{2})
+k_{2}u^{2}+k_{3}u]\partial_{x}^{l+1}m\|_{H^{s-k}}\nonumber\\
&&\leq \sum_{l=0}^{k-1}C_{k}^{l}
(\|k_{1}(u^{2}-u_{x}^{2})
+k_{2}u^{2}+k_{3}u\|_{H^{s-l+1}}\|\partial_{x}^{l}m\|_{L^{\infty}}\nonumber\\
&&\quad+\|\partial_{x}^{k-l}
[k_{1}(u^{2}-u_{x}^{2})
+k_{2}u^{2}+k_{3}u]\|_{L^{\infty}}\|m\|_{H^{s-k+l+1}}
)\nonumber\\
&&\leq \sum_{l=0}^{k-1}C_{k}^{l}
\bigg[(\|u\|_{L^{\infty}}+\|u_{x}\|_{L^{\infty}}+1)
(\|u\|_{H^{s-l+1}}+\|u_{x}\|_{H^{s-l+1}})
\|m\|_{H^{l+\frac{1}{2}+\epsilon_{0}}}\nonumber\\
&&\quad
+(\|u\|_{H^{k-l+\frac{1}{2}+\epsilon_{0}}}^{2}
+\|u_{x}\|_{H^{k-l+\frac{1}{2}+\epsilon_{0}}}^{2}
+1
)\|m\|_{H^{s-k+l+1}}
\bigg]\nonumber\\
&&\leq C\|m\|_{H^{k-\frac{1}{2}+\epsilon_{0}}}^{2}\|m\|_{H^{s}}
+C\|m\|_{H^{k-\frac{1}{2}+\epsilon_{0}}}\|m\|_{H^{s}}\nonumber\\
&&\leq C(\|m\|_{H^{k-\frac{1}{2}+\epsilon_{0}}}^{2}+1)\|m\|_{H^{s}},
\end{eqnarray}
where $\epsilon_{0}\in (0,1/8)$ so that $H^{\frac{1}{2}+\epsilon_{0}}(\mathbb{R}) \hookrightarrow L^{\infty}(\mathbb{R})$.

Utilizing similar method as the estimation of
(\ref{blowupk2}), we have
\begin{eqnarray}\label{blowupk3}
\|\partial_{x}^{k}(u_{x}m^{2})\|_{H^{s-k}}
&\leq& \|u_{x}m^{2}\|_{H^{s}}\leq \|m^{2}\|_{H^{s}}\|u_{x}\|_{L^{\infty}}
+\|u_{x}\|_{H^{s}}\|m^{2}\|_{L^{\infty}}\leq \|m\|_{L^{\infty}}^{2}\|m\|_{H^{s}},
\end{eqnarray}
\begin{eqnarray}\label{blowupk4}
\|\partial_{x}^{k}(uu_{x}m)\|_{H^{s-k}}
&\leq& \|uu_{x}m\|_{H^{s}}\leq \|u\|_{L^{\infty}}\|u_{x}m\|_{H^{s}}
+\|u\|_{H^{s}}\|u_{x}m\|_{L^{\infty}}\nonumber\\
&\leq& \|u\|_{L^{\infty}}(\|u_{x}\|_{L^{\infty}}
\|m\|_{H^{s}}
+\|u_{x}\|_{H^{s}}\|m\|_{L^{\infty}})+\|u\|_{H^{s}}\|u_{x}\|_{L^{\infty}}\|m\|_{L^{\infty}}\nonumber\\
&\leq& C\|m\|_{L^{\infty}}^{2}\|m\|_{H^{s}}
\end{eqnarray}
and
\begin{eqnarray}\label{blowupk5}
\|\partial_{x}^{k}(u_{x}m)\|_{H^{s-k}}
\leq \|u_{x}m\|_{H^{s}}
\leq \|m\|_{H^{s}}\|u_{x}\|_{L^{\infty}}
+\|u_{x}\|_{H^{s}}\|m\|_{L^{\infty}}
\leq C\|m\|_{L^{\infty}}\|m\|_{H^{s}}.
\end{eqnarray}

Combining Eqs.~(\ref{blowupk1})-(\ref{blowupk5}) deduces
\begin{eqnarray} \label{mk}
\|m\|_{H^{s}}
\leq \|m_{0}\|_{H^{s}}
+C\int_{0}^{t}
(\|m\|_{H^{k-\frac{1}{2}+\epsilon_{0}}}^{2}+1)\|m\|_{H^{s}}(\tau)
\mbox{d}\tau.
\end{eqnarray}
Then the Gronwall inequality implies
\begin{equation}\label{blowupk6}
\|m(t)\|_{H^{s}}
\leq\left\|m_{0}\right\|_{H^{s}}
 \exp \left\{C \int_{0}^{t}
 \left(\|m(\tau)\|_{H^{k-\frac{1}{2}+\epsilon_{0}}}^{2}+1\right) \mbox{d} \tau\right\}.
\end{equation}
Therefore, if the maximal existence time $T^{\star}<\infty$ satisfies
$\int_{0}^{T^{\star}}\|m(\tau)\|_{L^{\infty}}^{2} \mbox{d} \tau<\infty,$
then the uniqueness of the solution provided by
Theorem \ref{WaveBreakingThm1} enables us to conclude the
uniform boundedness of
$\|m(t)\|_{H^{k-\frac{1}{2}+\epsilon_{0}}}$ in
$t \in (0, T^{*})$ taking into account  the induction assumption, which along with (\ref{blowupk6}) indicates
$\limsup _{t \rightarrow T^{*}}\|m(t)\|_{H^{s}}<\infty,$ a contradiction. We thus complete the proof of Theorem \ref{WaveBreakingThm2}.

\section{Proof of Theorem \ref{blowupQuantityThm}}

To prove Theorem \ref{blowupQuantityThm}, we will first of all establish a kind of  global conservative property
associated to Eq.~(\ref{gmCH-Novikov}). Note that the corresponding  trajectory equation
emanating  from
$x$ reads
\begin{eqnarray}\label{conservative0}
\begin{cases}
\dfrac{d}{dt}q(t,x)=[k_{1}(u^{2}-u_{x}^{2})+k_{2}u^{2}+k_{3}u](t, q(t,x)),\quad (t, x)\in (0, T^*)\times\mathbb{R}, \v \\
q(0,x)=x,\quad x\in\mathbb{R}.
\end{cases}
\end{eqnarray}
for the flow generated by $k_{1}(u^{2}-u_{x}^{2})+k_{2}u^{2}+k_{3}u$.

The following lemma shows that,
the momentum $m(t,x)$ will not change sign
for all $t\in[0,T^{*})$
as long as the initial value
$m_{0}=(1-\partial_{x}^{2})u_{0}$
 not changing sign.
This kind of conservative property of the momentum
$m$ is significant  in proving the
wave-breaking Theorem \ref{WaveBreakingThm3}
in the  non-sign-changing case.

\begin{lem}\label{conLem1}
 Suppose $u_{0} \in H^{s}(\mathbb{R})$  with $s>\frac{5}{2}$.  Let $T^{*}>0$ be the maximal existence time of the
 solution $u$ to Eq.~(\ref{gmCH-Novikov}).  Then there exists a  unique solution $q \in C^{1}([0, T^{*}) \times \mathbb{R}; \mathbb{R})$  to Eq.~(\ref{conservative0}) satisfying
\begin{equation}\label{conservative2}
 q_{x}(t, x) =\exp \left( \int_{0}^{t} \left[(2k_{1} m +2k_{2} u  +k_{3})u_x\right](\tau, q(\tau, x)) \mbox{d} \tau\right)>0,
\end{equation}
namely, $q(t, \cdot)$ is an increasing diffeomorphism of $\mathbb{R}$. Moreover, the expression of $m(t, q(t, x))$ is
\begin{equation}\label{conservative21}
m(t, q(t, x))
=m_{0}(x)\exp \left( -\int_{0}^{t} \left[(2k_{1} m +3k_{2} u  +2k_{3})u_x\right](\tau, q(\tau, x)) \mbox{d} \tau\right),\,\, (t,x)\in [0, T^*)\times \mathbb{R}.
\end{equation}
\end{lem}

\begin{proof}

Applying $\partial_{x}$ to Eq.~(\ref{conservative0}) leads to
\begin{eqnarray}\label{conservative3}
\begin{cases}
\dfrac{d}{dt}q_{x}(t,x)=[(2k_{1}m+2k_{2}u+k_{3})u_x](t,q(t,x))q_{x}(t,x),\v\\
q_{x}(0,x)=1.
\end{cases}
\end{eqnarray}
Solving (\ref{conservative3}) produces the solution given by (\ref{conservative2}).
The Sobolev embedding inequality
gives for $T^{\prime}<T^{*}$
\begin{eqnarray*}
\sup_{(s,x)\in [0,T^{\prime})\times \mathbb{R}}
|(2k_{1}m+2k_{2}u+k_{3})u_{x}(s,x)|<\infty,
\end{eqnarray*}
which along with (\ref{conservative2}) yields $q_{x}(t,x)\geq \exp(-Ct),\, (t,x)\in [0,T^*)\times \mathbb{R}$ for some $C>0$, which states clearly that
$q(t, \cdot)$ is an increasing diffeomorphism of $\mathbb{R}$ before the blow-up time,
that is, (\ref{conservative2}) holds.

It follows from Eqs.~(\ref{gmCH-Novikov}) and (\ref{conservative0}) that
\begin{eqnarray}\label{conservative6}
\frac{d}{dt}m(t, q(t,x))
&=&m_{t}(t, q(t,x))
+m_{x}(t, q(t,x))q_{t}(t,x)\nonumber\\
&=&m_{t}(t, q(t,x))
+m_{x}(t, q(t,x))
[k_{1}(u^{2}-u_{x}^{2})
+k_{2}u^{2}+k_{3}u](t, q(t,x))\nonumber\\
&=&[m_{t}
+m_{x}(k_{1}(u^{2}-u_{x}^{2})
+k_{2}u^{2}+k_{3}u)](t, q(t,x))\nonumber\\
&=&[(-2k_{1}m-3k_{2}u-2k_{3})u_{x}](t,q(t,x))m(t, q(t,x)).
\end{eqnarray}
Solving Eq.~(\ref{conservative6}) yields Eq.~(\ref{conservative21}). We thus complete the proof of Lemma \ref{conLem1}.
\end{proof}

 {\it Proof of Theorem \ref{blowupQuantityThm}}
  From the expression of $m(t,q(t,x))$   in Eq.~(\ref{conservative21}),
  we conclude that
  if there exists a positive constant $K_{1}$
  such that
  \begin{equation*}
  \inf _{x \in \mathbb{R}}(2 k_{1} u_{x} m(t, x)+3 k_{2} u u_{x}(t, x)+2k_{3}u_{x}(t,x))
  \geq -K_{1}, \quad 0\leq t\leq T^{*},
  \end{equation*}
then
  \begin{equation*}
  \begin{aligned}
\|m(t)\|_{L^{\infty}}
&=\|m(t, q(t, x))\|_{L^{\infty}}= m_{0}(x)
\exp \left(- \int_{0}^{t}\left(2 k_{1}m+3 k_{2} u+2k_{3}\right)u_x(\tau, q(\tau, x)) \mbox{d}\tau\right) \\
& \leq e^{ K_{1} T^{*}}\|m_{0}\|_{L^{\infty}},
\end{aligned}
  \end{equation*}
  which combined with Theorem \ref{WaveBreakingThm2}
  implies that $m(t,x)$ will not blow up in a finite time.

  However, if (\ref{blowupQuantity}) holds true,
  then Theorem \ref{WaveBreakingThm1} and the Sobolev embedding
  ensure that $m(t,x)$ will  blow up in a finite time.
  We thus complete the proof
   of Theorem \ref{blowupQuantityThm}.

\section{Proof of Theorem \ref{finalThm}}

For the proof of Theorem \ref{finalThm}, we need to introduce the following Lemma\cite{ConstantinEscher1998AM}:

\begin{lem}\label{lemFinal}
 Let $T>0$ and $w(t,x) \in \mathcal{C}^{1}\left([0, T) ; H^{2}(\mathbb{R})\right) .$ Then for every $t \in[0, T)$ there exists at least one point $x=\zeta(t) \in \mathbb{R}$ such that
$W(t):=\inf _{x \in \mathbb{R}}\left[w_{x}(t, x)\right]=w_{x}(t, \zeta(t)).$
and the function $W(t)$ is almost everywhere differentiable on $(0, T)$ with
$$
\frac{\mbox{d} W(t)}{\mbox{d}t}=w_{t x}(t, \zeta(t)), \text { a.e. on }(0, T).
$$
\end{lem}

We next prove (\ref{MuniformUpperBound}). We only need to consider the case $s=3,$ and the general case  follows
 by the density argument. Let $M(t,x)=2 k_{1} u_{x} m(t, x)+3 k_{2} u u_{x}(t, x)+2k_{3}u_{x}(t,x)$. Firstly, Corollary \ref{cor1} shows that
$M\in \mathcal{C}([0,T);H^{s})\cap
\mathcal{C}^{1}([0,T);H^{s-1})$
with $T$ representing  the maximal existence
time of the corresponding solution.
 Given $t\in [0, T )$,
Lemma \ref{lemFinal} implies there exists some point
$x_{0}(t)\in \mathbb{R}$ such that
\begin{eqnarray}\label{Mmaximum}
&M\left(t, x_{0}(t)\right)=\sup _{x \in \mathbb{R}} M(t, x),\quad
\text {i.e.,}\quad M_{x}\left(t, x_{0}(t)\right)=0, \text { a.e. on }(0, T).
\end{eqnarray}
The condition $s>\frac{1}{2}$
indicates
 $H^{s}(\mathbb{R}) \hookrightarrow \mathcal{C}_{0}(\mathbb{R}),$
the space of all continuous functions on $\mathbb{R}$ vanishing as $|x| \rightarrow \infty$.
We deduce  in light of Corollary \ref{cor1} that
\begin{eqnarray}\label{Mpositive}
M\left(t, x_{0}(t)\right) \geq 0 \quad \text { for all } \quad t \in[0, T).
\end{eqnarray}

Recall that the map $q(t, \cdot)$ is an increasing diffeomorphism of $\mathbb{R}$ in Lemma~\ref{conLem1}, therefore, there exists some  $\eta_{0}=\eta_{0}(t) \in \mathbb{R}$ satisfying $q\left(t, \eta_{0}(t)\right)=x_{0}(t).$
At the point $\left(t, q\left(t, \eta_{0}\right)\right)=\left(t, x_{0}(t)\right),$ one has
\begin{eqnarray}\label{20210115}
\frac{d}{\mbox{d} t} m\left(t, x_{0}(t)\right)=-(m M)\left(t, x_{0}(t)\right)
\end{eqnarray}
after substituting $x=\eta_{0}(t)$ into Eq.~(\ref{conservative6}).

Integrating (\ref{20210115}) and using
 (\ref{Mpositive}) enable us to conclude
\begin{eqnarray}\label{mx0t1}
m(t, x_{0}(t))=m_{0}(x_{0}(0))
\exp \left(-\int_{0}^{t}M(\tau, x_{0}(\tau))
 \mbox{d}\tau\right)\leq m_{0}(x_{0}(0))
\leq \sup_{x\in \mathbb{R}}m_{0}(x).
\end{eqnarray}

Next, let's  recall the relations between $u$ and $m$:
\begin{equation}\label{character5}
  u(t, x)
  =p * m(t, x),\quad
  u_{x}(t, x)=p_{x} * m(t, x).
\end{equation}
On the other hand, the conditions  $m_{0} \geq 0$ and
$m_{0}(x_{0})>0$
imply $m(t, x) \geq 0$ and $m(t, q(t,x_{0}))>0$
in view of the expression of $m(t, q(t,x))$
in Eq.~(\ref{conservative21}).
So one derives from (\ref{character5}) that
 $u(t, x) \geq 0$
and $u(t, q(t,x_{0}))>0.$
Furthermore, it follows from (\ref{pPlusMinus}) and (\ref{character5}) that
\begin{equation}\label{character6}
  u\pm u_{x}=2 p_{\mp} * m \geq 0,
\end{equation}
which indicates that $|u_{x}(t,x)|\leq u(t,x)$.

 From (\ref{mx0t1}) and the inequality
$|u_{x}|\leq u$, one obtains
\begin{eqnarray}\label{Mover}
M(t, x_{0}(t))
&=&2 k_{1} u_{x} m(t, x_{0}(t))
+3 k_{2} u u_{x}(t, x_{0}(t))
+2k_{3}u_{x}(t,x_{0}(t))\nonumber\\
&\leq& 2 k_{1}\|u\|_{L^{\infty}} \sup_{x\in \mathbb{R}}m_{0}(x)
+3 k_{2} \|u\|_{L^{\infty}}^{2}
+2k_{3}\|u\|_{L^{\infty}}\nonumber\\
&\leq& 2 k_{1}\|u_{0}\|_{H^{1}} \sup_{x\in \mathbb{R}}m_{0}(x)
+3 k_{2} \|u_{0}\|_{H^{1}}^{2}
+2k_{3}\|u_{0}\|_{H^{1}},
\end{eqnarray}
which gives (\ref{MuniformUpperBound}). We thus complete the proof of Theorem \ref{finalThm}.

\section{Proof of Theorem \ref{WaveBreakingThm3}}

The main idea of the proof of Theorem \ref{WaveBreakingThm3} is due to Ref.~\cite{ChenGuoLiuQu2016JFA},
where the wave-breaking problem of the generalized mCH-CH equation was considered. We need to compute the  dynamics of some important quantities  along the characteristic $q(t,x)$ defined in (\ref{conservative0}).

We first calculate  the
dynamics of $u$  along the characteristic.
For this purpose,
applying $(1-\partial_{x}^{2})$
to $u_{t}
+(k_{1}(u^{2}-u_{x}^{2})
+k_{2}u^{2}+k_{3}u)u_{x}$
and employing equation (\ref{gmCH-Novikov})
lead to
\begin{eqnarray}\label{character1}
\begin{array}{rl}
&(1-\partial_{x}^{2})[u_{t}+(k_{1}(u^{2}-u_{x}^{2})+k_{2}u^{2}+k_{3}u)u_{x}] \v\\
&\qquad = m_{t}+(1-\partial_{x}^{2})[(k_{1}(u^{2}-u_{x}^{2})+k_{2}u^{2}+k_{3}u)u_{x}]\v\\
&\qquad=k_{1}(2u_{x}^{2}u_{xxx}+4u_{x}u_{xx}^{2}-2uu_{x}u_{xx}-2u^{2}u_{x}-2u_{x}^{3})\v\\
&\qquad\quad+k_{2}(
-3u^{2}u_{x}
-3uu_{x}u_{xx}
-2u_{x}^{3})
+k_{3}(-2uu_{x}
-u_{x}u_{xx}
).
\end{array}\end{eqnarray}
Using the definitions of $p_{\pm}$ defined in (\ref{pPlus}),
one derives from (\ref{character1}) that
\begin{eqnarray}\label{character2}
&&u_{t}
+(k_{1}(u^{2}-u_{x}^{2})
+k_{2}u^{2}+k_{3}u)u_{x}\nonumber\\
&&\quad=-\frac{2}{3}k_{1}u_{x}^{3}+\left(\frac{k_1}{3}+\frac{k_2}{2}\right)[p_{+}*(u-u_{x})^{3}-p_{-}*(u+u_{x})^{3}]
-k_{3}p_{x}*\left(u^{2}+\frac{1}{2}u_{x}^{2}\right),
\end{eqnarray}
which is just the dynamics of $u$
 along the characteristic.

We next operate $\partial_{x}$ to (\ref{character2})
to deduce
\begin{eqnarray}\label{character3}
&& u_{xt}
+(k_{1}(u^{2}-u_{x}^{2})
+k_{2}u^{2}+k_{3}u)u_{xx}\nonumber\\
&&\quad=k_{1}\left(\frac{1}{3}u^{3}-uu_{x}^{2}\right)
+k_{2}\left(\frac{1}{2}u^{3}-\frac{1}{2}uu_{x}^{2}\right)
+k_{3}\left(u^{2}-\frac{1}{2}u_{x}^{2}\right)\nonumber\\
&&\quad\quad -\left(\frac{k_1}{3}+\frac{k_2}{2}\right)[p_{+}*(u-u_{x})^{3}+p_{-}*(u+u_{x})^{3}]-k_{3}p*\left(u^{2}+\frac{1}{2}u_{x}^{2}\right),
\end{eqnarray}
which is the dynamics of $u_{x}$ along the characteristic.


The inequality  $|u_{x}(t,x)|\leq u(t,x)$ and
the conservative relation $\|u\|_{H^{1}}=\|u_{0}\|_{H^{1}}$
implies $u$ will not blow up.
Consequently,
it is enough to consider the quantity
$N=mu_{x}$ on account of Theorem \ref{blowupQuantityThm}.

Define
\begin{eqnarray*}
\begin{array}{ll}
\widehat{u}(t)
=u(t, q(t, x_{0})), &\widehat{u_{x}}(t)
=u_{x}(t, q(t, x_{0})), \v\\
\widehat{m}(t)
=m(t, q(t, x_{0})), &
\widehat{N}(t)=(m u_{x})(t, q(t, x_{0})).
\end{array}
\end{eqnarray*}
We will trace the dynamics of $N(t)$
 along the characteristics emanating from $x_{0}.$
It follows from Eqs.~(\ref{character3}) and (\ref{meq})  that
\begin{eqnarray}\label{character7}
\widehat{N}^{\prime}(t)
&=&\widehat{m}^{\prime}\widehat{u}_{x}
+\widehat{m}\widehat{u}_{x}^{\prime}\nonumber\\
&=&-(2k_{1}\widehat{m}^{2}
+3k_{2}\widehat{m}\widehat{u}
+2k_{3}\widehat{m})\widehat{u}_{x}^{2}\nonumber\\
&&+k_1\widehat{m}\left(\frac{1}{3}\widehat{u}^{3}
-\widehat{u}\widehat{u}_{x}^{2}\right)
+\frac{k_2}{2}\widehat{m}\left(\widehat{u}^{3}-\widehat{u}\widehat{u}_{x}^{2}\right)
+k_3\widehat{m}\left(\widehat{u}^{2}
-\frac{1}{2}\widehat{u}_{x}^{2}\right)\nonumber\\
&&-\left(\frac{k_1}{3}+\frac{k_2}{2}\right)\widehat{m}[p_{+}*(u-u_{x})^{3}+p_{-}*(u+u_{x})^{3}](t,q(t,x_{0}))\nonumber\\
&&-k_{3}\widehat{m}
p*\left(u^{2}+\frac{1}{2}u_{x}^{2}\right)(t,q(t,x_{0})),
\end{eqnarray}
where $``\prime"$ represents the derivative $\partial_{t}+[k_{1}(u^{2}-u_{x}^{2})+k_{2}u^{2}+k_{3}u]\partial_{x}$ along the
characteristic.

Since $p*(u^{2}+\frac{1}{2}u_{x}^{2})\geq \frac{1}{2}u^{2},$ Eq.~(\ref{character7}) can be controlled by
\begin{eqnarray}\label{character8}
\widehat{N}^{\prime}(t)\leq -2k_{1}\widehat{N}^{2}(t)
+\frac{k_1}{3}\widehat{m}\widehat{u}(\widehat{u}^{2}-3\widehat{u}_{x}^{2})+\frac{k_2}{2}\widehat{m}\widehat{u}
(\widehat{u}^{2}-7\widehat{u}_{x}^{2})+\frac{k_3}{2}\widehat{m}(\widehat{u}^{2}-5\widehat{u}_{x}^{2})
\end{eqnarray}
in view of Eq.~(\ref{character6}).

We want every of the last three terms on the right hand side of (\ref{character8})
to be no greater than zero to establish the Riccati-like inequality $\widehat{N}^{\prime}(t)\leq C_0\widehat{N}^{2}(t)$
with $C_{0}$ being negative,
so we consider the dynamics of the ratio
$\widehat{u}_{x}/(\widehat{u}+\gamma)$ along the characteristics with some  parameter $\gamma\geq 0$ to be settled later.

Direct computation leads to
\begin{eqnarray}\label{character9}
\left(\frac{\widehat{u}_{x}}{\widehat{u}+\gamma}\right)^{\prime}(t)
&=&\frac{1}{(\widehat{u}+\gamma)^{2}}\bigg[\frac{k_{1}}{3}(\widehat{u}^{2}-\widehat{u}_{x}^{2})(\widehat{u}^{2}-2\widehat{u}_{x}^{2})
+\frac{\gamma k_{1}}{3}\widehat{u}(\widehat{u}^{2}-3\widehat{u}_{x}^{2})\nonumber\\
&&+\frac{k_{2}}{2}\widehat{u}^{2}(\widehat{u}^{2}-\widehat{u}_{x}^{2})
+\frac{\gamma k_{2}}{2}\widehat{u}(\widehat{u}^{2}-\widehat{u}_{x}^{2})+k_{3}\left(\widehat{u}^{3}-\frac{1}{2}\widehat{u}\widehat{u}_{x}^{2}\right)
+\frac{\gamma k_{3}}{2}(2\widehat{u}^{2}-\widehat{u}_{x}^{2})\nonumber\\
&&-\left(\frac{k_1}{3}+\frac{k_2}{2}\right)(\widehat{u}+\widehat{u}_{x}+\gamma)p_{+}*(u-u_{x})^{3}(t,q(t,x_{0}))\nonumber\\
&&-\left(\frac{k_1}{3}+\frac{k_2}{2}\right)(\widehat{u}-\widehat{u}_{x}+\gamma)p_{-}*(u+u_{x})^{3}(t,q(t,x_{0}))\nonumber\\
&&-k_{3}[(\widehat{u}+\widehat{u}_{x}+\gamma)p_{+}+(\widehat{u}-\widehat{u}_{x}+\gamma)p_{-}]*\left(u^{2}+\frac{1}{2}u_{x}^{2}\right)(t,q(t,x_{0}))
\bigg].
\end{eqnarray}
Employing the inequality
$p_{\pm} *\left(u^{2}+\frac{1}{2} u_{x}^{2}\right) \geq \frac{1}{4} u^{2}$, one can dominate Eq.~(\ref{character9}) as
\begin{eqnarray}\label{character10}
\left(\frac{\widehat{u}_{x}}{\widehat{u}+\gamma}\right)^{\prime}(t)
&\leq& \frac{1}{(\widehat{u}+\gamma)^{2}}\bigg[\frac{k_{1}}{3}(\widehat{u}^{2}-\widehat{u}_{x}^{2})
(\widehat{u}^{2}-2\widehat{u}_{x}^{2})+\frac{k_{1}}{3}\gamma\widehat{u}(\widehat{u}^{2}-3\widehat{u}_{x}^{2})\nonumber\\
&&+\frac{k_{2}}{2}\widehat{u}^{2}(\widehat{u}^{2}-\widehat{u}_{x}^{2})
+\frac{k_{2}}{2}\gamma\widehat{u}(\hat{u}^{2}-\widehat{u}_{x}^{2})
+\frac{k_{3}}{2}(\widehat{u}^{3}-\widehat{u}\widehat{u}_{x}^{2})+\frac{k_{3}}{2}\gamma(\widehat{u}^{2}-\widehat{u}_{x}^{2})
\bigg]\nonumber\\
&=&\frac{k_1}{3(\widehat{u}+\gamma)^{2}}\bigg\{(\widehat{u}^{2}-\widehat{u}_{x}^{2})\bigg[
\frac{2k_{1}+3k_{2}}{2k_{1}}\left(\widehat{u}+\frac{3}{2}\cdot\frac{k_{2}\gamma+k_{3}}{2k_{1}+3k_{2}}
\right)^{2}\nonumber\\
&&+\frac{3}{2k_1}\left(  k_3\gamma-\frac{3(k_2\gamma+k_3)^2}{4(2k_{1}+3k_{2})}\right)-2\widehat{u}_{x}^{2}\bigg]
+\gamma\widehat{u}(\widehat{u}^{2}-3\widehat{u}_{x}^{2})
\bigg\},
\end{eqnarray}
where we want the constant term in the brackets to be zero, so the parameter $\gamma$ should satisfy
\begin{eqnarray}\label{expressionOfTheParameter}
k_3\gamma=\frac{3( k_2\gamma+k_3)^2}{4(2k_{1}+3k_{2})}.
\end{eqnarray}

In accordance with Eq.~(\ref{expressionOfTheParameter}), one should discuss the values of $k_{1}, k_{2}$ and $k_{3}$
in the following four cases (Note that the truly blow-up quantity
is $mu_{x}$, so $k_{1}$ should be always greater than zero):

{\it Case 1}: For the case $k_{1},\,k_{2},\,k_{3}>0$, we obtain from Eq.~(\ref{expressionOfTheParameter}) that the  parameter $\gamma$ is
\begin{eqnarray}\label{character11}
\gamma=\frac{k_3}{3k_2^2}\left(4k_1+3k_2\pm 2\sqrt{4k_1^2+6k_1k_2}\right)>0,
\end{eqnarray}
from which (\ref{character10}) then becomes
\begin{eqnarray}\label{character12}
\left(\frac{\widehat{u}_{x}}{\widehat{u}+\gamma}\right)^{\prime}(t)
\leq \frac{k_1}{3(\widehat{u}+\gamma)^{2}}\bigg\{
(\widehat{u}^{2}-\widehat{u}_{x}^{2})\bigg[
\frac{2k_{1}+3k_{2}}{2k_{1}}\left(
\widehat{u}+\frac{3}{2}\cdot\frac{k_{2}\gamma+k_{3}}{2k_{1}+3k_{2}}
\right)^{2}\!-\!2\widehat{u}_{x}^{2}\!\bigg]+\gamma\widehat{u}(\widehat{u}^{2}-3\widehat{u}_{x}^{2})
\bigg\}.
\end{eqnarray}
We need the terms in the brackets in Eq.~(\ref{character12})
be negative to ensure the non-increasing of
$\widehat{u}_{x}/(\widehat{u}+\gamma)$, namely
\begin{eqnarray}\label{character13}
\sqrt{\frac{2k_{1}+3k_{2}}{2k_{1}}}\left|
\widehat{u}+\frac{3}{2}\cdot\frac{k_{2}\gamma+k_{3}}{2k_{1}+3k_{2}}\right|
< \sqrt{2}|\widehat{u}_{x}|
\end{eqnarray}
as well as
\begin{eqnarray}\label{character14}
\hat{u}_{x}(0)
< -\frac{1}{\sqrt{2}}(\widehat{u}(0)+\gamma)
\end{eqnarray}
 to be used to ensure
the non-positivity  of the last three terms
on the right hand side of   (\ref{character8}).

In view of Eq.~(\ref{character13}), it seems the initial condition should be
\begin{eqnarray}\label{character15}
\widehat{u}_{x}(0)
< -\frac{1}{\sqrt{2}}
\sqrt{\frac{2k_{1}+3k_{2}}{2k_{1}}}\left[\widehat{u}(0)+\frac{3}{2}\cdot\frac{k_{2}\gamma+k_{3}}{2k_{1}+3k_{2}}\right].
\end{eqnarray}
However, this condition cannot satisfy Eq.~(\ref{character14}).

Therefore, one can consider multiply a factor $\alpha>1$
to be determined later on the right of Eq.~(\ref{character15}).
Then simple computation yields
\begin{eqnarray}\label{character16}
\widehat{u}_{x}(0)
&<& -\frac{\alpha}{\sqrt{2}}\sqrt{\frac{2k_{1}+3k_{2}}{2k_{1}}}\left[\widehat{u}(0)+\frac{1}{2}\cdot\frac{3k_{2}\gamma
+3k_{3}}{2k_{1}+3k_{2}}\right]\nonumber\\
&=& -\frac{\alpha}{\sqrt{2}}\sqrt{\frac{2k_{1}+3k_{2}}{2k_{1}}}\widehat{u}(0)-\frac{\alpha}{\sqrt{2}}\sqrt{\frac{2k_{1}+3k_{2}}{2k_{1}}}
\frac{1}{2}\cdot\frac{3k_{2}\gamma+3k_{3}}{2k_{1}+3k_{2}}\nonumber\\
&<& -\frac{1}{\sqrt{2}}\widehat{u}(0)-\frac{\alpha}{\sqrt{2}}\cdot\frac{1}{2}\cdot\frac{3k_{2}\gamma+3k_{3}}{2k_{1}+3k_{2}}\nonumber\\
&<&-\frac{1}{\sqrt{2}}\widehat{u}(0)-\frac{\alpha}{\sqrt{2}}\cdot\frac{1}{2}\cdot\frac{3k_{2}\gamma}{2k_{1}+3k_{2}}.
\end{eqnarray}
If the factor $\alpha$ satisfy
\begin{eqnarray*}
\frac{\alpha}{2}\cdot\frac{3k_{2}}{2k_{1}+3k_{2}}
>1,\quad
\text{or}\quad
\alpha> 2+\frac{4k_{1}}{3k_{2}},
\end{eqnarray*}
then (\ref{character16}) implies
\begin{eqnarray*}
\frac{\widehat{u}_{x}}{\widehat{u}+\gamma}(t)<
\frac{\widehat{u}_{x}}{\widehat{u}+\gamma}(0)< -\frac{1}{\sqrt{2}}.
\end{eqnarray*}

{\it Case 2}: For the case $k_{1}>0$,  $k_{2}=0$,  $k_{3}>0,$ Eq.~(\ref{gmCH-Novikov}) becomes the mCH-CH equation, and one  can find from Eq.~(\ref{expressionOfTheParameter}) that~\cite{ChenGuoLiuQu2016JFA}
\begin{eqnarray}\label{character11-1}
\gamma=\frac{3k_{3}}{8k_{1}},
\end{eqnarray}
and the inequality~(\ref{character13}) then becomes
\begin{eqnarray*}
\left|\widehat{u}+\frac{3k_{3}}{4k_{1}}\right|
< \sqrt{2}|\widehat{u}_{x}|,
\end{eqnarray*}
according to which the initial condition should be
\begin{eqnarray*}
\widehat{u}_{x}(0)
< -\frac{1}{\sqrt{2}}
\left[\widehat{u}(0)+\frac{3k_{3}}{4k_{1}}\right].
\end{eqnarray*}
This condition can satisfy (\ref{character14}) with $\gamma$ defined by (\ref{character11-1}).

\v {\it Case 3}: For the case $k_{1},\, k_{2}>0$,  $k_{3}=0,$ Eq.~(\ref{gmCH-Novikov}) becomes the mCH-Novikov equation, and it follows from (\ref{expressionOfTheParameter}) that $\gamma=0$. Then reasoning similarly as the above two cases,
we derive the initial condition should be
\begin{eqnarray*}
\widehat{u}_{x}(0)
< -\frac{1}{\sqrt{2}}
\sqrt{\frac{2k_{1}+3k_{2}}{2k_{1}}}
\widehat{u}(0).
\end{eqnarray*}

\v{\it Case 4}: For the case $k_{1}>0$,  $k_{2}=k_{3}=0$, Eq.~(\ref{gmCH-Novikov}) becomes the mCH equation, and one deduces~\cite{ChenGuoLiuQu2016JFA}
\begin{eqnarray*}
\gamma=0,\quad \quad \widehat{u}_{x}(0)
< -\frac{1}{\sqrt{2}}
\widehat{u}(0).
\end{eqnarray*}

In each of  the above four cases, one always has
\begin{eqnarray*}
\widehat{u}+\gamma+\sqrt{2}\widehat{u}_{x}<0
\end{eqnarray*}
and
\begin{eqnarray}\label{character17}
\widehat{u}(t)+\sqrt{3}\widehat{u}_{x}(t)<0,\quad
\widehat{u}(t)+\sqrt{5}\widehat{u}_{x}(t)<0,\quad
\widehat{u}(t)+\sqrt{7}\widehat{u}_{x}(t)<0.
\end{eqnarray}
On the other hand,
\begin{eqnarray}\label{character18}
\widehat{u}(t)-\sqrt{3}\widehat{u}_{x}(t)>0,\quad
\widehat{u}(t)-\sqrt{5}\widehat{u}_{x}(t)>0,\quad
\widehat{u}(t)-\sqrt{7}\widehat{u}_{x}(t)>0.
\end{eqnarray}

Substituting (\ref{character17}) and  (\ref{character18}) into (\ref{character8}) leads to
\begin{eqnarray}\label{character19}
\widehat{N}^{\prime}(t)\leq -2k_{1}\widehat{N}^{2}(t).
\end{eqnarray}
Hence,
$\widehat{N}(t)$ will blow up in finite time $T^{*}$ with
\begin{eqnarray*}
T^{*} \leq-\frac{1}{2 k_{1} m_{0}\left(x_{4}\right) u_{0, x}\left(x_{4}\right)},
\end{eqnarray*}
which is obtained after integrating (\ref{character19}). We thus complete the proof of Theorem \ref{WaveBreakingThm3}.

\section{Proof of Theorem \ref{WaveBreakingThm4}}

Theorem \ref{WaveBreakingThm3} just considers  a portion of the blow-up quantity, however we would like to  trace the whole blow-up quantity defined as
\begin{eqnarray}\label{M}
M(t,x)=(2 k_{1} m+3 k_{2} u+2k_{3})u_x
\end{eqnarray}
along the characteristics (Here the main idea is similar to \cite{MiLiuHuangGuo2020JDE,GuiLiuOlverQu2013CMP})
in this section.

Direct computation leads to 
\begin{eqnarray}\label{M1}
&& M_{t}
+[k_{1} (u^{2}-u_{x}^{2})+k_{2} u^{2}+k_{3} u] M_{x}\nonumber\\
&&=2k_{1} u_{x}
\{m_{t}
+[k_{1} (u^{2}-u_{x}^{2})+k_{2} u^{2}+k_{3} u] m_{x}
\}+3k_{2} u_{x}
\{u_{t}
+[k_{1} (u^{2}-u_{x}^{2})+k_{2} u^{2}+k_{3} u]u_{x}
\}\nonumber\\
&&\quad+(2k_{1}m+3k_{2} u+2k_{3})
\{u_{xt}
+[k_{1} (u^{2}-u_{x}^{2})+k_{2} u^{2}+k_{3} u]u_{xx}
\}\nonumber\\
&&  =A_{1}+A_{2},
\end{eqnarray}
where $A_{1}$ is defined as the following local terms
\begin{eqnarray}\label{A1}
\begin{array}{ll}
A_{1}=&\!\!-2k_{1}k_{2}u_{x}^{4}
-4k_{1}^{2}m^{2}u_{x}^{2}
-6k_{1}k_{2}muu_{x}^{2}
-4k_{1}k_{3}mu_{x}^{2} \v\\
&\!\!+(2k_{1}m+3k_{2} u+2k_{3})
\left[k_{1}\left(\frac{1}{3}u^{3}-uu_{x}^{2}\right)
+k_{2}\left(\frac{1}{2}u^{3}-\frac{1}{2}uu_{x}^{2}\right)
+k_{3}\left(u^{2}-\frac{1}{2}u_{x}^{2}\right)
\right]\nonumber
\end{array}
\end{eqnarray}
and $A_{2}$ is the remaining nonlocal terms.

Recall that $m\geq0$, $u\geq0$
and $|u_{x}|\leq u$,
so
we have
\begin{eqnarray}\label{ApM2}
&&A_{1}+M^{2}=(2k_{1}m+3k_{2} u+2k_{3})
\left[k_{1}\left(\frac{1}{3}u^{3}-uu_{x}^{2}\right)
+\frac{k_{2}}{2}\left(u^{3}-uu_{x}^{2}\right)
+k_{3}\left(u^{2}-\frac{1}{2}u_{x}^{2}\right)
\right]\nonumber\\
&&\quad\quad\quad\quad\quad\quad +4k_{1}k_{3}mu_{x}^{2}
+6k_{1}k_{2}muu_{x}^{2}
+9k_{2}^{2}u^{2}u_{x}^{2}
+4k_{3}^{2}u_{x}^{2}
+12k_{2}k_{3}uu_{x}^{2}
-2k_{1}k_{2}u_{x}^{4}\nonumber\\
&&\quad\quad\quad\quad\,\,\leq
(2k_{1}m+3k_{2} u+2k_{3})
\left[\left(\frac{4}{3}k_{1}+k_{2}\right)\|u\|_{L^{\infty}}^{3}
+\frac{3}{2}k_{3}\|u\|_{L^{\infty}}^{2}\right]\nonumber\\
&&\quad\quad\quad\quad\quad\quad+2k_{1}m\cdot 2k_{3}\|u\|_{L^{\infty}}^{2}
+2k_{1}m\cdot 3k_{2}\|u\|_{L^{\infty}}^{3}
+3k_{2}u\cdot 3k_{2}\|u\|_{L^{\infty}}^{3}
+2k_{3}\cdot 2k_{3}\|u\|_{L^{\infty}}^{2}\nonumber\\
&&\quad\quad\quad\quad\quad\quad+3k_{2}u\cdot 4k_{3}\|u\|_{L^{\infty}}^{2}
+3k_{2}u\cdot \frac{2}{3}k_{1}\|u\|_{L^{\infty}}^{3}\nonumber\\
&&\quad\quad\quad\quad\,\,\leq
(2k_{1}m+3k_{2} u+2k_{3})
\left[\left(2k_{1}+7k_{2}
\right)\|u\|_{L^{\infty}}^{3}+
\frac{19}{2}k_{3}\|u\|_{L^{\infty}}^{2}
\right].
\end{eqnarray}
We can similarly deal with
the nonlocal terms $A_{2}$
by noticing $\|p_{\pm}\|_{L_{1}}\leq 1$
after using the convolution-type Young inequality
\begin{eqnarray}\label{A2}
\begin{array}{ll}
A_{2}&\leq 3k_{2}u\left(\dfrac{16k_1}{3}\|u\|_{L^{\infty}}^{3}+8k_{2}\|u\|_{L^{\infty}}^{3}
+\dfrac{3k_3}{2}\|u\|_{L^{\infty}}^{2}\right)\v\\
&\qquad +(2k_{1}m+3k_{2} u+2k_{3})
\left[\dfrac{16k_1}{3}\|u\|_{L^{\infty}}^{3}+8k_{2}\|u\|_{L^{\infty}}^{3}+\dfrac{3k_3}{2}\|u\|_{L^{\infty}}^{2}
\right]\v\\
&\leq(2k_{1}m+3k_{2} u+2k_{3})
\left[\left(\dfrac{32}{3}k_{1}
+16k_{2}\right)\|u\|_{L^{\infty}}^{3}
+3k_{3}\|u\|_{L^{\infty}}^{2}
\right].
\end{array}
\end{eqnarray}
Substituting (\ref{ApM2})-(\ref{A2}) into (\ref{M1}) leads to
\begin{eqnarray}\label{M2}
 M_{t}
+[k_{1} (u^{2}-u_{x}^{2})+k_{2} u^{2}+k_{3} u] M_{x}\leq -M^{2}
+C_{1}(2k_{1}m+3k_{2} u+2k_{3})
\end{eqnarray}
with
\begin{eqnarray}\label{C1}
C_{1}=\left(\frac{38}{3}k_{1}+23k_{2}+\frac{25}{2}k_{3}\right)\left(\|u_{0}\|_{H^{1}}^{2}+\|u_{0}\|_{H^{1}}^{3}\right).
\end{eqnarray}

It follows from (\ref{M2}) that
\begin{eqnarray}\label{MMx0t}
\frac{\mbox{d}}{\mbox{d} t}
 M\left(t, q(t,x_{1})\right)
 \leq
-M^{2}\left(t, q(t,x_{1})\right)
+ C_{1}(2k_{1}m+3k_{2} u+2k_{3})\left(t, q(t,x_{1})\right).
\end{eqnarray}
Invoking Eq.~(\ref{meq}), one derives
\begin{eqnarray}\label{mx0t}
\frac{\mbox{d}}{\mbox{d} t} m\left(t, q(t,x_{1})\right)
=-(m M)\left(t, q(t,x_{1})\right).
\end{eqnarray}
Define $\overline{M}(t):= M\left(t, q\left(t, x_{1}\right)\right),$
$\overline{m}(t):= m(t, q(t, x_{1}))$ and
$\overline{u}(t):= u(t, q(t, x_{1})),$
then (\ref{MMx0t}) and (\ref{mx0t})
can be respectively recast as
\begin{eqnarray}\label{Moverline}
\frac{\mbox{d}}{\mbox{d} t}\overline{M}(t)
\leq -\overline{M}(t)^{2}
+C_{1} (2k_{1}\overline{m}+3k_{2} \overline{u}+2k_{3})
\end{eqnarray}
and
\begin{eqnarray}\label{moverline1}
\frac{\mbox{d}}{\mbox{d} t} \overline{m}(t)
=-\overline{m}(t) \overline{M}(t).
\end{eqnarray}

Consequently, one  derives
\begin{eqnarray}\label{mrecast}
&&\frac{\mbox{d}}{\mbox{d} t}\left(\frac{1}{\overline{m}^2(t)} \frac{\mbox{d}}{\mbox{d} t} \overline{m}(t)\right)=\frac{\mbox{d}}{\mbox{d} t}\left(-\frac{1}{\overline{m}(t)} \overline{M}(t)\right) \nonumber\\
&&=\frac{1}{\overline{m}^2(t)}\left(-\overline{m}(t) \frac{d}{d t} \overline{M}(t)+\overline{M}(t) \frac{d}{d t} \overline{m}(t)\right) \nonumber\\
&& \geq \frac{1}{\overline{m}^2(t)}
\left[\overline{m}(t)\left(\overline{M}^2(t)- C_{1} (2k_{1}\overline{m}(t)+3k_{2}\overline{u}(t)+2k_{3})\right)
-\overline{m}(t) \overline{M}^2(t)\right]\nonumber\\
&&=-2k_{1} C_{1}
-3k_{2} C_{1} \frac{\overline{u}(t)}{\overline{m}(t)}
-2k_{3} C_{1} \frac{1}{\overline{m}(t)}\nonumber\\
&& \geq -2k_{1} C_{1}
-(3k_{2}+2k_{3}) C_{1}\frac{\overline{u}(t)+1}{\overline{m}(t)}.
\end{eqnarray}
If $\frac{\overline{u}(t)+1}{\overline{m}(t)}\leq C_{2}$ holds initially,
then it follows from (\ref{mrecast}) that
\begin{eqnarray}\label{mrecast1}
\frac{\mbox{d}}{\mbox{d} t}\left(\frac{1}{\overline{m}^2(t)} \frac{\mbox{d}}{\mbox{d} t} \overline{m}(t)\right)
\geq  -[2k_{1}+
(3k_{2}+2k_{3})C_{2}] C_{1}.
\end{eqnarray}

Integrating (\ref{mrecast1}) yields
\begin{eqnarray}\label{afterIntegration}
\frac{1}{\bar{m}^2(t)} \frac{d}{d t} \bar{m}(t)\geq C_{0}-[2k_{1}+
(3k_{2}+2k_{3})C_{2}] C_{1} t
=C_{0}-2C_{3} t,
\end{eqnarray}
where
\begin{eqnarray}\label{C0}
C_{0}=
\frac{\overline{m}'(0)}{\overline{m}^2(0)}
=-\frac{\overline{M}(0)}{\overline{m}(0)}
=-\partial_{x}u_{0}(x_{1})\left(2k_{1}\!+\!\frac{3k_{2}u_{0}(x_{1})+2k_{3}}{m_{0}(x_{1})}\right),\,
C_{3}=[k_{1}\!+\!(3k_{2}/2\!+\!k_{3})C_{2}] C_{1}.
\end{eqnarray}

From (\ref{moverline1}) and (\ref{afterIntegration}), one obtains
\begin{eqnarray}\label{afterIntegration1}
\overline{M}(t)
=-\frac{1}{m(t)} \frac{\mbox{d}}{\mbox{d} t} \overline{m}(t) \leq-\overline{m}(t)\left(C_{0}-2C_{3} t\right).
\end{eqnarray}
Integrating (\ref{afterIntegration}) from 0 to $t$ yields
\begin{eqnarray}\label{afterIntegration2}
\frac{1}{\overline{m}(t)}-\frac{1}{\overline{m}(0)}
\leq C_{3} t^{2}-C_{0} t,
\end{eqnarray}
which leads to
\begin{eqnarray}\label{afterIntegration3}
\frac{1}{\overline{m}(t)}
\leq C_{3}\left(t^{2}-\frac{C_{0}}{C_{3}} t+\frac{1}{C_{3} \overline{m}(0)}\right)=C_{3}\left(t^{2}-\frac{C_{0}}{C_{3}} t+\frac{1}{ C_{3} m_{0}\left(x_{1}\right)}\right).
\end{eqnarray}

Solving the algebraic equation
\begin{eqnarray*}
t^{2}-\frac{C_{0}}{C_{3}} t+\frac{1}{ C_{3} m_{0}\left(x_{1}\right)}=0
\end{eqnarray*}
leads to
\begin{eqnarray}\label{root1}
t_-:=\frac{C_{0}}{2 C_{3}}-\frac{1}{2} \sqrt{\left(\frac{C_{0}}{C_{3}}\right)^{2}
-\frac{2}{C_{3} m_{0}\left(x_{1}\right)}},\quad
\quad t_+:=\frac{C_{0}}{2 C_{3}}
+\frac{1}{2} \sqrt{\left(\frac{C_{0}}{C_{3}}\right)^{2}
-\frac{2}{C_{3} m_{0}\left(x_{1}\right)}}.
\end{eqnarray}

According to the condition (\ref{finalCondition}), one derives
\begin{eqnarray}
\left(\frac{C_{0}}{C_{3}}\right)^{2}
>\frac{2}{C_{3} m_{0}\left(x_{1}\right)}, \quad  \quad 0<t_-<\frac{C_{0}}{2 C_{3}}<t_+
\end{eqnarray}
and consequently
\begin{eqnarray}\label{last1}
0 \leq \frac{1}{\overline{m}(t)} \leq C_{3}\left(t-t_-\right)\left(t-t_+\right),
\end{eqnarray}
from which we can find a time $0<T_{*} \leq t_-$ such that
\begin{eqnarray*}
m(t) \longrightarrow+\infty, \quad  M(t) \longrightarrow-\infty,\quad \text { as } \quad t \longrightarrow T_{*} \leq t_-.
\end{eqnarray*}

Thus one has
\begin{eqnarray*}
\inf _{x \in \mathbb{R}} M(t, x) \leq M(t) \longrightarrow-\infty, \quad \text { as } \quad t \longrightarrow T_{*} \leq t_-,
\end{eqnarray*}
implying the blow-up of $m(t,x)$
at time $T_{*}$ since $u$ and $u_{x}$ are both bounded.

We finally prove the blow-up rate (\ref{blowupRate}).  Combining (\ref{afterIntegration1}) and (\ref{last1}) leads to
\begin{eqnarray*}
2\left(T_{*}-t\right)
\inf _{x \in \mathbb{R}} M(t, x)\leq 2\left(T_{*}-t\right) M(t)\leq-\left(T_{*}-t\right) \overline{m}(t)\left(C_{0}-2 C_{3} t\right) \leq \frac{2\left(T_{*}-t\right)}{\left(t-t_-\right)\left(t-t_+\right)}\left(t-\frac{C_{0}}{2 C_{3}}\right),
\end{eqnarray*}
which gives (\ref{blowupRate})
when $T_{*}=t_-.$ We thus complete the proof of Theorem \ref{WaveBreakingThm4}.

\section{Proof of Proposition \ref{exsitenceOfPeakon}}

Without loss of generality, one can take $x_0=0$ in the peakon solution (\ref{singlePeakon}). Let $c$ be a unknown real parameter, and the initial value be  $u_{0,c}(x)=u_{a}(0, x)$, $x\in \mathbb{R}.$ We first recall that \cite{GuiLiuOlverQu2013CMP} for all $t \geq 0$
\begin{eqnarray}\label{exsitenceOfPeakon001}
\partial_{x} u_{a}(t, x)=-\operatorname{sgn}(x-c t) u_{a}(t, x),\quad
\partial_{t} u_{a}(t, x)=c \operatorname{sgn}(x-c t) u_{a}(t, x) \in L^{\infty}(\mathbb{R}),
\end{eqnarray}
in the sense of distribution $\mathcal{S}^{\prime}(\mathbb{R})$, and
\begin{eqnarray}\label{exsitenceOfPeakon002}
\lim _{t \rightarrow 0^{+}}\left\|u_{a}(t, \cdot)-u_{0, c}(\cdot)\right\|_{W^{1, \infty}}=0
\end{eqnarray}

In light of Eqs.~(\ref{exsitenceOfPeakon001})-(\ref{exsitenceOfPeakon002}),
one obtains for any  test function $\varphi(t, x) \in C_{c}^{\infty}([0,+\infty) \times \mathbb{R})$ that
\begin{eqnarray}\label{exsitenceOfPeakon004}
&&\int_{0}^{+\infty}\!\! \int_{\mathbb{R}}\left(u_{a} \partial_{t} \varphi
+\frac{k_{1}+k_{2}}{3} u_{a}^{3} \partial_{x}\varphi+\frac{k_1}{3}(\partial_{x} u_{a})^{3} \varphi
+\frac{k_3}{2}u_{a}^{2}\varphi_{x}\right) \mbox{d} x \mbox{d} t+\int_{\mathbb{R}} u_{a}(0, x) \varphi(0, x) \mbox{d} x\nonumber \\
&&=-\int_{0}^{+\infty}\!\! \int_{\mathbb{R}} \varphi\left(\partial_{t} u_{a}+(k_{1}+k_{2})u_{a}^{2} \partial_{x} u_{a}
+k_{3}u_{a} \partial_{x} u_{a}-\frac{k_1}{3}\left(\partial_{x} u_{a}\right)^{3}\right)
\mbox{d} x \mbox{d} t \nonumber\\
&&=-\int_{0}^{+\infty}\!\! \int_{\mathbb{R}} \varphi \operatorname{sgn}(x-c t) u_{a}\left[c-
\left(\frac{2k_1}{3}+k_{2}\right) u_{a}^{2}-k_{3}u_{a}\right] \mbox{d} x \mbox{d} t.
\end{eqnarray}

We next deal with  the nonlocal terms in Eq.~(\ref{defOfWeakSolution1}) with $u_{a}(t, x)$
instead of $u(t, x)$.
\begin{eqnarray}\label{exsitenceOfPeakon005}
&&\quad\int_{0}^{+\infty} \int_{\mathbb{R}}
\bigg\{p*\left[\left(\frac{2k_1}{3}+k_{2}\right)u_{a}^{3}+\left(k_{1}+\frac{3k_2}{2}\right)u_{a}(\partial_{x}u_{a})^{2}
+k_{3}u_{a}^{2}+\frac{k_3}{2}(\partial_{x}u_{a})^{2}\right]\partial_{x} \varphi\nonumber\\
&&\quad\quad-p*\left[\left(\frac{k_1}{3}+\frac{k_2}{2}\right)(\partial_{x}u_{a})^{3}\right]\varphi \bigg\}\mbox{d} x \mbox{d} t\nonumber\\
&&=-\int_{0}^{+\infty} \int_{\mathbb{R}}
\bigg\{\varphi p_{x}*\left[\left(k_{1}+\frac{3k_2}{2}\right)u_{a}(\partial_{x}u_{a})^{2}+\left(\frac{7k_1}{9}
+\frac{7k_2}{6}\right) u_{a}^{3}\right]
\bigg\}\mbox{d} x \mbox{d} t\nonumber\\
&&\quad-\int_{0}^{+\infty} \int_{\mathbb{R}}
\bigg\{\varphi p_{x}*\left[k_{3}u_{a}^{2}+\frac{k_3}{2}(\partial_{x}u_{a})^{2}\right]
\bigg\}\mbox{d} x \mbox{d} t,
\end{eqnarray}
where the equality  $\operatorname{sgn}(x-c t)u_{a}^{3}=-\frac{1}{3}\partial_{x}(u_{a}^{3})$
has been used to deduce
\begin{eqnarray}\label{exsitenceOfPeakon006}
&&(2k_{1}\!+\!3k_{2})u_{a}^{2}\partial_{x}u_{a}\!+\!\left(\!\frac{k_1}{3}\!+\!\frac{k_2}{2}\!\right)\!(\partial_{x}u_{a})^{3}\!
=\!\left(\!\frac{7k_1}{9}\!+\!\frac{7k_2}{6}\!\right) \partial_{x}(u_{a}^{3}).\qquad\quad
\end{eqnarray}

 Invoking
$\partial_{x} p(x)=-\frac{1}{2} \operatorname{sgn}(x) e^{-|x|}$
leads to
\begin{eqnarray}\label{exsitenceOfPeakon007}
&& p_{x}*\left[\left(k_{1}+\frac{3k_2}{2}\right)u_{a}(\partial_{x}u_{a})^{2}
+\left(\frac{7k_1}{9}+\frac{7k_2}{6}\right) u_{a}^{3}+k_{3}u_{a}^{2}+\frac{k_3}{2}(\partial_{x}u_{a})^{2}\right]\nonumber\\
&&=-\frac{1}{2}\int_{-\infty}^{+\infty}
 \operatorname{sgn}(x-y) e^{-|x-y|}\Big\{a^3\left[\left(\frac{7k_1}{9}+\frac{7k_2}{6}\right)+\left(k_{1}
+\frac{3k_2}{2}\right)\operatorname{sgn^{2}}(y-ct)\right]e^{-3|y-ct|} \nonumber\\
&& \qquad +a^2\left[k_{3}+\frac{k_3}{2}\operatorname{sgn^{2}}(y-ct)\right]e^{-2|y-ct|}\Big\}
\mbox{d} y .
\end{eqnarray}

For the case $x > ct$, the right hand side of Eq.~(\ref{exsitenceOfPeakon007}) can be separated into  three parts:
\begin{eqnarray}\label{exsitenceOfPeakon009}
&&p_{x}*\left[\left(k_{1}+\frac{3k_2}{2}\right)u_{a}(\partial_{x}u_{a})^{2}+\left(\frac{7k_1}{9}+\frac{7k_2}{6}\right) u_{a}^{3}+k_{3}u_{a}^{2}+\frac{k_3}{2}(\partial_{x}u_{a})^{2}
\right]\nonumber\\
&&=-\frac{1}{2}\left(\int_{-\infty}^{ct}+\int_{ct}^{x}+\int_{x}^{+\infty}\right)
 \operatorname{sgn}(x-y) e^{-|x-y|}\nonumber\\
&&\quad\quad\times\Big\{a^3\left[\left(\frac{7k_1}{9}+\frac{7k_2}{6}\right)+\left(k_{1}+\frac{3k_2}{2}\right)\operatorname{sgn^{2}}(y-ct)
\right]e^{-3|y-ct|} \nonumber\\
&&\qquad +a^2\left[k_{3}+\frac{k_3}{2}\operatorname{sgn^{2}}(y-ct)\right]e^{-2|y-ct|}\Big\}dy
\nonumber\\
&&=a^3\left(\frac{2k_1}{3}+k_{2}\right)(e^{3(ct-x)}-e^{ct-x})+a^2k_{3}(e^{2(ct-x)}-e^{ct-x}).
\end{eqnarray}

For the case $x \leq ct$, we split the right hand side of (\ref{exsitenceOfPeakon007}) into the following three parts:
\begin{eqnarray}\label{exsitenceOfPeakon0010}
&&p_{x}*\left[\left(k_{1}+\frac{3k_2}{2}\right)u_{a}(\partial_{x}u_{a})^{2}+\left(\frac{7k_1}{9}+\frac{7k_2}{6}\right) u_{a}^{3}+k_{3}u_{a}^{2}+\frac{k_3}{2}(\partial_{x}u_{a})^{2}
\right]\nonumber\\
&&=-\frac{1}{2}\left(\int_{-\infty}^{x}+\int_{x}^{ct}+\int_{ct}^{+\infty}\right)
 \operatorname{sgn}(x-y) e^{-|x-y|}\nonumber\\
&&\quad\quad\times\Big\{a^3\left[\left(\frac{7k_1}{9}+\frac{7k_2}{6}\right)+\left(k_{1}+\frac{3k_2}{2}\right)\operatorname{sgn^{2}}(y-ct)
\right]e^{-3|y-ct|} \nonumber\\
&&\qquad +a^2\left[k_{3}+\frac{k_3}{2}\operatorname{sgn^{2}}(y-ct)\right]e^{-2|y-ct|}\Big\}dy
\nonumber\\
&&=a^3\left(\frac{2k_1}{3}+k_{2}\right)(e^{x-ct}-e^{3(x-ct)})+a^2k_{3}(e^{x-ct}-e^{2(x-ct)}).
\end{eqnarray}

Thus, combining Eqs.~(\ref{exsitenceOfPeakon009}) and (\ref{exsitenceOfPeakon0010}), and according to the definition of $u_{a}$ yield
\begin{eqnarray*}
&&\operatorname{sgn}(x-c t) u_{a}\left[c-\left(\frac{2}{3}k_{1}+k_{2}\right) u_{a}^{2}-k_{3}u_{a}\right]\nonumber\\
&&\quad +p_{x}*\left[\left(k_{1}+\frac{3k)2}{2}\right)u_{a}(\partial_{x}u_{a})^{2}+\left(\frac{7k_1}{9}+\frac{7k_2}{6}\right) u_{a}^{3}
 +k_{3}u_{a}^{2}+\frac{k_3}{2}(\partial_{x}u_{a})^{2}\right]\nonumber\\
 &&\quad=
\left\{ \begin{array}{ll}
   \left[-\left(\frac{2}{3}k_{1}+k_{2}\right)a^{3}-k_{3}a^{2}+ac\right]e^{ct-x}, & {\rm as}\,\, x>ct, \v\\
 \left[\left(\frac{2}{3}k_{1}+k_{2}\right)a^{3}+k_{3}a^{2}-ac\right]e^{x-ct}, & {\rm as}\,\, x\leq ct,
  \end{array}\right.
\end{eqnarray*}
which is equal to zero if
\begin{eqnarray}
\left(\frac{2}{3}k_{1}+k_{2}
\right)a^{3}+k_{3}a^{2}-ac=0,
\end{eqnarray}
that is, $c$ satisfies Eq.~(\ref{solu-c}).
With the condition (\ref{solu-c}) satisfied by $c$, we have arrived at
\begin{eqnarray*}
&&\int_{0}^{+\infty} \int_{\mathbb{R}}
\bigg\{u_{a} \varphi_{t}
+\frac{k_{1}+k_{2}}{3} u_{a}^{3} \varphi_{x}+\frac{k_1}{3}(\partial_{x}u_{a})^{3} \varphi+\frac{k_3}{2}u_{a}^{2}\varphi_{x}\nonumber\\
&&\quad \quad +p*\left[\left(\frac{2k_1}{3}+k_{2}\right)u_{a}^{3}
+\left(k_{1}+\frac{3k_2}{2}\right)u_{a}(\partial_{x}u_{a})^{2}+k_{3}u^{2}+\frac{k_3}{2}(\partial_{x}u_{a})^{2}
\right]\partial_{x} \varphi\nonumber\\
&&\quad\quad-p*\left[\left(\frac{k_1}{3}+\frac{k_2}{2}\right)(\partial_{x}u_{a})^{3}
\right]\varphi \bigg\}\mbox{d} x \mbox{d} t
+\int_{\mathbb{R}} u_{a}(0,x) \varphi(0, x) \mbox{d} x=0
\end{eqnarray*}
for any smooth test function $\varphi(t, x) \in C_{c}^{\infty}([0, +\infty) \times \mathbb{R}).$  We thus complete the proof of Proposition \ref{exsitenceOfPeakon}.

\section{Proof of Proposition~\ref{multipeakon}}

In the last section, we have deduced  the single peakon solution of Eq.~(\ref{gmCH-Novikov}). It is found that Eq.~(\ref{gmCH-Novikov}) also admits the multi-peakon solutions. We are now in a position to derive them.
Suppose the multi-peakon solutions of Eq.~(\ref{gmCH-Novikov}) are of the form (\ref{mp}).
Then  we have
\begin{eqnarray}\label{multipeakon2}
u_{npm,t}=\sum_{i=1}^{N}[\dot{p}_{i}(t)+p_i(t)\dot{q}_{i}(t)\operatorname{sgn}(\xi_i)]e^{-|\xi_i|},\quad
u_{npm,x}=-\sum_{i=1}^{N}\operatorname{sgn}(\xi_i)p_{i}(t)e^{-|\xi_i|}.
\end{eqnarray}

Without loss of generality, we assume that the time-dependent positions are ordered, that is, $q_1(t)<q_2(t)<\cdots<q_{N-1}(t)<q_N(t)$. Since Eq.~(\ref{gmCH-Novikov}) has the non-periodic multi-peakon solution (\ref{mp}) in the sense of Definition \ref{defOfWeakSolution}, the multi-peakon solution $u_{npm}(x,t)$ given by Eq.~(\ref{mp}) satisfies ($m=u-u_{xx}$ and $u=u_{npm}(x,t)$) the following equation
\begin{eqnarray}\label{multi}
&&\int_{0}^{\infty} \int_{\mathbb{R}}
\bigg\{u_{t} (\psi-\psi_{xx})
+\frac{1}{3}(k_{1}+k_{2}) u^{3} \psi_{xxx}
+\frac{k_{1}}{3}u_{x}^{3}\psi_{xx}
-\left(k_{1}+\frac{4}{3}k_{2}\right)u^{3}\psi_{x}\nonumber\\
&&\qquad-\left(k_{1}+\frac{3}{2}k_{2}\right)u u_{x}^{2}\psi_{x}
+\frac{1}{2}k_{3}u^{2}\psi_{xxx}
-\frac{3}{2}k_{3}u^{2}\psi_{x}
-\frac{1}{2}k_{3}u_{x}^{2}\psi_{x}
+\frac{1}{2}k_{2}u_{x}^{3}\psi
 \bigg\}\mbox{d} x \mbox{d} t=0\qquad
\end{eqnarray}
after substituting $\psi=(1-\partial_{x}^{2})^{-1}\varphi$ in Eq.~(\ref{defOfWeakSolution1}) for $\varphi(t, x) \in C_{c}^{\infty}([0, \infty) \times \mathbb{R}).$

Note that it has been calculated  in \cite{GuiLiuOlverQu2013CMP}
that
\begin{eqnarray*}
\int_{0}^{\infty} \int_{\mathbb{R}}
u_{t} (\psi-\psi_{xx})\mbox{d} x \mbox{d} t
=2\int_{0}^{\infty}\sum_{j=1}^{N}
(p_{j}^{\prime}\psi(q_{j})+p_{j}q_{j}^{\prime}\psi_{x}(q_{j}))\mbox{d} t.
\end{eqnarray*}
We just need to  deal with the remaining terms in (\ref{multi}).
For the term
$\int_{0}^{\infty} \int_{\mathbb{R}}
u^{3} \psi_{xxx}
\mbox{d} x\mbox{d} t$,
we have
\begin{eqnarray*}
\int_{\mathbb{R}}
u^{3} \psi_{xxx}
\mbox{d} x
=\sum_{i,j,k=1}^{N}p_{i}p_{j}p_{k}
\int_{\mathbb{R}}e^{-|x-q_{i}|-|x-q_{j}|-|x-q_{k}|}\psi_{xxx}\mbox{d} x
=K_1+K_2+K_3,
\end{eqnarray*}
where
\begin{align*}
K_1=&\displaystyle\sum_{i=1}^{N}p_{i}^{3}\int_{\mathbb{R}}e^{-3|x-q_{i}|}\psi_{xxx}\mbox{d} x, \v\\
=& \displaystyle-6\sum_{i=1}^{N}p_{i}^{3}\psi_{x}(q_{i})
-27\sum_{i=1}^{N}p_{i}^{3}\left(\int_{-\infty}^{q_{i}}\psi e^{3x-3q_{i}}\mbox{d}x
 -\int_{q_{i}}^{+\infty}\psi e^{-3x+3q_{i}}\mbox{d} x\right), \\
K_2=&\displaystyle 3\sum_{i=1}^{N}p_i^2\left(\sum_{k>i}p_{k}\int_{\mathbb{R}}e^{-2|x-q_{i}|-|x-q_{k}|}\psi_{xxx}dx
+\sum_{k<i}p_{k}\int_{\mathbb{R}}e^{-2|x-q_{i}|-|x-q_{k}|}\psi_{xxx}dx\right)\v\\
=&\displaystyle 3\sum_{i=1}^{N}\sum_{k>i}p_{i}^{2}p_{k}\bigg\{
4[2\psi(q_{i})-\psi_{x}(q_{i})]e^{q_{i}-q_{k}}-2[\psi_{x}(q_{k})+4\psi(q_{k})]e^{2q_{i}-2q_{k}} \v\\
&\displaystyle -27\int_{-\infty}^{q_{i}}\psi e^{3x-2q_{i}-q_{k}}dx
+\int_{q_{i}}^{q_{k}}\psi e^{-x+2q_{i}-q_{k}}dx+27\int_{q_{k}}^{+\infty}\psi e^{-3x+2q_{i}+q_{k}}dx\bigg\} \v\\
&\displaystyle +3\sum_{i=1}^{N}\sum_{k<i}p_{i}^{2}p_{k}\bigg\{
-4[\psi_{x}(q_{i})+2\psi(q_{i})]e^{q_{k}-q_{i}}-2[\psi_{x}(q_{k})-4\psi(q_{k})]e^{2q_{k}-2q_{i}} \v\\
&\displaystyle -27\int_{-\infty}^{q_{k}}\psi e^{3x-2q_{i}-q_{k}}dx
-\int_{q_{k}}^{q_{i}}\psi e^{x-2q_{i}+q_{k}}dx
+27\int_{q_{i}}^{+\infty}\psi e^{-3x+2q_{i}+q_{k}}dx \bigg\}, \\
K_3=&\displaystyle 6\sum_{i<j<k}p_{i}p_{j}p_{k}\int_{\mathbb{R}}e^{-|x-q_{i}|-|x-q_{j}|-|x-q_{k}|}\psi_{xxx}dx \v\\
&\displaystyle
+\int_{q_{j}}^{q_{k}}e^{q_{i}+q_{j}-q_{k}-x}\psi_{xxx}\mbox{d} x
+\int_{q_{k}}^{+\infty}e^{q_{i}+q_{j}+q_{k}-3x}\psi_{xxx}\mbox{d} x
\bigg\}\v\\
=&\displaystyle6\sum_{i<j<k}p_{i}p_{j}p_{k}\bigg\{
-2\psi_{x}(q_{i})e^{2q_{i}-q_{j}-q_{k}}
-2\psi_{x}(q_{j})e^{q_{i}-q_{k}}
-2\psi_{x}(q_{k})e^{q_{i}+q_{j}-2q_{k}}\v\\
&\displaystyle +8\psi(q_{i})e^{2q_{i}-q_{j}-q_{k}}
-8\psi(q_{k})e^{q_{i}+q_{j}-2q_{k}}
-27\int_{-\infty}^{q_{i}}\psi e^{3x-q_{i}-q_{j}-q_{k}}\mbox{d} x
-\int_{q_{i}}^{q_{j}}\psi e^{x+q_{i}-q_{j}-q_{k}}\mbox{d} x\v\\
&\displaystyle
+\int_{q_{j}}^{q_{k}}\psi e^{-x+q_{i}+q_{j}-q_{k}}\mbox{d} x
+27\int_{q_{k}}^{+\infty}\psi e^{-3x+q_{i}+q_{j}+q_{k}}\mbox{d} x
\bigg\}.
\end{align*}

Applying  similar method to the remaining terms in (\ref{multi}) produces
\begin{align}
\int_{\mathbb{R}}u_{x}^{3}\psi_{xx}\mbox{d}x
=&-\sum_{i=1}^{N}p_{i}^{3}\bigg\{
-2\psi_{x}(q_{i})-9\int_{-\infty}^{q_{i}}\psi e^{3x-3q_{i}}dx+9\int_{q_{i}}^{+\infty}\psi e^{-3x+3q_{i}}\mbox{d}x
\bigg\}\qquad\qquad\no\\
& -3\sum_{i=1}^{N}\sum_{k>i}p_{i}^{2}p_{k}\bigg\{
-2[\psi_{x}(q_{k})+2\psi(q_{k})]e^{2q_{i}-2q_{k}}+4\psi(q_{i})e^{q_{i}-q_{k}}\no\\
&\qquad-9\int_{-\infty}^{q_{i}}\psi e^{3x-2q_{i}-q_{k}}\mbox{d}x
-\int_{q_{i}}^{q_{k}}\psi e^{-x+2q_{i}-q_{k}}\mbox{d}x
+9\int_{q_{k}}^{+\infty}\psi e^{-3x+2q_{i}+q_{k}}\mbox{d}x\bigg\}\no\\
&-3\sum_{i=1}^{N}\sum_{k<i}p_{i}^{2}p_{k}\bigg\{-2[\psi_{x}(q_{k})-2\psi(q_{k})]e^{2q_{k}-2q_{i}}
-4\psi(q_{i})e^{q_{k}-q_{i}}\no\\
&\qquad
-9\int_{-\infty}^{q_{k}}\psi e^{3x-2q_{i}-q_{k}}\mbox{d}x
+\int_{q_{k}}^{q_{i}}\psi e^{x-2q_{i}+q_{k}}\mbox{d}x
+9\int_{q_{i}}^{+\infty}\psi e^{-3x+2q_{i}+q_{k}}\mbox{d}x\bigg\}\no\\
&-6\sum_{i<j<k}p_{i}p_{j}p_{k}\bigg\{
-2\psi_{x}(q_{i})e^{2q_{i}-q_{j}-q_{k}}
+2\psi_{x}(q_{j})e^{q_{i}-q_{k}}
-2\psi_{x}(q_{k})e^{q_{i}+q_{j}-2q_{k}}\no\\
&\qquad
+4\psi(q_{i})e^{2q_{i}-q_{j}-q_{k}}
-4\psi(q_{k})e^{q_{i}+q_{j}-2q_{k}}
-9\int_{-\infty}^{q_{i}}\psi e^{3x-q_{i}-q_{j}-q_{k}}\mbox{d} x \no\\
&\qquad+\int_{q_{i}}^{q_{j}}\psi e^{x+q_{i}-q_{j}-q_{k}}\mbox{d} x
-\int_{q_{j}}^{q_{k}}\psi e^{-x+q_{i}+q_{j}-q_{k}}\mbox{d} x
+9\int_{q_{k}}^{+\infty}\psi e^{-3x+q_{i}+q_{j}+q_{k}}\mbox{d} x
\bigg\},\qquad
\end{align}
\begin{align}
\int_{\mathbb{R}}& u^{3}\psi_{x}dx
=\sum_{i=1}^{N}p_{i}^{3}\bigg\{
-3\int_{-\infty}^{q_{i}}\psi e^{3x-3q_{i}}\mbox{d}x+3\int_{q_{i}}^{+\infty}\psi e^{-3x+3q_{i}}\mbox{d}x
\bigg\} \no\\
&+3\sum_{i=1}^{N}\sum_{k>i}p_{i}^{2}p_{k}\bigg\{
-3\int_{-\infty}^{q_{i}}\psi e^{3x-2q_{i}-q_{k}}\mbox{d}x
+\int_{q_{i}}^{q_{k}}\psi e^{-x+2q_{i}-q_{k}}\mbox{d}x
+3\int_{q_{k}}^{+\infty}\psi e^{-3x+2q_{i}+q_{k}}\mbox{d}x
\bigg\}\no\\
&+3\sum_{i=1}^{N}\sum_{k<i}p_{i}^{2}p_{k}\bigg\{
-3\int_{-\infty}^{q_{k}}\psi e^{3x-2q_{i}-q_{k}}\mbox{d}x
-\int_{q_{k}}^{q_{i}}\psi e^{x-2q_{i}+q_{k}}\mbox{d}x
+3\int_{q_{i}}^{+\infty}\psi e^{-3x+2q_{i}+q_{k}}\mbox{d}x
\bigg\}\no\\
&+6\sum_{i<j<k}p_{i}p_{j}p_{k}\bigg\{
-3\int_{-\infty}^{q_{i}}\psi e^{3x-q_{i}-q_{j}-q_{k}}\mbox{d} x
-\int_{q_{i}}^{q_{j}}\psi e^{x+q_{i}-q_{j}-q_{k}}\mbox{d} x\no\\
&\quad\quad\quad\quad\quad\quad\quad+\int_{q_{j}}^{q_{k}}\psi e^{-x+q_{i}+q_{j}-q_{k}}\mbox{d} x
+3\int_{q_{k}}^{+\infty}\psi e^{-3x+q_{i}+q_{j}+q_{k}}\mbox{d} x
\bigg\},
\end{align}
\begin{align}
\int_{\mathbb{R}}u u_{x}^{2}\psi_{x}\mbox{d}x
=&\sum_{i=1}^{N}p_{i}^{3}\bigg\{
-3\int_{-\infty}^{q_{i}}\psi e^{3x}e^{-3q_{i}}\mbox{d}x
+3\int_{q_{i}}^{+\infty}\psi e^{-3x}e^{3q_{i}}\mbox{d}x
\bigg\}\no\\
&+2\sum_{i=1}^{N}\sum_{k>i}p_{i}^{2}p_{k}\bigg\{
2\psi(q_{i})e^{q_{i}-q_{k}}
-2\psi(q_{k})e^{2q_{i}-2q_{k}}-3\int_{-\infty}^{q_{i}}\psi e^{3x}e^{-2q_{i}-q_{k}}\mbox{d}x\no\\
&\quad\quad
-\int_{q_{i}}^{q_{k}}\psi e^{-x}e^{2q_{i}-q_{k}}\mbox{d}x
+3\int_{q_{k}}^{+\infty}\psi e^{-3x}e^{2q_{i}+q_{k}}\mbox{d}x
\bigg\}\no\\
&+\sum_{i=1}^{N}\!\sum_{k>i}p_{i}^{2}p_{k}\!\bigg\{\!-\!3\int_{-\infty}^{q_{i}}\psi e^{3x-2q_{i}-q_{k}}\mbox{d}x
\!+\!\int_{q_{i}}^{q_{k}}\psi e^{-x+2q_{i}-q_{k}}\mbox{d}x\!+\!
3\int_{q_{k}}^{+\infty}\psi e^{-3x+2q_{i}+q_{k}}\mbox{d}x
\bigg\}\no\\
&+2\sum_{i=1}^{N}\sum_{k<i}p_{i}^{2}p_{k}\bigg\{
-2\psi(q_{i})e^{q_{k}-q_{i}}
+2\psi(q_{k})e^{2q_{k}-2q_{i}}-3\int_{-\infty}^{q_{k}}\psi e^{3x}e^{-2q_{i}-q_{k}}\mbox{d}x\no\\
&\quad\quad
+\int_{q_{k}}^{q_{i}}\psi e^{x}e^{-2q_{i}+q_{k}}\mbox{d}x
+3\int_{q_{i}}^{+\infty}\psi e^{-3x}e^{2q_{i}+q_{k}}\mbox{d}x
\bigg\}\no\\
&+\sum_{i=1}^{N}\sum_{k<i}p_{i}^{2}p_{k}\!\bigg\{\!-\!3\int_{-\infty}^{q_{k}}\psi e^{3x-2q_{i}-q_{k}}\mbox{d}x
\!-\!\int_{q_{k}}^{q_{i}}\psi e^{x-2q_{i}+q_{k}}\mbox{d}x
\!+\!3\int_{q_{i}}^{+\infty}\psi e^{-3x+2q_{i}+q_{k}}\mbox{d}x
\bigg\}\no\\
&+2\sum_{i<j<k}p_{i}p_{j}p_{k}\bigg\{
4\psi(q_{i})e^{2q_{i}-q_{j}-q_{k}}
-4\psi(q_{k})e^{q_{i}+q_{j}-2q_{k}}
-9\int_{-\infty}^{q_{i}}\psi e^{3x}e^{-q_{i}-q_{j}-q_{k}}\mbox{d} x\no\\
&\quad\quad\quad
+\int_{q_{i}}^{q_{j}}\psi e^{x}e^{q_{i}-q_{j}-q_{k}}\mbox{d} x
-\int_{q_{j}}^{q_{k}}\psi e^{-x}e^{q_{i}+q_{j}-q_{k}}\mbox{d} x
+9\int_{q_{k}}^{+\infty}\psi e^{-3x}e^{q_{i}+q_{j}+q_{k}}\mbox{d} x
\bigg\},
\end{align}
\begin{align}
\int_{\mathbb{R}}u^{2}\psi_{xxx}\mbox{d}x
=&\sum_{i=1}^{N}p_{i}^{2}\bigg\{
-4\psi_{x}(q_{i})
-8\int_{-\infty}^{q_{i}}\psi e^{2x}e^{-2q_{i}}\mbox{d}x
+8\int_{q_{i}}^{+\infty}\psi e^{-2x}e^{2q_{i}}\mbox{d}x
\bigg\}\no\\
&+\sum_{i=1}^{N}\sum_{j>i}p_{i}p_{j}\bigg\{
-2\psi_{x}(q_{i})e^{q_{i}-q_{j}}
-2\psi_{x}(q_{j})e^{q_{i}-q_{j}}
+4\psi(q_{i})e^{q_{i}-q_{j}}
-4\psi(q_{j})e^{q_{i}-q_{j}}\no\\
&-8\int_{-\infty}^{q_{i}}\psi e^{2x}e^{-q_{i}-q_{j}}\mbox{d}x
+8\int_{q_{j}}^{+\infty}\psi e^{-2x}e^{q_{i}+q_{j}}\mbox{d}x
\bigg\}\no\\
&+\sum_{i=1}^{N}\sum_{j<i}p_{i}p_{j}\bigg\{
-2\psi_{x}(q_{i})e^{q_{j}-q_{i}}
-2\psi_{x}(q_{j})e^{q_{j}-q_{i}}
-4\psi(q_{i})e^{q_{j}-q_{i}}
+4\psi(q_{j})e^{q_{j}-q_{i}}\no\\
&-8\int_{-\infty}^{q_{j}}\psi e^{2x}e^{-q_{i}-q_{j}}\mbox{d}x
+8\int_{q_{i}}^{+\infty}\psi e^{-2x}e^{q_{i}+q_{j}}\mbox{d}x
\bigg\},
\end{align}
\begin{align}
\int_{\mathbb{R}}u^{2}\psi_{x}\mbox{d}x
=&2\sum_{i=1}^{N}p_{i}^{2}\bigg\{
-\int_{-\infty}^{q_{i}}\psi e^{2x}e^{-2q_{i}}\mbox{d}x
+\int_{q_{i}}^{+\infty}\psi e^{-2x}e^{2q_{i}}\mbox{d}x
\bigg\}\no\\
&\quad
+2\sum_{i=1}^{N}\sum_{j>i}p_{i}p_{j}\bigg\{
-\int_{-\infty}^{q_{i}}\psi e^{2x}e^{-q_{i}-q_{j}}\mbox{d}x
+\int_{q_{j}}^{+\infty}\psi e^{-2x}e^{q_{i}+q_{j}}\mbox{d}x
\bigg\}\no\\
&\quad
+2\sum_{i=1}^{N}\sum_{j<i}p_{i}p_{j}\bigg\{
-\int_{-\infty}^{q_{j}}\psi e^{2x}e^{-q_{i}-q_{j}}\mbox{d}x
+\int_{q_{i}}^{+\infty}\psi e^{-2x}e^{q_{i}+q_{j}}\mbox{d}x
\bigg\},
\end{align}
\begin{align}
\int_{\mathbb{R}}u_{x}^{2}\psi_{x}\mbox{d}x
=&2\sum_{i=1}^{N}p_{i}^{2}\bigg\{-\int_{-\infty}^{q_{i}}\psi e^{2x}e^{-2q_{i}}\mbox{d}x
+\int_{q_{i}}^{+\infty}\psi e^{-2x}e^{2q_{i}}\mbox{d}x\bigg\}\qquad\qquad\qquad \no\\
&+2\sum_{i=1}^{N}\sum_{j>i}p_{i}p_{j}\bigg\{[\psi(q_{i})-\psi(q_{j})]e^{q_{i}-q_{j}}
-\int_{-\infty}^{q_{i}}\psi e^{2x}e^{-q_{i}-q_{j}}\mbox{d}x
+\int_{q_{j}}^{+\infty}\psi e^{-2x}e^{q_{i}+q_{j}}\mbox{d}x
\bigg\}\no\\
&+2\sum_{i=1}^{N}\sum_{j<i}p_{i}p_{j}\bigg\{[\psi(q_{j})-\psi(q_{i})]e^{q_{j}-q_{i}}
-\int_{-\infty}^{q_{j}}\psi e^{2x}e^{-q_{i}-q_{j}}\mbox{d}x
+\int_{q_{i}}^{+\infty}\psi e^{-2x}e^{q_{i}+q_{j}}\mbox{d}x
\bigg\}
\end{align}
and
\begin{align}
\int_{\mathbb{R}}u_{x}^{3}\psi\mbox{d}x
=&\sum_{i=1}^{N}p_{i}^{3}\bigg\{
\int_{-\infty}^{q_{i}}\psi e^{3x}e^{-3q_{i}}\mbox{d}x
-\int_{q_{i}}^{+\infty}\psi e^{-3x}e^{3q_{i}}\mbox{d}x
\bigg\}\qquad\qquad\qquad\qquad\qquad\qquad\no\\
&-3\sum_{i=1}^{N}\sum_{k>i}p_{i}^{2}p_{k}\bigg\{
-\int_{-\infty}^{q_{i}}\psi e^{3x}e^{-2q_{i}-q_{k}}\mbox{d}x
-\int_{q_{i}}^{q_{k}}\psi e^{-x}e^{2q_{i}-q_{k}}\mbox{d}x
+\int_{q_{k}}^{+\infty}\psi e^{-3x}e^{2q_{i}+q_{k}}\mbox{d}x
\bigg\}\no\\
&-3\sum_{i=1}^{N}\sum_{k<i}p_{i}^{2}p_{k}\bigg\{
-\int_{-\infty}^{q_{k}}\psi e^{3x}e^{-2q_{i}-q_{k}}\mbox{d}x
+\int_{q_{k}}^{q_{i}}\psi e^{x}e^{-2q_{i}+q_{k}}\mbox{d}x
+\int_{q_{i}}^{+\infty}\psi e^{-3x}e^{2q_{i}+q_{k}}\mbox{d}x
\bigg\}\no\\
&-6\sum_{i<j<k}p_{i}p_{j}p_{k}\bigg\{
-\int_{-\infty}^{q_{i}}\psi e^{3x}e^{-q_{i}-q_{j}-q_{k}}\mbox{d}x
+\int_{q_{i}}^{q_{j}}\psi e^{x}e^{q_{i}-q_{j}-q_{k}}\mbox{d}x\no\\
&\quad
-\int_{q_{j}}^{q_{k}}\psi e^{-x}e^{q_{i}+q_{j}-q_{k}}\mbox{d}x
+\int_{q_{k}}^{+\infty}\psi e^{-3x}e^{q_{i}+q_{j}+q_{k}}\mbox{d}x
\bigg\}.
\end{align}

Substituting all the above-founded expressions into Eq.~(\ref{multi}), and setting the
coefficients of $\psi(q_{i})$  and  $\psi_{x}(q_{i})$
both to be zero, we finds the dynamical system (\ref{pDot}) satisfied by $p_{i}(t)$
and $q_{i}(t)$ for $i=1,2,\cdots,N.$ This completes the proof of this Proposition.

\section{Proof of Proposition \ref{exsitenceOfPeriodicPeakon}}

In the periodic case, we have
\begin{eqnarray}
u=(1-\partial_{x}^{2})m=G*m,\qquad G(x)=\frac12{\rm csch}(1/2)\cosh\left(x-1/2\right).
\end{eqnarray}

We will show the relation (\ref{c-con}) between the two parameters $a$ and $c.$ Similar to Ref.~\cite{QuLiuLiu2013CMP} for the mCH equation, let  $\mathbb{S}=[0,1)$, and $u_{c}$ be  periodic functions on $\mathbb{S}$ with period one. Notice that $u_{c}$ is continuous on $\mathbb{S}$ with peak at $x=0$.
Furthermore, $u_{c}$ is smooth on $(0,1)$
and for all $t \in \mathbb{R}^{+}$
there holds
\begin{eqnarray}\label{periodic1}
\partial_{x} u_{c}(x, t)=-a \sinh(\xi)
\in L^{\infty}(\mathbb{S})
\end{eqnarray}
in the sense of periodic distribution $\mathcal{P}^{\prime},$ where $\mathcal{P}^{\prime}$ is the dual space
of the space of $C_{c}^{\infty}$ functions on $\mathbb{S}$.

Denote by $u_{c, 0}(x)=u_{c}(x, 0), x \in \mathbb{S}.$
  Then we have
\begin{eqnarray}\label{periodic2}
\lim _{t \rightarrow 0^{+}}
\left\|u_{c}(\cdot, t)-u_{c, 0}(\cdot)\right\|_{W^{1,\infty}(\mathbb{S})}=0.
\end{eqnarray}
Also, one finds
\begin{eqnarray}\label{periodic3}
\partial_{t} u_{c}(x, t)=a c \sinh(\xi) \in L^{\infty}(\mathbb{S}), \quad (t \geq 0)
\end{eqnarray}
and
\begin{eqnarray}\label{periodic4}
u_{c}^{2} \partial_{x} u_{c}=-a^{3} [\sinh(\xi)+ \sinh^{3}(\xi)],\quad
u_{c} \partial_{x} u_{c}=-\frac{1}{2}a^{2}\sinh(2\xi).
\end{eqnarray}

Similar to the case of the
peakon solution presented in Section 10, we firstly consider the local terms in Eq.~(\ref{defOfWeakSolution1}) with $u_{c}$ instead of $u$ and
derive  by invoking Eqs.~(\ref{periodic1})-(\ref{periodic4})
that
\begin{eqnarray}
&&\int_{0}^{+\infty}
 \int_{\mathbb{S}}\left(u_{c} \partial_{t} \varphi
+\frac{k_{1}+k_{2}}{3}u_{c}^{3} \partial_{x}\varphi
+\frac{k_1}{3}(\partial_{x} u_{c})^{3} \varphi+\frac{k_3}{2}u_{c}^{2}\partial_{x}\varphi\right)
\mbox{d} x \mbox{d} t+\int_{\mathbb{S}} u_{c}(0, x) \varphi(0, x) \mbox{d} x\nonumber \\
&&\quad =-\int_{0}^{+\infty}\!\!\!\int_{\mathbb{S}}\varphi\!\sinh(\xi)\bigg[
ac-(k_{1}+k_{2})a^{3}\!-\!\left(\frac{2k_1}{3}+k_{2}\right)a^{3} \sinh^2(\xi)\!
 -\!k_3a^{2}\cosh(\xi)\bigg]
  \mbox{d} x \mbox{d} t\label{periodic5}
\end{eqnarray}
for $\varphi(t, x) \in C_{c}^{\infty}([0,+\infty) \times \mathbb{S})$.

We next deal with  the nonlocal terms in Eq.~(\ref{defOfWeakSolution1}) with $u_{c}$ and $G(x)$
instead of $u$ and $p(x)$, respectively.
We derive
\begin{eqnarray}\label{periodic6}
&&\int_{0}^{+\infty} \int_{\mathbb{S}}
\bigg\{
G*\bigg[\left(\frac{2}{3}k_{1}+k_{2}\right)u_{c}^{3}
+\left(k_{1}+\frac{3}{2}k_{2}\right)u_{c}(\partial_{x}u_{c})^{2}
+k_{3}u_{c}^{2}+\frac{1}{2}k_{3}(\partial_{x}u_{c})^{2}
\bigg]\partial_{x} \varphi\nonumber\\
&&\quad\quad\quad\quad-G*\left[\left(\frac{1}{3}k_{1}
+\frac{1}{2}k_{2}\right)(\partial_{x}u_{c})^{3}\right]\varphi \bigg\}\mbox{d} x \mbox{d} t\nonumber\\
&&=-\int_{0}^{+\infty} \int_{\mathbb{S}}
\varphi G_{x}*\left[\left(k_{1}+\frac{3}{2}k_{2}\right)a^{3}\cosh(\xi)
\sinh^{2}(\xi)\right]\mbox{d} x \mbox{d} t\nonumber\\
&&\quad-\int_{0}^{+\infty} \int_{\mathbb{S}}\varphi G_{x}*\left[k_{3}a^{2}\cosh^{2}(\xi)
+\frac{1}{2}k_{3}a^{2}\sinh^{2}(\xi)\right]\mbox{d} x \mbox{d} t\nonumber\\
&&\quad+\int_{0}^{+\infty} \int_{\mathbb{S}}\varphi G*\left[(2k_{1}+3k_{2})a^{3}\sinh(\xi)
+\left(\frac{7}{3}k_{1}
+\frac{7}{2}k_{2}\right)a^{3}\sinh^{3}(\xi)
\right]
\mbox{d} x \mbox{d} t.
\end{eqnarray}
The use of $\partial_{x} G(x)=1/2{\rm csch}(1/2)\sinh(x-1 / 2-[x])$, $\zeta(x,y)=0.5-(x-y)+[x-y]$ and $\eta(y,t)=1/2-(y-c t)+[y-c t]$ gives
\begin{eqnarray}\label{periodic10}
&&(\ref{periodic6})=\int_{0}^{+\infty} \int_{\mathbb{S}}
\varphi\left[\frac{(2k_1+3k_2)a^3}{4\sinh(0.5)} \int_{\mathbb{S}}
\sinh(\zeta)\cosh(\eta)\sinh^2(\eta)dy\right]dxdt
\nonumber\\
&&+ \int_{0}^{+\infty} \int_{\mathbb{S}}
\varphi\left[\frac{k_3a^2}{8\sinh(0.5)} \int_{\mathbb{S}} \sinh(\zeta) \bigg(1+3\cosh(2\eta)\bigg)dy\right]dxdt\nonumber\\
&&+ \int_{0}^{+\infty} \int_{\mathbb{S}}
\varphi\left\{\frac{a^3}{2\sinh(0.5)}  \int_{\mathbb{S}} \cosh(\zeta)
\bigg[\frac{2k_{1}+3k_2}{8}\sinh(\eta)+7\left(\frac{k_1}{12}+\frac{k_2}{8}\right)\sinh(3\eta)\bigg]
dy\right\}dxdt \quad\quad \nonumber\\
&&=\int_{0}^{+\infty} \int_{\mathbb{S}}\varphi \sinh(\xi)\left[k_3a^2(\cosh(0.5)-\cosh(\xi))\!+
\!\left(\frac{2k_1}{3}\!+\! k_{2}\right)a^3(\sinh^{2}(0.5)\!-\!\sinh^2(\xi))\right]dxdt.\quad\qquad
\end{eqnarray}

Using (\ref{periodic5}) and (\ref{periodic10}),
one derives
\begin{eqnarray}
&&\int_{0}^{+\infty}\int_{\mathbb{S}}\left(u_{c} \partial_{t} \varphi+\frac{1}{3}(k_{1}+k_{2})u_{c}^{3} \partial_{x}\varphi
+\frac{1}{3}k_{1}(\partial_{x} u_{c})^{3} \varphi+\frac{1}{2}k_{3}u_{c}^{2}\partial_{x}\varphi\right)
\mbox{d} x \mbox{d} t+\int_{\mathbb{S}} u_{c}(0, x) \varphi(0, x) \mbox{d} x\nonumber \\
&&+\int_{0}^{+\infty} \int_{\mathbb{S}}
\bigg\{
G*\bigg[\left(\frac{2}{3}k_{1}+k_{2}\right)u_{c}^{3}
+\left(k_{1}+\frac{3}{2}k_{2}\right)u_{c}(\partial_{x}u_{c})^{2}
+k_{3}u_{c}^{2}+\frac{1}{2}k_{3}(\partial_{x}u_{c})^{2}
\bigg]\partial_{x} \varphi\nonumber\\
&&\quad\quad\quad\quad-G*\left[\left(\frac{1}{3}k_{1}
+\frac{1}{2}k_{2}\right)(\partial_{x}u_{c})^{3}\right]\varphi \bigg\}\mbox{d} x \mbox{d} t\nonumber\\
&&=\int_{0}^{+\infty}\!\!\! \int_{\mathbb{S}}\varphi\!\sinh(\xi)\bigg[(k_{1}+k_{2})a^{3}-ac
+k_{3}a^{2}\cosh(0.5)+\left(\frac{2k_1}{3}\!+\!k_{2}\right)a^{3}\sinh^{2}(0.5)]dxdt=0
\end{eqnarray}
under the condition (\ref{c-con}) satisfied by $c$.
We thus complete the proof of Proposition
\ref{exsitenceOfPeriodicPeakon}.

\section{Proof of Proposition \ref{mppeakon}}

The periodic multi-peakon solutions (\ref{mpp}) will be considered in this subsection.
In this case, $u=u_{pm}(x,t)$ given by Eq.~(\ref{mpp}) should satisfy
\begin{eqnarray}\label{multi-p}
&&\int_{0}^{\infty} \int_{0}^{1}
\bigg\{u_{t} (\psi-\psi_{xx})
+\frac{k_{1}+k_{2}}{3} u^{3} \psi_{xxx}
+\frac{k_{1}}{3}u_{x}^{3}\psi_{xx}
-\left(k_{1}+\frac{4k_2}{3}\right)u^{3}\psi_{x}\nonumber\\
&&-\left(k_{1}+\frac{3k_2}{2}\right)u u_{x}^{2}\psi_{x}
+\frac{k_3}{2}u^{2}\psi_{xxx}
-\frac{3k_3}{2}u^{2}\psi_{x}
-\frac{k_3}{2}u_{x}^{2}\psi_{x}
+\frac{k_2}{2}u_{x}^{3}\psi
 \bigg\}\mbox{d} x \mbox{d} t
=0
\end{eqnarray}
 for $\psi \in C_{c}^{\infty}([0, \infty) \times (0,1)).$

Direct computation leads to
\begin{eqnarray}\label{ut-p}
\int_{0}^{\infty} \int_{0}^{1}
u_{t} (\psi-\psi_{xx})\mbox{d} x \mbox{d} t
=2\operatorname{sinh}(1/2)\int_{0}^{\infty}\sum_{j=1}^{N}\left[\dot{p}_{j}\psi(q_{j})+p_{j}\dot{q}_{j}\psi_{x}(q_{j})\right]dt.
\end{eqnarray}

The remaining terms in (\ref{multi-p})
are calculated by using integration by parts,
we list the results as follows:
\begin{align}\label{p1}
\int_{0}^{1}u^{3}\psi_{xxx}\mbox{d}x
=&\frac{3}{4}\sum_{i=1}^{N}p_{i}^{3}\bigg\{-2\psi_{x}(q_{i})\operatorname{sinh}(3/2)
-2\psi_{x}(q_{i})\operatorname{sinh}(1/2) \no\\
&\quad\quad
+\int_{0}^{q_{i}}\psi [9\operatorname{sinh}(-3/2-3x+3q_{i})+\operatorname{sinh}(-1/2-x+q_{i})]dx \no\\
&\quad\quad
+\int_{q_{i}}^{1}\psi[9 \operatorname{sinh}(3/2-3x+3q_{i})+\operatorname{sinh}(1/2-x+q_{i})]dx\bigg\} \no\\
&+\frac{3}{4}\sum_{i=1}^{N}\sum_{k>i}p_{i}^{2}p_{k}\bigg\{
-4\psi_{x}(q_{i})[\operatorname{sinh}(3/2+q_{i}-q_{k})-\operatorname{sinh}(-1/2+q_{i}-q_{k})] \no\\
&\quad\quad
+2\psi_{x}(q_{k})[\operatorname{sinh}(1/2+2q_{i}-2q_{k})-\operatorname{sinh}(3/2+2q_{i}-2q_{k})] \no\\
&\quad\quad
+8\psi(q_{i})[\operatorname{cosh}(3/2+q_{i}-q_{k})-\operatorname{cosh}(1/2-q_{i}+q_{k})] \no\\
&\quad\quad
+8\psi(q_{k})[\operatorname{cosh}(1/2+2q_{i}-2q_{k})-\operatorname{cosh}(3/2+2q_{i}-2q_{k})] \no\\
&\quad\quad
+\int_{0}^{q_{i}}\psi [27\operatorname{sinh}(-3/2-3x+2q_{i}+q_{k})+\operatorname{sinh}(-1/2-x+2q_{i}-q_{k})]dx \no\\
&\quad\quad
+\int_{q_{i}}^{q_{k}}\psi[27 \operatorname{sinh}(1/2-3x+2q_{i}+q_{k})+\operatorname{sinh}(3/2-x+2q_{i}-q_{k})]dx \no\\
&\quad\quad
+\int_{q_{k}}^{1}\psi [27\operatorname{sinh}(3/2-3x+2q_{i}+q_{k})+ \operatorname{sinh}(1/2-x+2q_{i}-q_{k})]dx
\bigg\} \no\\
&+\frac{3}{2}\sum_{i=1}^{N}\sum_{k>i}p_{i}^{2}p_{k}\bigg\{
-2\psi_{x}(q_{k})\operatorname{sinh}(1/2)
+\int_{0}^{q_{i}}\psi \operatorname{sinh}(-1/2-x+q_{k})\mbox{d}x \no\\
&\quad\quad
+\int_{q_{i}}^{q_{k}}\psi \operatorname{sinh}(-1/2-x+q_{k})\mbox{d}x
+\int_{q_{k}}^{1}\psi \operatorname{sinh}(1/2-x+q_{k})\mbox{d}x
\bigg\}\no\\
&+\frac{3}{4}\sum_{i=1}^{N}\sum_{k<i}p_{i}^{2}p_{k}\bigg\{
4\psi_{x}(q_{i})[\operatorname{sinh}(-3/2+q_{i}-q_{k})+\operatorname{sinh}(-1/2-q_{i}+q_{k})]\no\\
&\quad\quad
-2\psi_{x}(q_{k})[\operatorname{sinh}(-1/2+2q_{i}-2q_{k})-\operatorname{sinh}(-3/2+2q_{i}-2q_{k})]\no\\
&\quad\quad
+8\psi(q_{i})[\operatorname{cosh}(1/2+q_{i}-q_{k})-\operatorname{cosh}(3/2-q_{i}+q_{k})]\no\\
&\quad\quad
+8\psi(q_{k})[\operatorname{cosh}(-3/2+2q_{i}-2q_{k})-\operatorname{cosh}(-1/2+2q_{i}-2q_{k})]\no\\
&\quad\quad
+\int_{0}^{q_{k}}\psi [27\operatorname{sinh}(-3/2-3x+2q_{i}+q_{k})+\operatorname{sinh}(-1/2-x+2q_{i}-q_{k})]dx\no\\
&\quad\quad
+\int_{q_{k}}^{q_{i}}\psi[27 \operatorname{sinh}(-1/2-3x+2q_{i}+q_{k})+\operatorname{sinh}(-3/2-x+2q_{i}-q_{k})]dx\no\\
&\quad\quad
+\int_{q_{i}}^{1}\psi[27 \operatorname{sinh}(3/2-3x+2q_{i}+q_{k})+\operatorname{sinh}(1/2-x+2q_{i}-q_{k})]dx
\bigg\}\no\\
&+\frac{3}{2}\sum_{i=1}^{N}\sum_{k<i}p_{i}^{2}p_{k}\bigg\{-2\psi_{x}(q_{k})\operatorname{sinh}(1/2)
+\int_{0}^{q_{k}}\psi \operatorname{sinh}(-1/2-x+q_{k})\mbox{d}x \no\\
&\quad\quad
+\int_{q_{k}}^{q_{i}}\psi \operatorname{sinh}(1/2-x+q_{k})\mbox{d}x
+\int_{q_{i}}^{1}\psi \operatorname{sinh}(1/2-x+q_{k})\mbox{d}x\bigg\}\no\\
&+\frac{3}{2}\!\sum_{i<j<k}\!\!p_{i}p_{j}p_{k}\!\Big\{\!
2\psi_{x}(q_{i})[\operatorname{sinh}(1/2+2q_{i}-q_{j}-q_{k})-\operatorname{sinh}(3/2+2q_{i}-q_{j}-q_{k})\no\\
&\quad\quad
+\operatorname{sinh}(-1/2+q_{j}-q_{k})+\operatorname{sinh}(-1/2-q_{j}+q_{k})]\no\\
&\quad\quad
+8\psi(q_{i})[\operatorname{cosh}(3/2+2q_{i}-q_{j}-q_{k})-\operatorname{cosh}(1/2+2q_{i}-q_{j}-q_{k})]\no\\
&\quad\quad
+2\psi_{x}(q_{j})[\operatorname{sinh}(-1/2+q_{i}-2q_{j}+q_{k})-\operatorname{sinh}(1/2+q_{i}-2q_{j}+q_{k})\no\\
&\quad\quad
+\operatorname{sinh}(1/2+q_{i}-q_{k})-\operatorname{sinh}(3/2+q_{i}-q_{k})]\no\\
&\quad\quad
+8\psi(q_{j})[\operatorname{cosh}(-1/2+q_{i}-2q_{j}+q_{k})-\operatorname{cosh}(1/2+q_{i}-2q_{j}+q_{k})]\no\no\\
&\quad\quad
+2\psi_{x}(q_{k})[\operatorname{sinh}(1/2+q_{i}+q_{j}-2q_{k})-\operatorname{sinh}(3/2+q_{i}+q_{j}-2q_{k}) \no\\
&\quad\quad
+\operatorname{sinh}(-1/2+q_{i}-q_{j})-\operatorname{sinh}(1/2+q_{i}-q_{j})]\no\\
&\quad\quad
+8\psi(q_{k})[\operatorname{cosh}(1/2+q_{i}+q_{j}-2q_{k})-\operatorname{cosh}(3/2+q_{i}+q_{j}-2q_{k})]
\Big\} \no\\
&+\frac{3}{2}\!\!\sum_{i<j<k}p_{i}p_{j}p_{k}\bigg\{\!\!\int_{0}^{q_{i}}\!\!\!
\psi [27 \operatorname{sinh}(-3/2\!-\!3x\!+\!q_{i}\!+\!q_{j}\!+\!q_{k})
\!+\!\operatorname{sinh}(-1/2\!-\!x\!+\!q_{i}\!+\!q_{j}\!-\!q_{k}) \no\\
&\quad\quad\quad\quad
+\operatorname{sinh}(-1/2-x+q_{i}-q_{j}+q_{k})-\operatorname{sinh}(1/2+x+q_{i}-q_{j}-q_{k})]dx\no\\
&\quad\quad
+\int_{q_{i}}^{q_{j}}\psi [27 \operatorname{sinh}(-1/2-3x+q_{i}+q_{j}+q_{k})+
\operatorname{sinh}(1/2-x+q_{i}+q_{j}-q_{k})\no\\
&\quad\quad\quad\quad\quad\quad
+\operatorname{sinh}(1/2-x+q_{i}-q_{j}+q_{k})-\operatorname{sinh}(3/2+x+q_{i}-q_{j}-q_{k})]dx\no\\
&\quad\quad
+\int_{q_{j}}^{q_{k}}\psi [27 \operatorname{sinh}(1/2-3x+q_{i}+q_{j}+q_{k})+
 \operatorname{sinh}(3/2-x+q_{i}+q_{j}-q_{k})\no\\
&\quad\quad\quad\quad\quad\quad
+\operatorname{sinh}(-1/2-x+q_{i}-q_{j}+q_{k})-\operatorname{sinh}(1/2+x+q_{i}-q_{j}-q_{k})]dx\no\\
&\quad\quad
+\int_{q_{k}}^{1}\psi[27 \operatorname{sinh}(3/2-3x+q_{i}+q_{j}+q_{k})+
\operatorname{sinh}(1/2-x+q_{i}+q_{j}-q_{k})\no\\
&\quad\quad\quad\quad\quad\quad
+\operatorname{sinh}(1/2-x+q_{i}-q_{j}+q_{k})-\operatorname{sinh}(-1/2+x+q_{i}-q_{j}-q_{k})]dx
\bigg\},
\end{align}
\begin{align}
\int_{0}^{1}u_{x}^{3}\psi_{xx}\mbox{d}x
=&-\frac{1}{4}\sum_{i=1}^{N}p_{i}^{3}\bigg\{
2\psi_{x}(q_{i})[3\operatorname{sinh}(1/2)-\operatorname{sinh}(3/2)] \no\\
&\quad\quad
+3\int_{0}^{q_{i}}\psi [3\operatorname{sinh}(-3/2-3x+3q_{i})-\operatorname{sinh}(-1/2-x+q_{i})]dx \no\\
&\quad\quad
+3\int_{q_{i}}^{1}\psi [3\operatorname{sinh}(3/2-3x+3q_{i})-\operatorname{sinh}(1/2-x+q_{i})]\mbox{d}x
\bigg\}\no\\
&-\frac{3}{4}\sum_{i=1}^{N}\sum_{k>i}p_{i}^{2}p_{k}\bigg\{
2\psi_{x}(q_{k})[\operatorname{sinh}(1/2+2q_{i}-2q_{k})-
\operatorname{sinh}(3/2+2q_{i}-2q_{k})]\no\\
&\quad\quad +4\psi(q_{i})[\operatorname{cosh}(3/2+q_{i}-q_{k})
-\operatorname{cosh}(1/2-q_{i}+q_{k})]\no\\
&\quad\quad
+4\psi(q_{k})[\operatorname{cosh}(1/2+2q_{i}-2q_{k})
-\operatorname{cosh}(3/2+2q_{i}-2q_{k})] \no\\
&\quad\quad
+\int_{0}^{q_{i}}\psi [9\operatorname{sinh}(-3/2-3x+2q_{i}+q_{k})
-\operatorname{sinh}(-1/2-x+2q_{i}-q_{k})]\mbox{d}x \no\\
&\quad\quad
+\int_{q_{i}}^{q_{k}}\psi[9\operatorname{sinh}(1/2-3x+2q_{i}+q_{k})
- \operatorname{sinh}(3/2-x+2q_{i}-q_{k})]\mbox{d}x\no\\
&\quad\quad
+\int_{q_{k}}^{1}\psi [9\operatorname{sinh}(3/2-3x+2q_{i}+q_{k})
- \operatorname{sinh}(1/2-x+2q_{i}-q_{k})]\mbox{d}x\bigg\}\no\\
&+\frac{3}{2}\sum_{i=1}^{N}\sum_{k>i}p_{i}^{2}p_{k}\bigg\{
-2\psi_{x}(q_{k})\operatorname{sinh}(1/2)
+\int_{0}^{q_{i}}\psi \operatorname{sinh}(-1/2-x+q_{k})\mbox{d}x\no\\
&\quad\quad
+\int_{q_{i}}^{q_{k}}\psi \operatorname{sinh}(-1/2-x+q_{k})\mbox{d}x
+\int_{q_{k}}^{1}\psi \operatorname{sinh}(1/2-x+q_{k})\mbox{d}x
\bigg\} \no\\
&-\frac{3}{4}\sum_{i=1}^{N}\sum_{k<i}p_{i}^{2}p_{k}\bigg\{2\psi_{x}(q_{k})[\operatorname{sinh}(-3/2+2q_{i}-2q_{k})
-\operatorname{sinh}(-1/2+2q_{i}-2q_{k})]\no\\
&\quad\quad
+4\psi(q_{i})[\operatorname{cosh}(1/2+q_{i}-q_{k})-\operatorname{cosh}(3/2-q_{i}+q_{k})]\no\\
&\quad\quad
+4\psi(q_{k})[\operatorname{cosh}(-3/2+2q_{i}-2q_{k})-\operatorname{cosh}(-1/2+2q_{i}-2q_{k})]\no\\
&\quad\quad
+\int_{0}^{q_{k}}\psi [9 \operatorname{sinh}(-3/2-3x+2q_{i}+q_{k})-\operatorname{sinh}(-1/2-x+2q_{i}-q_{k})]\mbox{d}x\no\\
&\quad\quad
+\int_{q_{k}}^{q_{i}}\psi [9\operatorname{sinh}(-1/2-3x+2q_{i}+q_{k})-\operatorname{sinh}(-3/2-x+2q_{i}-q_{k})]\mbox{d}x\no\\
&\quad\quad
+\int_{q_{i}}^{1}\psi [9 \operatorname{sinh}(3/2-3x+2q_{i}+q_{k})-\operatorname{sinh}(1/2-x+2q_{i}-q_{k})]\mbox{d}x
\bigg\}\no\\
&+\frac{3}{2}\sum_{i=1}^{N}\sum_{k<i}p_{i}^{2}p_{k}\bigg\{
-2\psi_{x}(q_{k})\operatorname{sinh}(1/2)
+\int_{0}^{q_{k}}\psi \operatorname{sinh}(-1/2-x+q_{k})\mbox{d}x\no\\
&\quad\quad
+\int_{q_{k}}^{q_{i}}\psi \operatorname{sinh}(1/2-x+q_{k})\mbox{d}x
+\int_{q_{i}}^{1}\psi \operatorname{sinh}(1/2-x+q_{k})\mbox{d}x
\bigg\}\no\\
&-\frac{3}{2}\sum_{i<j<k}p_{i}p_{j}p_{k}\bigg\{
2\psi_{x}(q_{i})[\operatorname{sinh}(1/2+2q_{i}-q_{j}-q_{k})
-\operatorname{sinh}(3/2+2q_{i}-q_{j}-q_{k})\no\\
&\quad\quad
+\operatorname{sinh}(1/2-q_{j}+q_{k})
+\operatorname{sinh}(1/2+q_{j}-q_{k})]\no\\
&\quad\quad
+4\psi(q_{i})[\operatorname{cosh}(3/2+2q_{i}-q_{j}-q_{k})
-\operatorname{cosh}(1/2+2q_{i}-q_{j}-q_{k})]\no\\
&\quad\quad
+2\psi_{x}(q_{j})[\operatorname{sinh}(-1/2+q_{i}-2q_{j}+q_{k})
-\operatorname{sinh}(1/2+q_{i}-2q_{j}+q_{k})\no\\
&\quad\quad
+\operatorname{sinh}(3/2+q_{i}-q_{k})
-\operatorname{sinh}(1/2+q_{i}-q_{k})]\no\\
&\quad\quad
+4\psi(q_{j})[\operatorname{cosh}(-1/2+q_{i}-2q_{j}+q_{k})
-\operatorname{cosh}(1/2+q_{i}-2q_{j}+q_{k})]\no\\
&\quad\quad
+2\psi_{x}(q_{k})[\operatorname{sinh}(1/2+q_{i}+q_{j}-2q_{k})
-\operatorname{sinh}(3/2+q_{i}+q_{j}-2q_{k})\no\\
&\quad\quad
+\operatorname{sinh}(1/2+q_{i}-q_{j})
-\operatorname{sinh}(-1/2+q_{i}-q_{j})]\no\\
&\quad\quad
+4\psi(q_{k})[\operatorname{cosh}(1/2+q_{i}+q_{j}-2q_{k})
-\operatorname{cosh}(3/2+q_{i}+q_{j}-2q_{k})]
\bigg\}\no\\
&-\frac{3}{2}\!\sum_{i<j<k}p_{i}p_{j}p_{k}\!\bigg\{\!\!\int_{0}^{q_{i}}\!\!\psi [9\operatorname{sinh}(-3/2\!-\!3x\!+\!q_{i}\!+\!q_{j}\!+\!q_{k})\!-\!\operatorname{sinh}(-1/2\!-\!x\!+\!q_{i}\!+\!q_{j}\!-\!q_{k})
\no\\
&\quad\quad\quad\quad
-\operatorname{sinh}(-1/2-x+q_{i}-q_{j}+q_{k})+\operatorname{sinh}(1/2+x+q_{i}-q_{j}-q_{k})]\mbox{d}x\no\\
&\quad\quad +\int_{q_{i}}^{q_{j}}\psi [9\operatorname{sinh}(-1/2-3x+q_{i}+q_{j}+q_{k})
-\operatorname{sinh}(1/2-x+q_{i}+q_{j}-q_{k})\no\\
&\quad\quad\quad\quad
 -\operatorname{sinh}(1/2-x+q_{i}-q_{j}+q_{k})+\operatorname{sinh}(3/2+x+q_{i}-q_{j}-q_{k})]\mbox{d}x\no\\
&\quad\quad
+\int_{q_{j}}^{q_{k}}\psi [9 \operatorname{sinh}(1/2-3x+q_{i}+q_{j}+q_{k})
-\operatorname{sinh}(3/2-x+q_{i}+q_{j}-q_{k})\no\\
&\quad\quad\quad\quad
-\operatorname{sinh}(-1/2-x+q_{i}-q_{j}+q_{k})+\operatorname{sinh}(1/2+x+q_{i}-q_{j}-q_{k})]\mbox{d}x\no\\
&\quad\quad
+\int_{q_{k}}^{1}\psi [9 \operatorname{sinh}(3/2-3x+q_{i}+q_{j}+q_{k})
-\operatorname{sinh}(1/2-x+q_{i}+q_{j}-q_{k})\no\\
&\quad\quad\quad\quad
-\operatorname{sinh}(1/2-x+q_{i}-q_{j}+q_{k})+\operatorname{sinh}(-1/2+x+q_{i}-q_{j}-q_{k})]\mbox{d}x\Big\},
\end{align}
\begin{align}
\int_{0}^{1}u^{3}\psi_{x}\mbox{d}x
=&\frac{1}{4}\sum_{i=1}^{N}p_{i}^{3}\bigg\{3\int_{0}^{q_{i}}\psi [\operatorname{sinh}(-3/2-3x+3q_{i})
+\operatorname{sinh}(-1/2-x+q_{i})]\mbox{d}x \no\\
&\quad\quad
+3\int_{q_{i}}^{1}\psi [\operatorname{sinh}(3/2-3x+3q_{i})+ \operatorname{sinh}(1/2-x+q_{i})]\mbox{d}x
\bigg\}\no\\
& +\frac{3}{4}\sum_{i=1}^{N}\sum_{k>i}p_{i}^{2}p_{k}\bigg\{
\int_{0}^{q_{i}}\psi [3\operatorname{sinh}(-3/2-3x+2q_{i}+q_{k})+\operatorname{sinh}(-1/2-x+2q_{i}-q_{k})]\mbox{d}x \no\\
&\quad\quad
+\int_{q_{i}}^{q_{k}}\psi[3 \operatorname{sinh}(1/2-3x+2q_{i}+q_{k})
+ \operatorname{sinh}(3/2-x+2q_{i}-q_{k})]\mbox{d}x\no\\
&\quad\quad
+\int_{q_{k}}^{1}\psi [3\operatorname{sinh}(3/2-3x+2q_{i}+q_{k})\mbox{d}x
+\operatorname{sinh}(1/2-x+2q_{i}-q_{k})]\mbox{d}x\bigg\}\no\\
&+\frac{3}{2}\sum_{i=1}^{N}\sum_{k>i}p_{i}^{2}p_{k}\bigg\{
\left(\int_{0}^{q_{i}}+\int_{q_{i}}^{q_{k}}\right)\psi \operatorname{sinh}(-1/2-x+q_{k})\mbox{d}x
+\int_{q_{k}}^{1}\psi \operatorname{sinh}(1/2-x+q_{k})\mbox{d}x\bigg\}\no\\
&+\frac{3}{4}\sum_{i=1}^{N}\sum_{k<i}p_{i}^{2}p_{k}\bigg\{
\int_{0}^{q_{k}}\psi [3\operatorname{sinh}(-3/2-3x+2q_{i}+q_{k})+\operatorname{sinh}(-1/2-x+2q_{i}-q_{k})]\mbox{d}x\no\\
&\quad\quad
+\int_{q_{k}}^{q_{i}}\psi[3 \operatorname{sinh}(-1/2-3x+2q_{i}+q_{k})+\operatorname{sinh}(-3/2-x+2q_{i}-q_{k})]\mbox{d}x\no\\
&\quad\quad
+\int_{q_{i}}^{1}\psi [3\operatorname{sinh}(3/2-3x+2q_{i}+q_{k})+\operatorname{sinh}(1/2-x+2q_{i}-q_{k})]\mbox{d}x\bigg\}\no\\
&+\frac{3}{2}\sum_{i=1}^{N}\sum_{k<i}p_{i}^{2}p_{k}\bigg\{
\int_{0}^{q_{k}}\psi \operatorname{sinh}(-1/2-x+q_{k})\mbox{d}x
+\left(\int_{q_{k}}^{q_{i}}+\int_{q_{i}}^{1}\right)\psi \operatorname{sinh}(1/2-x+q_{k})\mbox{d}x
\bigg\}\no\\
&+\frac{3}{2}\sum_{i<j<k}p_{i}p_{j}p_{k}\bigg\{
\int_{0}^{q_{i}}\psi [3\operatorname{sinh}(-3/2-3x+q_{i}+q_{j}+q_{k})+\operatorname{sinh}(-1/2-x+q_{i}+q_{j}-q_{k})\no\\
&\quad\quad\quad\quad
+\operatorname{sinh}(-1/2-x+q_{i}-q_{j}+q_{k})-\operatorname{sinh}(1/2+x+q_{i}-q_{j}-q_{k})]\mbox{d}x\no\\
&\quad\quad
+\int_{q_{i}}^{q_{j}}\psi [3\operatorname{sinh}(-1/2-3x+q_{i}+q_{j}+q_{k})+\operatorname{sinh}(1/2-x+q_{i}+q_{j}-q_{k})\no\\
&\quad\quad\quad\quad
+\operatorname{sinh}(1/2-x+q_{i}-q_{j}+q_{k})-\operatorname{sinh}(3/2+x+q_{i}-q_{j}-q_{k})]\mbox{d}x\no\\
&\quad\quad
+\int_{q_{j}}^{q_{k}}\psi [3\operatorname{sinh}(1/2-3x+q_{i}+q_{j}+q_{k})
+ \operatorname{sinh}(3/2-x+q_{i}+q_{j}-q_{k})\no\\
&\quad\quad\quad\quad
+ \operatorname{sinh}(-1/2-x+q_{i}-q_{j}+q_{k})-\operatorname{sinh}(1/2+x+q_{i}-q_{j}-q_{k})]\mbox{d}x\no\\
&\quad\quad
+\int_{q_{k}}^{1}\psi[3 \operatorname{sinh}(3/2-3x+q_{i}+q_{j}+q_{k})+\operatorname{sinh}(1/2-x+q_{i}+q_{j}-q_{k})\no\\
&\quad\quad\quad\quad
+ \operatorname{sinh}(1/2-x+q_{i}-q_{j}+q_{k})-\operatorname{sinh}(-1/2+x+q_{i}-q_{j}-q_{k})]\mbox{d}x\bigg\},
\end{align}
\begin{align}
\int_{0}^{1}u u_{x}^{2}\psi_{x}\mbox{d}x
=&\frac{1}{4}\sum_{i=1}^{N}p_{i}^{3}\bigg\{
\int_{0}^{q_{i}}\psi[3 \operatorname{sinh}(-3/2-3x+3q_{i})-\operatorname{sinh}(-1/2-x+q_{i})]\mbox{d}x \no\\
&\quad\quad
+\int_{q_{i}}^{1}\psi[3 \operatorname{sinh}(3/2-3x+3q_{i})-\operatorname{sinh}(1/2-x+q_{i})]\mbox{d}x\bigg\}\no\\
&+\frac{1}{2}\sum_{i=1}^{N}\sum_{k>i}p_{i}^{2}p_{k}\bigg\{
2\psi(q_{i})[\operatorname{cosh}(3/2+q_{i}-q_{k})-\operatorname{cosh}(1/2-q_{i}+q_{k})]\no\\
&\quad\quad
+2\psi(q_{k})[\operatorname{cosh}(1/2+2q_{i}-2q_{k})-\operatorname{cosh}(3/2+2q_{i}-2q_{k})]\no\\
&\quad\quad
+\int_{0}^{q_{i}}\psi [3\operatorname{sinh}(-3/2-3x+2q_{i}+q_{k})-\operatorname{sinh}(-1/2-x+2q_{i}-q_{k})]\mbox{d}x\no\\
&\quad\quad
+\int_{q_{i}}^{q_{k}}\psi [3\operatorname{sinh}(1/2-3x+2q_{i}+q_{k})- \operatorname{sinh}(3/2-x+2q_{i}-q_{k})]\mbox{d}x\no\\
&\quad\quad
+\int_{q_{k}}^{1}\psi [3\operatorname{sinh}(3/2-3x+2q_{i}+q_{k})-\operatorname{sinh}(1/2-x+2q_{i}-q_{k})]\mbox{d}x\bigg\}\no\\
&+\frac{1}{2}\sum_{i=1}^{N}\sum_{k<i}p_{i}^{2}p_{k}\bigg\{
2\psi(q_{i})[\operatorname{cosh}(1/2+q_{i}-q_{k})-\operatorname{cosh}(3/2-q_{i}+q_{k})]\no\\
&\quad\quad
+2\psi(q_{k})[\operatorname{cosh}(-3/2+2q_{i}-2q_{k})-\operatorname{cosh}(-1/2+2q_{i}-2q_{k})]\no\\
&\quad\quad
+\int_{0}^{q_{k}}\psi[3 \operatorname{sinh}(-3/2-3x+2q_{i}+q_{k})-\operatorname{sinh}(-1/2-x+2q_{i}-q_{k})]\mbox{d}x\no\\
&\quad\quad
+\int_{q_{k}}^{q_{i}}\psi [3\operatorname{sinh}(-1/2-3x+2q_{i}+q_{k})-\operatorname{sinh}(-3/2-x+2q_{i}-q_{k})]\mbox{d}x\no\\
&\quad\quad
+\int_{q_{i}}^{1}\psi[3 \operatorname{sinh}(3/2-3x+2q_{i}+q_{k})-\operatorname{sinh}(1/2-x+2q_{i}-q_{k})]\mbox{d}x\bigg\}\no\\
&+\frac{1}{4}\sum_{j=1}^{N}\sum_{i>j}p_{j}^{2}p_{i}\bigg\{
\int_{0}^{q_{j}}\psi [3\operatorname{sinh}(-3/2-3x+2q_{j}+q_{i})+\operatorname{sinh}(-1/2-x+2q_{j}-q_{i})]\mbox{d}x\no\\
&\quad\quad
+\int_{q_{j}}^{q_{i}}\psi [3\operatorname{sinh}(1/2-3x+2q_{j}+q_{i})+\operatorname{sinh}(3/2-x+2q_{j}-q_{i})]\mbox{d}x\no\\
&\quad\quad
+\int_{q_{i}}^{1}\psi [3\operatorname{sinh}(3/2-3x+2q_{j}+q_{i})+\operatorname{sinh}(1/2-x+2q_{j}-q_{i})]\mbox{d}x\bigg\}\no\\
&-\frac{1}{2}\sum_{j=1}^{N}\sum_{i>j}p_{j}^{2}p_{i}\bigg\{
\left(\int_{0}^{q_{j}}\!+\!\int_{q_{j}}^{q_{i}}\right)\psi \operatorname{sinh}(-1/2-x+q_{i})\mbox{d}x
+\int_{q_{i}}^{1}\psi \operatorname{sinh}(1/2-x+q_{i})\mbox{d}x\bigg\}\no\\
&+\frac{1}{4}\sum_{j=1}^{N}\sum_{i<j}p_{j}^{2}p_{i}\bigg\{
\int_{0}^{q_{i}}\psi[3 \operatorname{sinh}(-3/2-3x+2q_{j}+q_{i})+\operatorname{sinh}(-1/2-x+2q_{j}-q_{i})]\mbox{d}x\no\\
&\quad\quad
+\int_{q_{i}}^{q_{j}}\psi[3 \operatorname{sinh}(-1/2-3x+2q_{j}+q_{i})
+ \operatorname{sinh}(-3/2-x+2q_{j}-q_{i})]\mbox{d}x\no\\
&\quad\quad
+\int_{q_{j}}^{1}\psi [3\operatorname{sinh}(3/2-3x+2q_{j}+q_{i})+\operatorname{sinh}(1/2-x+2q_{j}-q_{i})]\mbox{d}x\bigg\}\no\\
&-\frac{1}{2}\sum_{j=1}^{N}\sum_{i<j}p_{j}^{2}p_{i}\bigg\{
\int_{0}^{q_{i}}\psi \operatorname{sinh}(-1/2-x+q_{i})\mbox{d}x
+\left(\int_{q_{i}}^{q_{j}}+\int_{q_{j}}^{1}\right)\psi \operatorname{sinh}(1/2-x+q_{i})\mbox{d}x\bigg\}\no\\
&+\frac{1}{2}\sum_{i<j<k}p_{i}p_{j}p_{k}\bigg\{
4\psi(q_{i})[\operatorname{cosh}(-3/2+q_{j}+q_{k}-2q_{i})-\operatorname{cosh}(-1/2+q_{j}+q_{k}-2q_{i})]\no\\
&\quad\quad
+4\psi(q_{j})[\operatorname{cosh}(-1/2+q_{i}+q_{k}-2q_{j})-\operatorname{cosh}(1/2+q_{i}+q_{k}-2q_{j})]\no\\
&\quad\quad
+4\psi(q_{k})[\operatorname{cosh}(1/2+q_{j}+q_{i}-2q_{k})-\operatorname{cosh}(3/2+q_{j}+q_{i}-2q_{k})]\no\\
&\quad\quad
+\int_{0}^{q_{i}}\psi[9 \operatorname{sinh}(-3/2-3x+q_{i}+q_{j}+q_{k})
-\operatorname{sinh}(-1/2-x+q_{i}+q_{j}-q_{k})\no\\
&\quad\quad
-\operatorname{sinh}(-1/2-x+q_{i}-q_{j}+q_{k})+\operatorname{sinh}(1/2+x+q_{i}-q_{j}-q_{k})]\mbox{d}x\no\\
&\quad\quad
+\int_{q_{i}}^{q_{j}}\psi [9 \operatorname{sinh}(-1/2-3x+q_{i}+q_{j}+q_{k})
-\operatorname{sinh}(1/2-x+q_{i}+q_{j}-q_{k})\no\\
&\quad\quad
-\operatorname{sinh}(1/2-x+q_{i}-q_{j}+q_{k})+\operatorname{sinh}(3/2+x+q_{i}-q_{j}-q_{k})]\mbox{d}x\no\\
&\quad\quad
+\int_{q_{j}}^{q_{k}}\psi [9\operatorname{sinh}(1/2-3x+q_{i}+q_{j}+q_{k})
-\operatorname{sinh}(3/2-x+q_{i}+q_{j}-q_{k})\no\\
&\quad\quad
-\operatorname{sinh}(-1/2-x+q_{i}-q_{j}+q_{k})+\operatorname{sinh}(1/2+x+q_{i}-q_{j}-q_{k})]\mbox{d}x\no\\
&\quad\quad
+\int_{q_{k}}^{1}\psi[9 \operatorname{sinh}(3/2-3x+q_{i}+q_{j}+q_{k})
- \operatorname{sinh}(1/2-x+q_{i}+q_{j}-q_{k}) \no\\
&\quad\quad
-\operatorname{sinh}(1/2-x+q_{i}-q_{j}+q_{k})+\operatorname{sinh}(-1/2+x+q_{i}-q_{j}-q_{k})]\mbox{d}x
\bigg\},
\end{align}
\begin{align}
\int_{0}^{1}u^{2}\psi_{xxx}\mbox{d}x
=&2\sum_{i=1}^{N}p_{i}^{2}\bigg\{-\psi_{x}(q_{i})\operatorname{sinh}(1)
\!+\!2\!\int_{0}^{q_{i}}\!\!\psi\operatorname{sinh}(-1-2x+2q_{i})\mbox{d}x
\!+\!2\!\int_{q_{i}}^{1}\!\!\psi \operatorname{sinh}(1-2x+2q_{i})\mbox{d}x\bigg\}\no\\
&+\sum_{i=1}^{N}\sum_{j>i}p_{i}p_{j}\bigg\{
[\psi_{x}(q_{j})+\psi_{x}(q_{i})][\operatorname{sinh}(q_{i}-q_{j})-\operatorname{sinh}(1+q_{i}-q_{j})]\no\\
&\quad\quad
+2[\psi(q_{j})-\psi(q_{i})][\operatorname{cosh}(q_{i}-q_{j})-\operatorname{cosh}(1+q_{i}-q_{j})]\no\\
&\quad\quad
+4\int_{0}^{q_{i}}\psi \operatorname{sinh}(-1-2x+q_{i}+q_{j})\mbox{d}x
+4\int_{q_{i}}^{q_{j}}\psi \operatorname{sinh}(-2x+q_{i}+q_{j})\mbox{d}x\no\\
&\quad\quad
+4\int_{q_{j}}^{1}\psi \operatorname{sinh}(1-2x+q_{i}+q_{j})\mbox{d}x\bigg\} \no\\
&+\sum_{i=1}^{N}\sum_{j<i}p_{i}p_{j}\bigg\{
[\psi_{x}(q_{i})+\psi_{x}(q_{j})][\operatorname{sinh}(-q_{i}+q_{j})-\operatorname{sinh}(1-q_{i}+q_{j})]\no\\
&\quad\quad
+2[\psi(q_{i})-\psi(q_{j})][\operatorname{cosh}(-q_{i}+q_{j})-\operatorname{cosh}(1-q_{i}+q_{j})]\no\\
&\quad\quad
+4\int_{0}^{q_{j}}\psi \operatorname{sinh}(-1-2x+q_{i}+q_{j})\mbox{d}x
+4\int_{q_{j}}^{q_{i}}\psi \operatorname{sinh}(-2x+q_{i}+q_{j})\mbox{d}x\no\\
&\quad\quad
+4\int_{q_{i}}^{1}\psi \operatorname{sinh}(1-2x+q_{i}+q_{j})\mbox{d}x
\bigg\},
\end{align}
\begin{align}
\int_{0}^{1}u^{2}\psi_{x}\mbox{d}x
=&\sum_{i=1}^{N}p_{i}^{2}\bigg\{\int_{0}^{q_{i}}\psi \operatorname{sinh}(-1-2x+2q_{i})\mbox{d}x
+\int_{q_{i}}^{1}\psi \operatorname{sinh}(1-2x+2q_{i})\mbox{d}x\bigg\}\no\\
&
+\sum_{i=1}^{N}\sum_{j>i}p_{i}p_{j}\bigg\{
\int_{0}^{q_{i}}\psi \operatorname{sinh}(-1-2x+q_{i}+q_{j})\mbox{d}x
+\int_{q_{i}}^{q_{j}}\psi \operatorname{sinh}(-2x+q_{i}+q_{j})\mbox{d}x\no\\
&\quad\quad
+\int_{q_{j}}^{1}\psi \operatorname{sinh}(1-2x+q_{i}+q_{j})\mbox{d}x
\bigg\}\no\\
&+\sum_{i=1}^{N}\sum_{j<i}p_{i}p_{j}\bigg\{
\int_{0}^{q_{j}}\psi \operatorname{sinh}(-1-2x+q_{i}+q_{j})\mbox{d}x
+\int_{q_{j}}^{q_{i}}\psi \operatorname{sinh}(-2x+q_{i}+q_{j})\mbox{d}x\no\\
&\quad\quad
+\int_{q_{i}}^{1}\psi \operatorname{sinh}(1-2x+q_{i}+q_{j})\mbox{d}x
\bigg\},
\end{align}
\begin{align}
\int_{0}^{1}u_{x}^{2}\psi_{x}\mbox{d}x
=&\sum_{i=1}^{N}p_{i}^{2}\bigg\{\int_{0}^{q_{i}}\psi \operatorname{sinh}(-1-2x+2q_{i})\mbox{d}x
+\int_{q_{i}}^{1}\psi \operatorname{sinh}(1-2x+2q_{i})\mbox{d}x
\bigg\}\no\\
&+\sum_{i=1}^{N}\sum_{j>i}p_{i}p_{j}\bigg\{
[\psi(q_{i})-\psi(q_{j})][\operatorname{cosh}(1+q_{i}-q_{j})-\operatorname{cosh}(q_{i}-q_{j})]\no\\
&\quad\quad
+\int_{0}^{q_{i}}\psi \operatorname{sinh}(-1-2x+q_{i}+q_{j})\mbox{d}x
+\int_{q_{i}}^{q_{j}}\psi \operatorname{sinh}(-2x+q_{i}+q_{j})\mbox{d}x\no\\
&\quad\quad
+\int_{q_{j}}^{1}\psi \operatorname{sinh}(1-2x+q_{i}+q_{j})\mbox{d}x
\bigg\}\no\\
&+\sum_{i=1}^{N}\sum_{j<i}p_{i}p_{j}\bigg\{
[\psi(q_{i})-\psi(q_{j})][\operatorname{cosh}(q_{i}-q_{j})-\operatorname{cosh}(1-q_{i}+q_{j})]\no\\
&\quad\quad
+\int_{0}^{q_{j}}\psi \operatorname{sinh}(-1-2x+q_{i}+q_{j})\mbox{d}x
+\int_{q_{j}}^{q_{i}}\psi \operatorname{sinh}(-2x+q_{i}+q_{j})\mbox{d}x \no\\
&\quad\quad
+\int_{q_{i}}^{1}\psi \operatorname{sinh}(1-2x+q_{i}+q_{j})\mbox{d}x
\bigg\},
\end{align}
\begin{align}\label{px}
\int_{0}^{1}u_{x}^{3}\psi\mbox{d}x
=&-\frac{1}{4}\sum_{i=1}^{N}p_{i}^{3}\bigg\{
\int_{0}^{q_{i}}\psi [\operatorname{sinh}(-3/2-3x+3q_{i})-3\operatorname{sinh}(-1/2-x+q_{i})]\mbox{d}x\no\\
&\quad\quad
+\int_{q_{i}}^{1}\psi [\operatorname{sinh}(3/2-3x+3q_{i})-3\operatorname{sinh}(1/2-x+q_{i})]\mbox{d}x\bigg\}\no\\
&-\frac{3}{4}\sum_{i=1}^{N}\sum_{k>i}p_{i}^{2}p_{k}\bigg\{
\int_{0}^{q_{i}}\psi[ \operatorname{sinh}(-3/2-3x+2q_{i}+q_{k})-\operatorname{sinh}(-1/2-x+2q_{i}-q_{k})]\mbox{d}x\no\\
&\quad\quad
+\int_{q_{i}}^{q_{k}}\psi [\operatorname{sinh}(-1/2-3x+2q_{i}+q_{k})
-\operatorname{sinh}(3/2-x+2q_{i}-q_{k})]\mbox{d}x\no\\
&\quad\quad
+\int_{q_{k}}^{1}\psi [\operatorname{sinh}(3/2-3x+2q_{i}+q_{k})
- \operatorname{sinh}(1/2-x+2q_{i}-q_{k})]\mbox{d}x\bigg\}\no\\
&+\frac{3}{2}\sum_{i=1}^{N}\sum_{k>i}p_{i}^{2}p_{k}\bigg\{
\left(\int_{0}^{q_{i}}+\int_{q_{i}}^{q_{k}}\right)\psi \operatorname{sinh}(-1/2-x+q_{k})\mbox{d}x
+\int_{q_{k}}^{1}\psi \operatorname{sinh}(1/2-x+q_{k})\mbox{d}x
\bigg\}\no\\
&-\frac{3}{4}\sum_{i=1}^{N}\sum_{k<i}p_{i}^{2}p_{k}\bigg\{
\int_{0}^{q_{k}}\psi [\operatorname{sinh}(-3/2-3x+2q_{i}+q_{k})-\operatorname{sinh}(-1/2-x+2q_{i}-q_{k})]\mbox{d}x\no\\
&\quad\quad
+\int_{q_{k}}^{q_{i}}\psi [\operatorname{sinh}(-1/2-3x+2q_{i}+q_{k})
-\operatorname{sinh}(-3/2-x+2q_{i}-q_{k})]\mbox{d}x\no\\
&\quad\quad
+\int_{q_{i}}^{1}\psi [\operatorname{sinh}(3/2-3x+2q_{i}+q_{k})-\operatorname{sinh}(1/2-x+2q_{i}-q_{k})]\mbox{d}x
\bigg\}\no\\
&
+\frac{3}{2}\sum_{i=1}^{N}\sum_{k<i}p_{i}^{2}p_{k}\bigg\{
\int_{0}^{q_{k}}\psi \operatorname{sinh}(-1/2-x+q_{k})\mbox{d}x
+\left(\int_{q_{k}}^{q_{i}}+\int_{q_{i}}^{1}\right)\psi \operatorname{sinh}(1/2-x+q_{k})\mbox{d}x
\bigg\}\no\\
&
-\frac{3}{2}\sum_{i<j<k}p_{i}p_{j}p_{k}\bigg\{
\int_{0}^{q_{i}}\psi [\operatorname{sinh}(-3/2-3x+q_{i}+q_{j}+q_{k})
-\operatorname{sinh}(-1/2-x+q_{i}+q_{j}-q_{k})\no\\
&\quad\quad\quad\quad
- \operatorname{sinh}(-1/2-x+q_{i}-q_{j}+q_{k})
+\operatorname{sinh}(1/2+x+q_{i}-q_{j}-q_{k})]\mbox{d}x\no\\
&\quad\quad
+\int_{q_{i}}^{q_{j}}\psi [\operatorname{sinh}(-1/2-3x+q_{i}+q_{j}+q_{k})
- \operatorname{sinh}(1/2-x+q_{i}+q_{j}-q_{k})\no\\
&\quad\quad\quad\quad
-\operatorname{sinh}(1/2-x+q_{i}-q_{j}+q_{k})
+\operatorname{sinh}(3/2+x+q_{i}-q_{j}-q_{k})]\mbox{d}x\no\\
&\quad\quad
+\int_{q_{j}}^{q_{k}}\psi [\operatorname{sinh}(1/2-3x+q_{i}+q_{j}+q_{k})
-\operatorname{sinh}(3/2-x+q_{i}+q_{j}-q_{k})\no\\
&\quad\quad\quad\quad
-\operatorname{sinh}(-1/2-x+q_{i}-q_{j}+q_{k})
+\operatorname{sinh}(1/2+x+q_{i}-q_{j}-q_{k})]\mbox{d}x\no\\
&\quad\quad
+\int_{q_{k}}^{1}\psi [\operatorname{sinh}(3/2-3x+q_{i}+q_{j}+q_{k})
-\operatorname{sinh}(1/2-x+q_{i}+q_{j}-q_{k})\no\\
&\quad\quad\quad\quad
-\operatorname{sinh}(1/2-x+q_{i}-q_{j}+q_{k})
+\operatorname{sinh}(-1/2+x+q_{i}-q_{j}-q_{k})]\mbox{d}x
\bigg\}.
\end{align}

Substituting Eqs.~(\ref{ut-p})-(\ref{px}) into Eq.~(\ref{multi-p}), one finds the coefficient of the term $\psi(q_{m})$ are
\begin{align}\label{coe1}
&2\operatorname{sinh}(1/2)\dot{p}_{m}+\frac{k_2}{2}\bigg\{\sum_{k>m}p_{m}^{2}p_{k}[
\operatorname{cosh}(3/2+q_{m}-q_{k})-\operatorname{cosh}(1/2-q_{m}+q_{k})] \no\\
&+\sum_{i<m}p_{i}^{2}p_{m}[
\operatorname{cosh}(1/2+2q_{i}-2q_{m})
-\operatorname{cosh}(3/2+2q_{i}-2q_{m})]\no\\
&+\sum_{k<m}p_{m}^{2}p_{k}[
\operatorname{cosh}(1/2+q_{m}-q_{k})
-\operatorname{cosh}(3/2-q_{m}+q_{k})]\no\\
&+\sum_{i>m}p_{i}^{2}p_{m}[
\operatorname{cosh}(-3/2+2q_{i}-2q_{m})
-\operatorname{cosh}(-1/2+2q_{i}-2q_{m})]\no\\
&+2\sum_{m<j<k}p_{m}p_{j}p_{k}[
\operatorname{cosh}(3/2+2q_{m}-q_{j}-q_{k})
-\operatorname{cosh}(1/2+2q_{m}-q_{j}-q_{k})]\no\\
&
+2\sum_{i<m<k}p_{i}p_{m}p_{k}[
\operatorname{cosh}(-1/2+q_{i}-2q_{m}+q_{k})
-\operatorname{cosh}(1/2+q_{i}-2q_{m}+q_{k})
]\no\\
&
+2\sum_{i<j<m}p_{i}p_{j}p_{m}[
\operatorname{cosh}(1/2+q_{i}+q_{j}-2q_{m})
-\operatorname{cosh}(3/2+q_{i}+q_{j}-2q_{m})
]\bigg\}\no\\
&
+k_3\bigg\{
\sum_{j>m}p_{m}p_{j}[
\operatorname{cosh}(1+q_{m}-q_{j})
-\operatorname{cosh}(q_{m}-q_{j})
]
+\sum_{j<m}p_{m}p_{j}[
\operatorname{cosh}(q_{j}-q_{m})
-\operatorname{cosh}(1+q_{j}-q_{m})
]
\bigg\}
\end{align}
and the coefficient of $\psi_{x}(q_{m})$ are
\begin{align}\label{coe2}
&2\operatorname{sinh}(1/2)p_{m}\dot{q}_{m}
-\left(\frac{k_1}{3}+\frac{k_2}{2}\right)p_{m}^{3}
[\operatorname{sinh}(3/2)+\operatorname{sinh}(1/2)]
-\frac{2k_1}{3}p_{m}^{3}\operatorname{sinh}(1/2) \no\\
&
-k_{3}p_{m}^{2}\operatorname{sinh}(1)
+k_{1}\bigg\{
\sum_{k>m}p_{m}^{2}p_{k}[
-\operatorname{sinh}(3/2+q_{m}-q_{k})
+\operatorname{sinh}(-1/2+q_{m}-q_{k})
]\no\\
&
-2\operatorname{sinh}(1/2)\sum_{i<m}p_{i}^{2}p_{m}
+\sum_{k<m}p_{m}^{2}p_{k}[
\operatorname{sinh}(-1/2-q_{m}+q_{k})
+\operatorname{sinh}(-3/2+q_{m}-q_{k})
]\no\\
&
-2\operatorname{sinh}(1/2)\sum_{i>m}p_{i}^{2}p_{m}
+2\sum_{m<j<k}p_{m}p_{j}p_{k}[
\operatorname{sinh}(-1/2+q_{j}-q_{k})
-\operatorname{sinh}(1/2+q_{j}-q_{k})
]\no\\
&
+2\sum_{i<m<k}p_{i}p_{m}p_{k}[
\operatorname{sinh}(1/2+q_{i}-q_{k})
-\operatorname{sinh}(3/2+q_{i}-q_{k})
]\no\\
&
+2\sum_{i<j<m}p_{i}p_{j}p_{m}[
\operatorname{sinh}(-1/2+q_{i}-q_{j})
-\operatorname{sinh}(1/2+q_{i}-q_{j})
]
\bigg\}\no\\
&
+k_{2}\bigg\{\sum_{k>m}p_{m}^{2}p_{k}[-\operatorname{sinh}(3/2+q_{m}-q_{k})+\operatorname{sinh}(-1/2+q_{m}-q_{k})
]\no\\
&+\frac12\sum_{i<m}p_{i}^{2}p_{m}[\operatorname{sinh}(1/2+2q_{i}-2q_{m})-\operatorname{sinh}(3/2+2q_{i}-2q_{m})
-2\operatorname{sinh}(1/2)]\no \\
&+\sum_{k<m}p_{m}^{2}p_{k}[
\operatorname{sinh}(-1/2-q_{m}+q_{k})
+\operatorname{sinh}(-3/2+q_{m}-q_{k})
]\no\\
&
+\frac12\sum_{i>m}p_{i}^{2}p_{m}[\operatorname{sinh}(-3/2+2q_{i}-2q_{m})-\operatorname{sinh}(-1/2+2q_{i}-2q_{m})
-2\operatorname{sinh}(1/2)]\no\\
&
+\sum_{m<j<k}p_{m}p_{j}p_{k}[
-\operatorname{sinh}(3/2+2q_{m}-q_{j}-q_{k})
+\operatorname{sinh}(1/2+2q_{m}-q_{j}-q_{k})\no\\
&\quad\quad
+\operatorname{sinh}(-1/2+q_{j}-q_{k})
+\operatorname{sinh}(-1/2-q_{j}+q_{k})
]\no\\
&
+\sum_{i<m<k}p_{i}p_{m}p_{k}[
\operatorname{sinh}(-1/2+q_{i}-2q_{m}+q_{k})
-\operatorname{sinh}(1/2+q_{i}-2q_{m}+q_{k})\no\\
&\quad\quad
+\operatorname{sinh}(1/2+q_{i}-q_{k})
-\operatorname{sinh}(3/2+q_{i}-q_{k})
]\no\\
&
+\sum_{i<j<m}p_{i}p_{j}p_{m}[
\operatorname{sinh}(1/2+q_{i}+q_{j}-2q_{m})
-\operatorname{sinh}(3/2+q_{i}+q_{j}-2q_{m})\no\\
&\quad\quad
+\operatorname{sinh}(-1/2+q_{i}-q_{j})
-\operatorname{sinh}(1/2+q_{i}-q_{j})
]
\bigg\}\no\\
&+k_{3}\bigg\{\sum_{j>m}p_{m}p_{j}[\operatorname{sinh}(-1\!+\!q_{j}\!-\!q_{m})\!-\!\operatorname{sinh}(q_j\!-\!q_{m})]
\!+\!\sum_{j<m}p_{m}p_{j}[\operatorname{sinh}(q_j\!-\!q_{m})\!-\!\operatorname{sinh}(1\!+\!q_{j}\!-\!q_{m})]\bigg\}.
\end{align}

According to the  arbitrariness
of the function $\psi(t,x)$,
one derives the dynamical system (\ref{multi-pp}) satisfied by  $p_{m}(t)$ and $q_{m}(t)\, (m=1,2,...,N)$
after setting the coefficients given by Eqs.~(\ref{coe1}) and (\ref{coe2}) both to be zero.
 This completes the proof of Proposition~\ref{mppeakon}.

\vspace{0.15in}
\noindent \textbf{Acknowledgement} \vspace{0.05in}

This work was partially supported by the NSF of China under Grants Nos. 11925108 and 11731014.

\vspace{0.15in}
\noindent {\bf Appendix A.\, Basics properties of the  Littlewood-Paley theory} \\

\setcounter{equation}{0}
\renewcommand\theequation{A.\arabic{equation}}

Let $B(x_{0},r)$ be the open ball centered at $x_{0}$ with
radius $r,$ $\mathcal{C}\equiv \{\xi\in \mathbb{R}^{d} | 4/3\leq|\xi|\leq 8/3\}$,
 and $\mathcal{\tilde{C}}\equiv B(0,2/3)+\mathcal{C}.$ Then there are two radial functions $\chi\in \mathcal{D}(B(0,4/3))$
 and $\varphi\in \mathcal{D}(\mathcal{C})$
 satisfying
\begin{equation*}
\left\{\begin{array}{l}
\chi(\xi)+\sum_{q \geq 0} \varphi(2^{-q} \xi)=1,\quad
1/3 \leq \chi^2(\xi)+\sum_{q \geq 0} \varphi^2(2^{-q} \xi) \leq 1 \quad (\forall \xi \in \mathbb{R}^{d}),\v\\
|q-q^{\prime}| \geq 2 \Rightarrow \operatorname{Supp} \varphi(2^{-q} \cdot) \cap \operatorname{Supp} \varphi(2^{-q^{\prime}} \cdot)=\varnothing, \v \\
q \geq 1 \Rightarrow \operatorname{Supp} \chi(\cdot) \cap \operatorname{Supp} \varphi(2^{-q^{\prime}} \cdot)=\varnothing,\quad
|q-q^{\prime}| \geq 5 \Rightarrow 2^{q^{\prime}} \widetilde{\mathcal{C}} \cap 2^{q} \mathcal{C}=\varnothing.
\end{array}\right.
\end{equation*}

The dyadic operators $\Delta_{q}$ and $S_{q}$ acting on $u(t,x)\in S'(\mathbb{R}^d)$ are defined as
\begin{equation*}
\begin{array}{l}
\Delta_{q} u=\left\{
 \begin{array}{ll}
   0, & q \leq-2, \v \\
 \chi(D) u=\int_{\mathbb{R}^{d}} \tilde{h}(y) u(x-y)dy, & q=-1, \v \\
 \varphi\left(2^{-q} D\right) u=2^{q d} \int_{\mathbb{R}^{d}} h\left(2^{q} y\right) u(x-y) dy, & q \geq 0,
 \end{array}\right. \v \\
S_{q} u=\sum_{q^{\prime} \leq q-1} \Delta_{q^{\prime}} u,
\end{array}
\end{equation*}
where
$h = \mathcal{F}^{-1} \varphi$ and $ \tilde{h} = \mathcal{F}^{-1} \chi$ with $\mathcal{F}^{-1}$ denoting the inverse Fourier transform.

The Besov spaces is
$B_{p, r}^{s}(\mathbb{R}^{d})=\left\{u \in S^{\prime}\, \big|\, \|u\|_{B_{p, r}^{s}(\mathbb{R}^{d})}=\big(\sum_{j \geq-1} 2^{r j s}\|\Delta_{j} u\|_{L^{p}(\mathbb{R}^{d})}^{r}\big)^{1/r}
<\infty\right\}.$
With the above-defined Besov spaces, we next recall some of their properties.

\begin{lem} {\rm (Embedding property)}\cite{BahouriCheminDanchin2011,Chemin2004CRM,Danchin2003}\label{21lem1}
Suppose $1 \leq p_{1} \leq p_{2} \leq \infty$, $1 \leq r_{1} \leq r_{2} \leq \infty$ and $s$ be real.
Then it holds that $B_{p_{1}, r_{1}}^{s}(\mathbb{R}^{d}) \hookrightarrow B_{p_{2}, r_{2}}^{s-d(1/p_1-1/p_2)}
(\mathbb{R}^{d})$. If $s>d/p$ or $s=d/p,\, r=1,$ then there holds $B_{p, r}^{s}(\mathbb{R}^{d}) \hookrightarrow L^{\infty}(\mathbb{R}^{d})$.
\end{lem}

\begin{lem}{\rm (Interpolation)}\cite{BahouriCheminDanchin2011,Chemin2004CRM,Danchin2003} \label{21lem2}
Let $s_1,\, s_2$  be real numbers with $s_{1}<s_{2}$ and $\theta \in(0,1).$
 Then  there exists a constant $C$ such that
\begin{eqnarray*}
\|u\|_{B_{p, r}^{\theta s_{1}+(1-\theta) s_{2}}}
\leq\|u\|_{B_{p, r}^{s_{1}}}^{\theta}\|u\|_{B_{p, r}^{s_{2}}}^{(1-\theta)},\quad
\|u\|_{B_{p, 1}^{\theta s_{1}+(1-\theta) s_{2}}}
\leq \frac{C}{s_{2}-s_{1}}\frac{1}{\theta (1-\theta)}
\|u\|_{B_{p, \infty}^{s_{1}}}^{\theta}
\|u\|_{B_{p, \infty}^{s_{2}}}^{(1-\theta)},
\end{eqnarray*}
where $(p, r) \in[1, \infty]^{2}$.
\end{lem}

\begin{lem}{\rm (Product law)}\cite{BahouriCheminDanchin2011,Chemin2004CRM,Danchin2003}\label{21lem3}
Let  $(p, r)\in [1, \infty]^{2}$
and $s$ be real. Then
$
\|u v\|_{B_{p, r}^{s}(\mathbb{R}^{d})} \leq C(\|u\|_{L^{\infty}(\mathbb{R}^{d})}\|v\|_{B_{p, r}^{s}(\mathbb{R}^{d})}
+\|u\|_{B_{p, r}^{s}(\mathbb{R}^{d})}
\|v\|_{L^{\infty}(\mathbb{R}^{d})}),
$
namely,  the space $L^{\infty}(\mathbb{R}^{d}) \cap B_{p, r}^{s}(\mathbb{R}^{d})$ is an algebra.
Moreover, if  $s>d/p$ or $s=d/p,\, r=1,$ then there holds
$
\|u v\|_{B_{p, r}^{s}\left(\mathbb{R}^{d}\right)} \leq C\|u\|_{B_{p, r}^{s}\left(\mathbb{R}^{d}\right)}\|v\|_{B_{p, r}^{s}\left(\mathbb{R}^{d}\right)}.
$
\end{lem}

\begin{lem}{\rm (Moser-type estimates)}\cite{BahouriCheminDanchin2011,Danchin2001DIE}         \label{21lem4}
Let $s>\max \{d/p,\, d/2\}$
and $(p, r)\in[1, \infty]^{2}$.
Then, for any $a \in B_{p, r}^{s-1}(\mathbb{R}^{d})$ and $b \in B_{p, r}^{s}(\mathbb{R}^{d}),$ there holds
$
\|a b\|_{B_{p, r}^{s-1}(\mathbb{R}^{d})} \leq C\|a\|_{B_{p, r}^{s-1}(\mathbb{R}^{d})}\|b\|_{B_{p, r}^{s}(\mathbb{R}^{d})}.
$
\end{lem}

The following  Lemma is useful  for proving the blow-up criterion.
\begin{lem}{\rm (Moser-type estimates)}\cite{GuiLiu2010JFA,
GuiLiuOlverQu2013CMP}  \label{21lem5} Let $s \geq 0.$ Then one has
$$
\begin{aligned}
\|f g\|_{H^{s}(\mathbb{R})}
& \leq C(\|f\|_{H^{s}(\mathbb{R})}\|g\|_{L^{\infty}(\mathbb{R})}
+\|f\|_{L^{\infty}(\mathbb{R})}\|g\|_{H^{s}(\mathbb{R})}),\v \\
\|f \partial_{x} g\|_{H^{s}(\mathbb{R})}
 & \leq C(\|f\|_{H^{s+1}(\mathbb{R})}\|g\|_{L^{\infty}(\mathbb{R})}
 +\|f\|_{L^{\infty}(\mathbb{R})}\|\partial_{x}g\|_{H^{s}(\mathbb{R})}),
\end{aligned}
$$
where $C$'s are constants independent of $f$ and $g$.
\end{lem}

\begin{lem}\cite{BahouriCheminDanchin2011}\label{besovProperty}
 Let $s \in \mathbb{R}$ and $1 \leq p, r \leq \infty$.
Then the Besov spaces have the following properties:
\begin{itemize}
\item{}$B_{p, r}^{s}(\mathbb{R}^{d})$ is a Banach space and continuously embedding into $\mathcal{S}^{\prime}(\mathbb{R}^{d}),$ where $\mathcal{S}^{\prime}(\mathbb{R}^{d})$ is the dual space of the Schwartz space $\mathcal{S}(\mathbb{R}^{d})$;

\item{}If $p, r<\infty,$ then $\mathcal{S}(\mathbb{R}^{d})$ is dense in $B_{p, r}^{s}(\mathbb{R}^{d})$;

\item {} If $u_{n}$ is a bounded sequence of $B_{p, r}^{s}(\mathbb{R}^{d}),$ then an element $u \in B_{p, r}^{s}(\mathbb{R}^{d})$ and a subsequence $u_{n_{k}}$ exist such that
$\lim _{k \rightarrow \infty} u_{n_{k}}=u \text { in } \mathcal{S}^{\prime}(\mathbb{R}^{d}) \text { and }\|u\|_{B_{p, r}^{s}(\mathbb{R}^{d})} \leq C \liminf _{k \rightarrow \infty}\|u_{n_{k}}\|_{B_{p, r}^{s}(\mathbb{R}^{d})}.$
\end{itemize}
\end{lem}

\vspace{0.15in}
\noindent {\bf Appendix B.\,  Some lemmas in the theory of the transport equation} \\

\renewcommand\theequation{B.\arabic{equation}}

We recall some a priori estimates\cite{BahouriCheminDanchin2011,
Danchin2001DIE} for the following transport equation
\begin{equation}\label{transport}
\phi_{t}+v\cdot \nabla \phi=\omega, \quad \left.\phi\right|_{t=0}=\phi_{0}.
\end{equation}

\begin{lem}\cite{BahouriCheminDanchin2011,
Danchin2001DIE}\label{22lem1} Let $1 \leq p \leq p_{1} \leq \infty,$ $1 \leq r \leq \infty$
and $s \geq -d \min (1/p_1,\, 1-1/p)$.
Let  $\phi_{0} \in B_{p, r}^{s}(\mathbb{R}^{d})$.
$\omega \in L^{1}([0, T] ; B_{p, r}^{s}(\mathbb{R}^{d}))$ and $\nabla v \in L^{1}([0, T] ; B_{p, r}^{s}(\mathbb{R}^{d}) \cap L^{\infty}(\mathbb{R}^{d})),$ then there exists a unique solution
 $\phi \in L^{\infty}([0, T] ; B_{p, r}^{s}(\mathbb{R}^{d}))$  to Eq.~(\ref{transport})  satisfying:
\begin{eqnarray}
\|\phi\|_{B_{p, r}^{s}(\mathbb{R}^{d})} \leq\|\phi_{0}\|_{B_{p, r}^{s}(\mathbb{R}^{d})}
+\int_{0}^{t}[\|\omega(t^{\prime})\|_{B_{p, r}^{s}(\mathbb{R}^{d})}+C U_{p_{1}}(t^{\prime})\|\phi(t^{\prime})\|_{B_{p, r}^{s}
(\mathbb{R}^{d})}] d t^{\prime}, \label{priori1}\\
\|\phi\|_{B_{p, r}^{s}(\mathbb{R}^{d})}
\leq\left[\|\phi_{0}\|_{B_{p, r}^{s}(\mathbb{R}^{d})}+\int_{0}^{t}\|\omega(t^{\prime})\|_{B_{p, r}^{s}
(\mathbb{R}^{d})}e^{-C U_{p_{1}}(t^{\prime})} d t^{\prime}\right]e^{C U_{p_1}(t)},\label{priori2}
\end{eqnarray}
where $U_{p_{1}}(t)=\int_{0}^{t}\|\nabla v\|_{B_{p_{1}, \infty}^{d/p_{1}}(\mathbb{R}^{d}) \cap L^{\infty}(\mathbb{R}^{d})} d t^{\prime}$ if $s<1+d/p_{1}, U_{p_{1}}(t)=\int_{0}^{t}\|\nabla v\|_{B_{p_{1}, r}^{s-1}\left(\mathbb{R}^{d}\right)} d t^{\prime}$
if $s>1+d/p_{1}$ or $s=1+d/p_{1},\, r=1,$ and $C$ is a constant depending only on $s, p, p_{1}$, and $r$.
\end{lem}

\begin{lem}\cite{BahouriCheminDanchin2011}\label{22lem4} Let $s \geq -d \min (1/p_1, 1-1/p).$
Let $\phi_{0} \in B_{p, r}^{s}(\mathbb{R}^{d})$, $\omega \in L^{1}([0, T] ; B_{p, r}^{s}(\mathbb{R}^{d}))$  and $v \in L^{\rho}([0, T]; B_{\infty, \infty}^{-M}(\mathbb{R}^{d}))$  for some $\rho>1$ and $M>0$  be a time-dependent vector field satisfying
$$
\nabla v\in\left\{\begin{array}{ll}
L^{1}([0, T] ; B_{p_{1}, \infty}^{d/p}(\mathbb{R}^{d})), &  {\rm if }\,\, s<1+d/p_{1}, \v \\
L^{1}\left([0, T] ; B_{p_{1}, \infty}^{s-1}\left(\mathbb{R}^{d}\right)\right), & {\rm if }\,\, s>1+d/p_1\,\,  {\rm or } \,\, s=1+d/p_1 \,\, {\rm and } \,\, r=1.
\end{array}\right.
$$
Then, Eq.~(\ref{transport}) has a unique solution $\phi\in \mathcal{C}([0, T] ; B_{p, r}^{s}(\mathbb{R}^{d}))$ for $r<\infty$, or
$\phi\in (\bigcap_{s^{\prime}<s} \mathcal{C}([0, T]; B_{p, \infty}^{s^{\prime}}(\mathbb{R}^{d}))) \cap \mathcal{C}_{w}([0, T] ; B_{p, \infty}^{s}(\mathbb{R}^{d})))$ for $r=\infty$. Furthermore, the inequalities
(\ref{priori1})-(\ref{priori2}) hold.
\end{lem}

\begin{lem} {\rm (A priori estimate in the Sobolev spaces)}\cite{BahouriCheminDanchin2011,
GuiLiu2010JFA}\label{22lem5} Let $0\leq\sigma<1$. Let $\phi_{0} \in H^{\sigma},\, \omega \in L^{1}(0, T ; H^{\sigma})$ and $\partial_{x} v \in L^{1}(0, T ; L^{\infty}).$ Then the solution $\phi$
to Eq.~(\ref{transport}) belongs to $C([0, T] ; H^{\sigma}).$
  More precisely,
 there is a constant $C$
depending only on $\sigma$ such that
$$
\|\phi\|_{H^{\sigma}}\leq\|\phi_{0}\|_{H^{\sigma}}
+\int_{0}^{t}[\|\omega(\tau)\|_{H^{\sigma}}+CU^{\prime}(\tau)\|\phi(\tau)\|_{H^{\sigma}}] d\tau,\quad U(t)=\int_{0}^{t}\|\partial_{x} v(\tau)\|_{L^{\infty}} \mbox{d} \tau.
$$
\end{lem}

\end{document}